\documentclass[12pt]{amsart}
\usepackage{blindtext}
\usepackage{url}
\usepackage{tabularx}
\usepackage[utf8]{inputenc}
\usepackage[margin=2.5cm]{geometry}
\usepackage[utf8]{inputenc}
\usepackage{amsmath}
\usepackage{amsfonts}
\usepackage{color}
\usepackage{slashbox}
\usepackage{amssymb}
\usepackage{amsthm}
\usepackage{array}
\usepackage{rotating}
\usepackage{physics}
\usepackage{tikz}
\usepackage{braids}
\usepackage{tikz-cd}
\tikzset{%
	symbol/.style={%
		draw=none,
		every to/.append style={%
			edge node={node [sloped, allow upside down, auto=false]{$#1$}}}
	}
}
\usepackage[all]{xy}
\usepackage{cite}
\usepackage{mathdots}
\usepackage{pstricks}
\usepackage{mathrsfs}
\usepackage{setspace}
\usepackage{hyperref}
\theoremstyle{theorem}
\newtheorem{theorem}{Theorem}
\newtheorem{diagram}[theorem]{Diagram}
\numberwithin{theorem}{subsection}
\newtheorem{lemma}[theorem]{Lemma}
\newtheorem{proposition}[theorem]{Proposition}
\newtheorem{corollary}[theorem]{Corollary}
\theoremstyle{definition}
\newtheorem{definition}[theorem]{Definition}
\newtheorem{example}[theorem]{Example}
\theoremstyle{remark}
\newtheorem{remark}[theorem]{Remark}
\theoremstyle{notation}
\newtheorem{notation}[theorem]{Notation}
\title{A semi-strictly generated closed structure on $\mathbf{Gray}$-$\mathbf{Cat}$}
\author{Adrian Miranda}
\thanks{The material is based on work supported by Macquarie University under MQRES PhD Scholarship 20192497 and by EPSRC under grant EP/V002325/2. Except for Section 5, this material is based on my PhD thesis \cite{Miranda PhD}. I am grateful to my supervisor, Steve Lack, for his guidance. I would also like to thank Nicola Gambino for advice while I was preparing this paper, and John Bourke and Gabriele Lobbia for useful conversations about the relationship between the closed structure of this paper and the closed, skew monoidal structure they defined in \cite{Bourke Lobbia Skew Approach to Enrichment for Gray-categories}.}
\address{Department of Mathematics, University of Manchester, 
	United Kingdom}
\email{adrian.miranda@manchester.ac.uk}

\begin{document}

\maketitle

\onehalfspacing

\begin{abstract}
	\noindent We show that the semi-strictly generated internal homs of $\mathbf{Gray}$-categories $[\mathfrak{A}, \mathfrak{B}]_\text{ssg}$ defined in \cite{Miranda strictifying operational coherences} underlie a closed structure on the category $\mathbf{Gray}$-$\mathbf{Cat}$ of $\mathbf{Gray}$-categories and $\mathbf{Gray}$-functors. The morphisms of $[\mathfrak{A}, \mathfrak{B}]_\text{ssg}$ are composites of those trinatural transformations which satisfy the unit and composition conditions for pseudonatural transformations on the nose rather than up to an invertible $3$-cell. Such trinatural transformations leverage three-dimensional strictification \cite{Miranda strictifying operational coherences} while overcoming the challenges posed by failure of middle four interchange to hold in $\mathbf{Gray}$-categories \cite{Bourke Gurski Cocategorical Obstructions to a Tensor Product of Gray Categories}. As a result we obtain a closed structure that is only partially monoidal with respect to \cite{crans tensor of gray categories}. As a corollary we obtain a slight strengthening of strictification results for braided monoidal bicategories \cite{Gurski Loop Spaces}, which will be improved further in a forthcoming paper \cite{Miranda weak interchange 4-categories}.
\end{abstract}

\tableofcontents

\section{Introduction}

\subsection{Motivation and background}

\noindent One way to describe $\left(n+1\right)$-dimensional categories is via enrichment over $n$-dimensional categories. Indeed, strict $n$-categories are defined via iterated enrichment over cartesian closed categories, starting with $\mathbf{Set}$. However, some amount of weakness is needed in order to describe categories occurring in mathematical practice, such as truncations of fundamental $\infty$-groupoids of topological spaces, even up to suitable equivalence \cite{GPS tricategory}. In dimension three, one such semi-strict notion keeps the middle four interchange law weak, but is otherwise strict. Such three-dimensional categories are called $\mathbf{Gray}$-categories, and they may be described via enrichment over a symmetric monoidal closed structure on $2$-$\mathbf{Cat}$. The closed fragment of this structure consists of $2$-functors, pseudonatural transformations and modifications. These higher dimensional maps between $2$-categories are characterised by preserving operations, such as composition and identities, on the nose, but only satisfying naturality like conditions up to coherent higher dimensional constraint data.
\\
\\
\noindent In a previous paper \cite{Miranda strictifying operational coherences} we identified those maps between $\mathbf{Gray}$-categories which preserve operations strictly but naturality like conditions weakly. For functors, this recovers the $\mathbf{Gray}$-enriched notion while for transformations it gives a semi-strict notion of trinatural transformation. Unfortunately, semi-strict trinatural transformations fail to be closed under composition due to failure of middle four interchange to hold strictly in a $\mathbf{Gray}$-category. Despite a symmetric monoidal product of $\mathbf{Gray}$-categories having been defined in \cite{crans tensor of gray categories}, it fails to have a corresponding closed structure for related reasons, as made precise in \cite{Bourke Gurski Cocategorical Obstructions to a Tensor Product of Gray Categories}. Recently \cite{Bourke Lobbia Skew Approach to Enrichment for Gray-categories}, closed skew monoidal structures on the category $\mathbf{Gray}$-$\mathbf{Cat}$ of $\mathbf{Gray}$-categories and $\mathbf{Gray}$-functors are being pursued as an alternative setting over which to enrich for a notion of four-dimensional category. The morphisms appearing in the internal homs of such skew structures are transformations which respect the identity structure on $1$-cell of their domain $\mathbf{Gray}$-categories strictly, but only respect their $1$-cell composition structure weakly. 

\subsection{Our contributions}

\noindent In Theorem \ref{Gray cat closed structure final theorem} we prove that the semi-strictly generated homs $[\mathfrak{A}, \mathfrak{B}]_\text{ssg}$, as defined in \cite{Miranda strictifying operational coherences} and recalled in Definition \ref{semi-strictly generated hom of Gray categories}, equip the category $\mathbf{Gray}$-$\mathbf{Cat}$ with a closed structure in the sense of Eilenberg and Kelly \cite{Eilenberg Kelly Closed Categories}, \cite{Manzyuk}, \cite{Laplaza}. This demonstrates that failure of semi-strict transformations to be closed under composition is the \emph{only} obstruction preventing them from featuring in a closed structure on $\mathbf{Gray}$-$\mathbf{Cat}$. As a corollary of our results we obtain in Theorem \ref{coherence for braided Gray monoids} a slight strengthening of strictification results for braided monoidal bicategories \cite{Baez Neuchl Braided Monoidal 2-categories}, \cite{Crans Generalised Centers}, \cite{Gurski Loop Spaces}. This will be further strengthened in the forthcoming paper \cite{Miranda weak interchange 4-categories} where a better behaved base for four-dimensional categories will be described. Finally, we relate our work to the tensor product of \cite{crans tensor of gray categories} and to the closed structures of \cite{Bourke Lobbia Skew Approach to Enrichment for Gray-categories}. Categories enriched over the closed structure of this paper will be described in detail in \cite{Miranda weak interchange 4-categories}. As we discuss in Remark \ref{undesirable aspects of enrichment}, they are unlikely to be models of an essentially algebraic theory \cite{Adamek Rosicky Locally Presentable and Accessible Categories}. As models of four-dimensional categories, the categories enriched over the closed structure of this paper are weaker than the $4$-teisi of \cite{Crans Braidings Syllapses and Symmetries}, but stricter than those enriched over the skew contexts considered in \cite{Bourke Lobbia Skew Approach to Enrichment for Gray-categories}.
\\
\\
\noindent Our proofs in Sections \ref{Data of the closed structure} and \ref{Axioms for closed structure} do not assume the results of \cite{Bourke Lobbia Skew Approach to Enrichment for Gray-categories}. Were we to assume their results, we would only need to check that their internal whiskering operators respect semi-strict decomposability and hence restrict to $\mathbf{Gray}$-functors between the semi-strictly generated homs. Our work is independent of theirs, and our techniques are more similar to those in the Appendices of \cite{Buhne PhD}. We give explicit proofs that various interchangers between trinatural transformations and trimodifications are well-defined, using pasting diagrams in hom-$2$-categories rather than the calculus of skew multicategories as in \cite{Bourke Lobbia Skew Approach to Enrichment for Gray-categories} or string diagrams as in \cite{Buhne PhD}. In Remark \ref{proof strategy pasting diagram chases} we outline a strategy via which these explicit proofs can be constructed. These calculations serve to re-prove almost all aspects of the sharp closed structure of \cite{Bourke Lobbia Skew Approach to Enrichment for Gray-categories}; we only need the semi-strictly generated property for the proof of the fragment of the associativity axiom considered in Proposition \ref{assoc condition underlying categories}. It is then straightforward to restrict these results to ones about semi-strictly generated trinatural transformations and establish the semi-strictly generated closed structure of this paper. We will elaborate further on the relationship between our techniques and those of \cite{Bourke Lobbia Skew Approach to Enrichment for Gray-categories} in Remark \ref{comparing techniques to Bourke Lobbia}.

\subsection{Structure of the paper}

\begin{itemize}
	\item Subsection \ref{Subsection recall closed structure} recalls the notion of a closed structure on a category, and of closed functors.
	\item Remark \ref{higher cells between Gray categories as functions on globular input data} describes $\mathbf{Gray}$-functors, trinatural transformations, trimodifications and perturbations as functions on globular input data. This perspective will clarify proofs in later sections of this paper.
	\item Definition \ref{semi-strictly generated hom of Gray categories} recalls the semi-strictly generated hom of $\mathbf{Gray}$-categories from \cite{Miranda strictifying operational coherences}. This will feature as the internal hom $[-, ?]: \mathbf{Gray}\text{-}\mathbf{Cat}^\text{op} \times \mathbf{Gray}\text{-}\mathbf{Cat} \rightarrow \mathbf{Gray}\text{-}\mathbf{Cat}$ of the closed structure on $\mathbf{Gray}$-$\mathbf{Cat}$. Its functoriality is established in Proposition \ref{coherent hom on Gray Cat}.
	\item Proposition \ref{identity assigners in closed structure} establishes the identity assigners $i_\mathfrak{A}:\mathbf{1} \rightarrow [\mathfrak{A}, \mathfrak{A}]$, and Proposition \ref{constant assigners closed structure} establishes the constant map assigners $c_\mathfrak{A}: \mathfrak{A} \rightarrow [\mathbf{1}, \mathfrak{A}]$.
	\item The most complicated fragment of the closed structure are the internal whiskering operators $[\mathfrak{A}, -]: [\mathfrak{B}, \mathfrak{C}] \rightarrow [[\mathfrak{A}, \mathfrak{B}], [\mathfrak{A}, \mathfrak{C}]]$. Subsection \ref{Subsection underlying functoriality of internal whiskering operators} establishes their functoriality at the level of underlying categories, as stated in Corollary \ref{underlying functor of internal whiskering operator in closed structure}. Subsection \ref{Subsection Gray functoriality of internal whiskering operators} treats $\mathbf{Gray}$-functoriality, establishing their $2$-functoriality between hom-$2$-categories in Proposition \ref{local 2 functoriality of whiskering}, and enrichment over $\mathbf{Gray}$ in Theorem \ref{internal whiskering operators closed structure on Gray cat}. By this stage, all of the data involved in the closed structure will have been established, and it will remain to check the axioms.
	\item Section \ref{Axioms for closed structure} treats the axioms for the closed structure. All axioms except for associativity are addressed in Proposition \ref{Non-associativity axioms for closed structure on Gray cat}. The associativity axiom is shown one dimension at a time, over the course of Proposition \ref{associativity condition underlying sets}, Proposition \ref{assoc condition underlying categories}, Proposition \ref{associativity condition for closed structure underlying sesquicategories}, and finally Theorem \ref{Gray cat closed structure final theorem}, in which we state the closed structure in full.
	\item Section \ref{subsection interaction with the skew structure} relates the semi-strictly generated closed structure to the closed structure $\mathbf{Ps}\left(-, ?\right)_\text{s}$ on $\mathbf{Gray}$-$\mathbf{Cat}$ introduced in \cite{Bourke Lobbia Skew Approach to Enrichment for Gray-categories}.
	\item Section \ref{Section Centres of braided Gray monoids} relates the various interchangers between trinatural transformations and trimodifications developed in Section \ref{Data of the closed structure} to the generalised centre construction for $\mathbf{Gray}$-monoids of \cite{Crans Generalised Centers}. As a corollary, we leverage the results of \cite{Miranda strictifying operational coherences} to improve the strictification for braided monoidal bicategories described in \cite{Gurski Coherence in Three Dimensional Category Theory}.
	\item This paper relies heavily on explicit descriptions of the various notions of weak higher cells between $\mathbf{Gray}$-functors, and the structure of the $\mathbf{Gray}$-categories that these data form. These have been described in Section 3 of \cite{Miranda strictifying operational coherences}. So that this paper can be self-contained, we also recall these in the Appendices \ref{Appendices}.
\end{itemize}

\begin{notation}
	Many of the proofs in this paper involve pasting diagram chases in hom-$2$-categories of $\mathbf{Gray}$-categories. We will use colour to draw the reader's attention to new data appearing in each step of such proofs. For each step, there will also be some explanation of the axioms for trinatural transformations, trimodifications or perturbations that are being used. These explanations will refer to properties as they are described in Appendix \ref{subsection weak higher maps between Gray functors}. On the other hand, $\mathbf{Gray}$-category axioms and $\mathbf{Gray}$-functoriality will be used without comment. Given a pair of $2$-cells in a $\mathbf{Gray}$-category $\mathfrak{A}$ as displayed below left, their interchanger is an invertible $2$-cell in the hom-$2$-category $\mathfrak{A}(X, Z)$, and will be denoted as displayed below right.

$$\begin{tikzcd}
	X \arrow[rr, bend left = 30, "f"name=A]\arrow[rr, bend right = 30, "f'"'name=B] && Y \arrow[rr, bend left = 30, "g"name=C]\arrow[rr, bend right = 30, "g'"'name=D] &&Z\arrow[from=A, to=B, Rightarrow, shorten = 10, "\phi"]\arrow[from=C, to=D, Rightarrow, shorten = 10, "\psi"]
	&{}
\end{tikzcd}\begin{tikzcd}
gf\arrow[rr, "g.\phi"] \arrow[dd, "\psi.f"']&{}\arrow[dd, Rightarrow, shorten = 10, "\psi_{\phi}"]&gf'\arrow[dd, "\psi.f'"]
\\
\\
g'f\arrow[rr, "g'.\phi"']&{}&g'f'
\end{tikzcd}$$
\end{notation}

\section{Preliminaries}

\subsection{Closed structures}\label{Subsection recall closed structure}

\begin{definition}\label{definition of a closed structure on a category}
	A \emph{closed structure} on a category $\mathcal{E}$ consists of the following data subject to the following axioms.
	\\
	\\
	\noindent DATA
	\\
	\begin{itemize}
		\item A functor $[-, ?]: \mathcal{E}^\text{op} \times \mathcal{E} \rightarrow \mathcal{E}$ called the \emph{internal hom}.
		\item An object $\mathbf{1} \in \mathcal{E}$ called the \emph{unit}.
		\item A family of maps $i_\mathfrak{A}:\mathbf{1} \rightarrow [\mathfrak{A}, \mathfrak{A}]$ called the \emph{identity assigners}, varying extra-naturally in $\mathfrak{A} \in \mathcal{E}$.
		\item A family of isomorphisms $c_\mathfrak{A} : \mathfrak{A}\rightarrow [\mathbf{1}, \mathfrak{A}]$ called the \emph{constant assigners}, varying naturally in $\mathfrak{A}$.
		\item A family of maps $[\mathfrak{A}, -]: [\mathfrak{B}, \mathfrak{C}] \rightarrow \big{[}[\mathfrak{A}, \mathfrak{B}], [\mathfrak{A}, \mathfrak{C}]\big{]}$ called the \emph{internal whiskering operators}, varying naturally in $\mathfrak{B}$ and $\mathfrak{C}$ and extra-naturally in $\mathfrak{A}$.
	\end{itemize}
	\noindent AXIOMS
	\\
	\begin{enumerate}
		\item The function assigning a morphism $F \in \mathcal{E}\left(\mathfrak{A}, \mathfrak{B}\right)$ to the map \begin{tikzcd}[column sep = 15]
			\mathbf{1} \arrow[rr, "i_{\mathfrak{A}}"] && {[}\mathfrak{A}{,}\mathfrak{A}{]} \arrow[rr, "{[}\text{id}{,}F{]}"] &&  {[}\mathfrak{A}{,}\mathfrak{B}{]}
		\end{tikzcd} in $\mathcal{E}\left(\mathbf{1}, [\mathfrak{A}, \mathfrak{B}]\right)$ is a bijection.
		\item \emph{(Left unit)} For any $\mathfrak{A}$, $\mathfrak{B} \in \mathcal{E}$, the diagram below left commutes.
		\item \emph{(Right unit)} For any $\mathfrak{A}$, $\mathfrak{B} \in \mathcal{E}$, the diagram below centre commutes.
		\item \emph{(Constant maps)} For any $\mathfrak{A}$, $\mathfrak{B} \in \mathcal{E}$, the diagram below right commutes. 
		
		$$\begin{tikzcd}[column sep = 12, row sep = 15, font=\fontsize{9}{6}]
			\mathbf{1} \arrow[rrr, "i_{\mathfrak{B}}"]
			\arrow[rrrddd, "i_{{[}\mathfrak{A}{,}\mathfrak{B}{]}}"']
			&&& {[}\mathfrak{B}{,}\mathfrak{B}{]}
			\arrow[ddd, "{[}\mathfrak{A}{,}-{]}"]
			\\
			\\
			\\
			&&& {[}{[}\mathfrak{A}{,}\mathfrak{B}{]}{,}{[}\mathfrak{A}{,}\mathfrak{B}{]}{]}&{}
		\end{tikzcd}\begin{tikzcd}[column sep = 12, font=\fontsize{9}{6}, row sep = 15]
			{[}\mathfrak{A}{,}\mathfrak{B}{]} \arrow[rrr, "{[}\mathfrak{A}{,}-{]}"]
			\arrow[rrrddd, "c_{{[}\mathfrak{A}{,}\mathfrak{B}{]}}"']
			&&& {[}{[}\mathfrak{A}{,}\mathfrak{A}{]}{,}{[}\mathfrak{A}{,}\mathfrak{B}{]}{]}
			\arrow[ddd, "{[}i_{\mathfrak{A}}{,}\text{id}{]}"]
			\\
			\\
			\\
			&&& {[}\mathbf{1}{,}{[}\mathfrak{A}{,}\mathfrak{B}{]}{]}&{}
		\end{tikzcd}\begin{tikzcd}[column sep = 12, font=\fontsize{9}{6}, row sep = 15]
			{[}\mathfrak{A}{,}\mathfrak{B}{]} \arrow[rrr, "{[}\mathbf{1}{,}-{]}"]
			\arrow[rrrddd, "{{[}\text{id}{,}c_{\mathfrak{B}}{]}}"']
			&&& {[}{[}\mathbf{1}{,}\mathfrak{A}{]}{,}{[}\mathbf{1}{,}\mathfrak{B}{]}{]}
			\arrow[ddd, "{[}c_{\mathfrak{A}}{,}\text{id}{]}"]
			\\
			\\
			\\
			&&& {[}\mathfrak{A}{,}{[}\mathbf{1}{,}\mathfrak{B}{]}{]}
		\end{tikzcd}$$
		
		\item \emph{(Associativity)} For any $\mathfrak{A}$, $\mathfrak{B}$, $\mathfrak{C}$, $\mathfrak{D} \in \mathcal{E}$, the following diagram commutes.
	\end{enumerate}

	$$\begin{tikzcd}[column sep = 1, font=\fontsize{9}{6}, row sep = 15]
		&&
		\big{[}{{[}\mathfrak{B}{,}\mathfrak{C}{]}}{,}{{[}\mathfrak{B}{,}\mathfrak{D}{]}}\big{]}
		\arrow[rrdd, "\big{[}\text{id}{,}{[}\mathfrak{A}{,}-{]}\big{]}"]
		\\
		\\
		{[}\mathfrak{C}{,}\mathfrak{D}{]}
		\arrow[rruu, "{[}\mathfrak{B}{,}-{]}"]
		\arrow[ddr, "{[}\mathfrak{A}{,}-{]}"']
		&&&&
		\bigg{[}{[}\mathfrak{B}{,}\mathfrak{C}{]}{,}\big{[}{[}\mathfrak{A}{,}\mathfrak{B}{]}{,}{[}\mathfrak{A}{,}\mathfrak{D}{]}\big{]}\bigg{]}
		\\
		\\
		&
		\big{[}{{[}\mathfrak{A}{,}\mathfrak{C}{]}}{,}{{[}\mathfrak{A}{,}\mathfrak{D}{]}}\big{]}
		\arrow[rr, "\big{[}{[}\mathfrak{A}{,}\mathfrak{B}{]}{,}-\big{]}"']
		&&
		\bigg{[}\big{[}{[}\mathfrak{A}{,}\mathfrak{B}{]}{,}{[}\mathfrak{A}{,}\mathfrak{C}{]}\big{]}{,}\big{[}{[}\mathfrak{A}{,}\mathfrak{B}{]}{,}{[}\mathfrak{A}{,}\mathfrak{D}{]}\big{]}\bigg{]}
		\arrow[ruu,"\big{[}{[}\mathfrak{A}{,}{-}{]}{,}{\text{id}}\big{]}"']
	\end{tikzcd}$$
\end{definition}

\begin{definition}\label{Definition normal skew closed functor}
	Let $\left(\mathcal{E}, [-, ?], I, i, c\right)$ and $\left(\mathcal{D}, \{-, ?\}, J, j, k\right)$ be closed categories as in Definition \ref{definition of a closed structure on a category}. A \emph{normal closed functor} from $\left(\mathcal{E}, [-, ?], I, i, c\right)$ to $\left(\mathcal{D}, \{-, ?\}, J, j, k\right)$ consists of the following data subject to the following axioms.
	\\
	\\
	\noindent DATA
	\\
	\begin{itemize}
		\item A functor $F: \mathcal{E} \rightarrow \mathcal{D}$,
		\item A natural transformation as depicted below.
	\end{itemize}

$$\begin{tikzcd}
	\mathcal{E}^\text{op}\times \mathcal{E}
	\arrow[rr, "{[}-{,}?{]}"]
	\arrow[dd, "F^\text{op}\times F"'] &{}\arrow[dd, shorten = 10, Rightarrow, "\chi"]& \mathcal{E}\arrow[dd, "F"]
	\\
	\\
	\mathcal{D}^\text{op}\times \mathcal{D}
	\arrow[rr, "\{-{,}?\}"'] &{}& \mathcal{D}
\end{tikzcd}$$

\noindent AXIOMS
\\
\begin{enumerate}
	\item \emph{(Normality)} $FI = J$,
	\item \emph{(Compatibility with internal whiskering operators)}
	
	$$\begin{tikzcd}
		F{[}\mathfrak{B}{,}\mathfrak{C}{]} \arrow[rr, "F{[}\mathfrak{A}{,}-{]}"]
		\arrow[dd, "\chi_{\mathfrak{A}{,}\mathfrak{B}}"'] && F{[}{[}\mathfrak{A}{,}\mathfrak{B}{]}{,}{[}\mathfrak{A}{,}\mathfrak{C}{]}{]} \arrow[rr, " \chi_{{[}\mathfrak{A}{,}\mathfrak{B}{]}{,}{[}\mathfrak{A}{,}\mathfrak{C}{]}}"] &&  \{F{[}\mathfrak{A}{,}\mathfrak{B}{]}{,}F{[}\mathfrak{A}{,}\mathfrak{C}{]}\}
		\arrow[dd, "\{1{,}\chi_{\mathfrak{A}{,}\mathfrak{C}}\}"]
		\\
		\\
		\{F\mathfrak{B}{,}F\mathfrak{C}\} \arrow[rr, "\{F\mathfrak{A}{,}-\}"'] && \{\{F\mathfrak{A}{,}F\mathfrak{B}\}{,}\{F\mathfrak{A}{,}F\mathfrak{C}\}\} \arrow[rr, "\{\chi_{\mathfrak{A}{,}\mathfrak{B}}{,}1\}"'] && \{F{[}\mathfrak{A}{,}\mathfrak{B}{]}{,}\{F\mathfrak{A}{,}F\mathfrak{C}\}\}
	\end{tikzcd}$$
	
	\item \emph{(Compatibility with identity assigners)}
	
	$$\begin{tikzcd}
		FI \arrow[rr, "i_{\mathfrak{A}}"]\arrow[dd, equal] && {[}\mathfrak{A}{,}\mathfrak{A}{]}
		\arrow[dd, "\chi_{\mathfrak{A}{,}\mathfrak{A}}"]
		\\
		\\
		J \arrow[rr,"j_{\mathfrak{A}}"'] &&\{\mathfrak{A}{,}\mathfrak{A}\}
	\end{tikzcd}$$
	
	\item \emph{(Compatibility with constant map assigners)}.
	
	$$\begin{tikzcd}
		F\mathfrak{A} \arrow[dd, "k_{\mathfrak{A}}"']
		\arrow[rr, "Fc_{\mathfrak{A}}"] && F{[}I{,}\mathfrak{A}{]}
		\arrow[dd, "\chi_{I{,}\mathfrak{A}}"]
		\\
		\\
		\{J{,}F\mathfrak{A}\}\arrow[rr, equal]&& \{FI{,}F\mathfrak{A}\}
	\end{tikzcd}$$
\end{enumerate}
\end{definition}

\subsection{Transfors between $\mathbf{Gray}$-categories}\label{Subsection transfors between Gray-categories}

\begin{definition}\label{semi-strictly generated hom of Gray categories}
	Let $\mathfrak{A}$ and $\mathfrak{B}$ be $\mathbf{Gray}$-categories. The \emph{semi-strictly generated hom} $[\mathfrak{A}, \mathfrak{B}]_\text{ssg}$ is defined (Definition 5.3.1, \cite{Miranda strictifying operational coherences}) to be the locally full sub-$\mathbf{Gray}$-category of the $\mathbf{Gray}$-category $\mathbf{Tricat}_\text{s}\left(\mathfrak{A}, \mathfrak{B}\right)$ described in Remark \ref{higher cells between Gray categories as functions on globular input data} consisting of those trinatural transformations $p$ admitting some decomposition $p= p_{n}\circ ... \circ p_{1}$ such that each $p_{i}$ is semi-strict in the sense of Definition \ref{degrees of strictness trinatural transformations}.
\end{definition}

\begin{notation}
	For brevity, for the rest of this paper we will denote $[\mathfrak{A}, \mathfrak{B}]_\text{ssg}$ by just $[\mathfrak{A}, \mathfrak{B}]$.
\end{notation}

\begin{remark}\label{higher cells between Gray categories as functions on globular input data}
	In this paper we will make use of an explicit description of the $\mathbf{Gray}$-category $\mathbf{Tricat}_{s}\left(\mathfrak{A}, \mathfrak{B}\right)$ consisting of $\mathbf{Gray}$-functors, trinatural transformations, trimodifications, and perturbations, for $\mathbf{Gray}$-categories $\mathfrak{A}$ and $\mathfrak{B}$. This explicit description is deferred to the Appendix, to be given in Remark \ref{Gray Functor Gray category with weak higher cells}. The sub-$\mathbf{Gray}$-category on unital trinatural transformations is denoted $\mathbf{Ps}\left(\mathfrak{A}, \mathfrak{B}\right)_\text{s}$ in \cite{Bourke Lobbia Skew Approach to Enrichment for Gray-categories}. We describe the $m$-cells in $\mathbf{Tricat}_{s}\left(\mathfrak{A}, \mathfrak{B}\right)$, as sending $n$-cells in $\mathfrak{A}$ to certain data in $\mathfrak{B}$. In this description, we ignore the axioms which $m$-cells in $\mathbf{Tricat}_{s}\left(\mathfrak{A}, \mathfrak{B}\right)$ must satisfy and instead focus on the families of functions which can be used to determine when two well-defined $m$-cells are equal. This explicit description will be used to clarify the proofs in this paper.
	\begin{itemize}
		\item When $m=0$, $m$-cells are just $\mathbf{Gray}$-functors, and they send $n$-cells in $\mathfrak{A}$ to $n$-cells in $\mathfrak{B}$. Thus there are functions $F_{n}: \mathfrak{A}_{n} \rightarrow \mathfrak{B}_{n}$ between sets of $n$-globes for $n = 0, 1, 2, 3$. Moreover, two $\mathbf{Gray}$-functors $F: \mathfrak{A} \rightarrow \mathfrak{B}$ are equal if and only if they agree on $3$-cells, i.e. if and only if $F_{3} = G_{3}$.
		\item When $m = 1$, if we are already given $F, G: \mathfrak{A} \rightarrow \mathfrak{B}$ then an $m$-cell $p: F \rightarrow G$ is a trinatural transformation. These have the following behaviour on $n$-cells $\Phi \in \mathfrak{A}$. \begin{itemize}
			\item For $n = 0, 1, 2$ they send $\Phi$ to $m+n$-cells $p_\Phi$ in $\mathfrak{B}$. Thus there are functions $p_{n}:\mathfrak{A}_{n} \rightarrow \mathfrak{B}_{n + 1}$.
			\item When $n = 1$, the $2$-cell $p_{\Phi}$ is the left adjoint part of an adjoint equivalence. If we write $\mathfrak{B}_\text{eq}$ for the set of adjoint equivalences in hom-$2$-categories of $\mathfrak{B}$ and $l_\mathfrak{B}: \mathfrak{B}_\text{eq} \rightarrow \mathfrak{B}_{2}$ for the function which sends such an adjoint equivalence to its left adjoint, seen as a $2$-cell in $\mathfrak{B}$, then we may say that $p_{1}: \mathfrak{A}_{1} \rightarrow \mathfrak{B}_{2}$ factors through $l_\mathfrak{B}$ via some specified $p_{1}'$ which is part of the structure of $p$.
			\item If $p$ is not semi-strict then there is more data. Specifically, for an object $X$ in $\mathfrak{A}$ there is a $3$-cell in $\mathfrak{B}$ given by the unitor $p^{X}$, and for a composable pair \begin{tikzcd}
				X \arrow[r, "f"] & Y \arrow[r, "g"] & Z
			\end{tikzcd} in $\mathfrak{A}$ there is another $3$-cell in $\mathfrak{B}$ given by the compositor $p_{g, f}$. We will write these as functions $p_{u}: \mathfrak{A}_{0} \rightarrow \mathfrak{B}_{3}$ and $p_{c}: \mathfrak{A}_{1, 1} \rightarrow \mathfrak{B}_{3}$ respectively, where $\mathfrak{A}_{1, 1}$ denotes the set of composable pairs of $1$-cells in $\mathfrak{A}$.
		\end{itemize}
		\noindent Moreover, two trinatural transformations $p$ and $q$ are equal if and only if they agree on all of these data. In the particular case where $p$ and $q$ are semi-strict then to check $p = q$ it suffices to check that the functions $p_{n}, q_{n}: \mathfrak{A}_{n} \rightarrow \mathfrak{B}_{n+1}$ are equal for $n \neq 1$, and the functions $p_{1}', q_{1}': \mathfrak{A}_{1} \rightarrow \mathfrak{B}_\text{eq}$ are equal.
		\item When $n = 2$, if we are already given trinatural transformations $p, q: F \rightarrow G$ then a $2$-cell $\sigma: p \Rightarrow q$ is a trimodification. These send $n$-cells in $\mathfrak{A}$ to $n+2$-cells in $\mathfrak{B}$ for $n = 0, 1$. The equality of a pair of trimodifications $\sigma, \tau: p \Rightarrow q$ can be checked by checking the equality of functions $\sigma_{n}, \tau_{n}: \mathfrak{A}_{n} \rightarrow \mathfrak{B}_{n+2}$ for $n = 0, 1$.
		\item When $n = 3$, if we are already given trimodifications $\sigma, \tau: p \Rightarrow q$ then a $3$-cell between them is a perturbation. These have data given by a single function $\mathfrak{A}_{0} \rightarrow \mathfrak{B}_{3}$, and equality of perturbations is tantamount to equality of these functions.
	\end{itemize} 
	\noindent In summary, if we are already given two $\left(m - 1\right)$-cells in $\mathbf{Tricat}_{s}\left(\mathfrak{A}, \mathfrak{B}\right)$ then an $m$-cell $\Phi$ between them is determined by the functions $\Phi_{n}: \mathfrak{A}_{n} \rightarrow \mathfrak{B}_{n+m}$ for $m + n \leq 3$, unless $m = n = 1$. In this case $\Phi_{1}$ is a function from $\mathfrak{A}_{1}$ to $\mathfrak{B}_\text{eq}$, and if $m = 1$ and $\Phi$ is not semi-strict then there are extra functions $\Phi_{u}: \mathfrak{A}_{0} \rightarrow \mathfrak{B}_{3}$ and $\Phi_{c}: \mathfrak{A}_{1, 1} \rightarrow \mathfrak{B}_{3}$. There are various axioms that all these families of functions must satisfy to give a well-defined $m$-cell, but if we are already given two parallel $m$-cells which we know are well-defined and we need to check if they are equal, then it suffices to check the equality of each of the functions in the relevant families. Finally, given two other $\mathbf{Gray}$-categories $\mathfrak{C}$ and $\mathfrak{D}$, an $l$-cell $\Psi \in \mathbf{Tricat}_\text{s}(\mathfrak{B}, \mathfrak{C})$ and a $k$-cell $\Theta \in \mathbf{Tricat}_\text{s}(\mathfrak{C}, \mathfrak{D})$, there are analogous families of functions $\Psi_{n}: \mathfrak{B}_{n} \to \mathfrak{C}_{n+l}$ and $\Theta_{n}: \mathfrak{C}_{n} \to \mathfrak{D}_{n+k}$ for $0 \leq n \leq 3$, as well as functions $\Psi_{u}: \mathfrak{B}_{0} \to \mathfrak{C}_{3}$, $\Theta_{u}: \mathfrak{C}_{0} \to \mathfrak{D}_{3}$, $\Psi_{c}: \mathfrak{B}_{1, 1} \to \mathfrak{C}_{3}$ and $\Theta_{c}: \mathfrak{C}_{1, 1} \to \mathfrak{D}_{3}$. Associativity of the compositional structure of functions of this form, such as displayed below, is a key tool in the proofs of Theorem \ref{internal whiskering operators closed structure on Gray cat}, Proposition \ref{associativity condition underlying sets}, Proposition \ref{associativity condition for closed structure underlying sesquicategories} and Theorem \ref{Gray cat closed structure final theorem}.
	
	$$\begin{tikzcd}
		\mathfrak{A}_{n} \arrow[r, "\Phi_{m}"] & \mathfrak{B}_{n+m} \arrow[r, "\Psi_{l}"] & \mathfrak{C}_{n+m+l} \arrow[r, "\Theta_{k}"] & \mathfrak{D}_{n+m+l+k}
	\end{tikzcd}$$
\end{remark}

\noindent Although trinatural transformations are common in three-dimensional category theory, it may be difficult in practise to determine whether a given trinatural transformation can be expressed as a composite of semi-strict ones. Moreover, there may be many ways in which a given trinatural transformation can be decomposed in this way. We end this section by briefly describing two sources of examples of semi-strictly decomposable trinatural transformations. The first is from two-dimensional monad theory while the second is from three-dimensional monad theory.

\begin{example}\label{Examples of semi-strictly generated trinatural transformations}
	\noindent Recall from \cite{Coherent Approach to Pseudomonads} that there is a cofibrant $\mathbf{Gray}$-category $\mathbf{Psmnd}$ with the property that $\mathbf{Gray}$-functors $\mathbf{Psmnd} \to \mathbf{Gray}$ precisely correspond to pseudomonads on $2$-categories whose underlying endomorphisms are $2$-functors. An explicit presentation of $\mathbf{Psmnd}$ is given in Example 2.2.13 of \cite{Miranda PhD}. From this presentation we can deduce that there is a bijection between the following data:
	
	\begin{enumerate}
		\item Semi-strict trinatural transformations \begin{tikzcd}
			\mathbf{Psmnd} \arrow[r, bend left = 30, "(\mathcal{A}{,}S)"name=A]\arrow[r, bend right = 30, "(\mathcal{B}{,}T)"'name=B]& \mathbf{Gray} \arrow[from=A, to = B, Rightarrow, shift left = 1, shorten = 5]
		\end{tikzcd}.
		\item Morphisms of pseudomonads $(F, \phi): (\mathcal{A}, A) \to (\mathcal{B}, T)$ in the sense of Definition 2.1 of \cite{Formal Theory of Pseudomonads} whose underlying pseudonatural transformations $\phi: TF \rightarrow FS$ are part of specified adjoint equivalences.
	\end{enumerate} 

\noindent The component on the unique object $\star \in \mathbf{Psmnd}$ corresponds to the $2$-functor $F: \mathcal{A} \to \mathcal{B}$, while the component at the generating morphism $P: \star \to \star$ corresponds to $\phi$ and the components at the generating $2$-cells $\eta: 1_\star \Rightarrow P$ and $\mu: P\circ P \Rightarrow P$ in $\mathbf{Psmnd}$ respectively correspond to the invertible modifications $\tilde{\phi}$ and $\overline{\phi}$ of Definition 2.1 of \cite{Formal Theory of Pseudomonads}, and the axioms listed there are naturality in the generating $3$-cells $\alpha$ and $\lambda$, while the equation given in Proposition 2.2 of \cite{Formal Theory of Pseudomonads} is naturality in the generating $3$-cell $\rho$. In this way it is easy to verify the conditions in Lemma 4.3.1 part (1) of \cite{Miranda strictifying operational coherences} and see that there is indeed a bijection as claimed.
\\
\\
\noindent The composite of such trinatural transformations is of course semi-strictly decomposable, but it is typically not semi-strict. As such it is different to the composition structure of morphisms of pseudomonads described on page 29 of \cite{Formal Theory of Pseudomonads}. Instead, the latter corresponds to the alternative composition structure of semi-strict trinatural transformations out of cofibrant $\mathbf{Gray}$-categories, defined in Proposition 5.2.1 of \cite{Miranda strictifying operational coherences}.
\\
\\
\noindent We note that there are analogous sources of semi-strictly decomposable trinatural transformations with $\mathbf{Psmnd}$ replaced with any other cofibrant $\mathbf{Gray}$-category.
\end{example}

\begin{example}\label{triadjunction associated to a Gray-monad}
	Recall from \cite{Three dimensional monad theory} and Chapter 13 of \cite{Gurski Coherence in Three Dimensional Category Theory} that a $\mathbf{Gray}$-monad is a monad enriched over $\mathbf{Gray}$, and hence consists of a $\mathbf{Gray}$-category $\mathfrak{A}$, a $\mathbf{Gray}$-functor $T: \mathfrak{A} \to \mathfrak{A}$, and two $\mathbf{Gray}$-natural transformations \begin{tikzcd}
		1_\mathfrak{A} \arrow[r, "\eta"] &T & T^2 \arrow[l, "\mu"']
	\end{tikzcd} satisfying the usual monad laws on the nose. Recall also that to such structure one can associate a $\mathbf{Gray}$-category of pseudoalgebras $\mathfrak{A}^T$, and a $\mathbf{Gray}$-functor $U^T: \mathfrak{A}^T \to \mathfrak{A}$ which forgets the extra structure in a pseudoalgebra, pseudomorphism or pseudotransformation. It is shown in Proposition 6.2.3 of \cite{Miranda PhD} that $U^T$ is the right adjoint part of a three-dimensional adjunction, or triadjunction, in which the counit $\varepsilon$ is semi-strict. Indeed, unitality and compositionality of $\varepsilon$ is easy to see from the description of the $\mathbf{Gray}$-category structure of $\mathfrak{A}^T$ described in Section 13.3 of \cite{Gurski Coherence in Three Dimensional Category Theory}. Meanwhile, as a $\mathbf{Gray}$-natural transformation, the unit $\eta$ of this triadjunction is in particular a semi-strict trinatural transformation. An inspection of Lemma \ref{tautological condition for failure of semi-stricts to be closed under composition} reveals that this means that the target of the trimodification equivalence $u: 1_{U} \Rightarrow U\varepsilon. \eta_{U}$ is also a semi-strictly generated trinatural transformation. However, the triadjunction also involves an invertible perturbation mediating what is often called the `left tetrahedral equation' in the context of pseudoadjunctions.

$$\begin{tikzcd}[font=\fontsize{9}{6}, row sep = 15, column sep = 15]
	&&{}\arrow[dd, Rightarrow, shorten = 5, "Fu"]
	&&FU
	\arrow[rrdd, "\varepsilon"]
	\arrow[dddd, Rightarrow, shorten = 20, "\varepsilon_{\varepsilon}"]
	\\
	\\
	FU
	\arrow[rr, "F\eta U"]
	\arrow[rrrruu, bend left, "1_{FU}"]
	\arrow[rrrrdd, bend right, "1_{FU}"']
	&& FUFU
	\arrow[rrdd, "\varepsilon_{FU}"]
	\arrow[rruu, "FU\varepsilon"'] &&&&\mathfrak{A}^T \cong 1_{\varepsilon}
	\\
	&&=
	\\
	&&&&FU
	\arrow[rruu, "\varepsilon"']
\end{tikzcd}$$

 \noindent The trimodification $\varepsilon_{\varepsilon}$ is constructed as we will describe in Proposition \ref{interchanging tritransformations} part (1). Note that the trimodification appearing in the boundary of this perturbation is constructed as a composite of trimodifications involving trinatural transformations $\varepsilon.FU\varepsilon.F\eta U$ and $\varepsilon.\varepsilon FU.F\eta U$. These trinatural transformations are not semi-strict but are semi-strictly decomposable. Similarly, the left triadjunction axiom described in Example 3.3.3 of \cite{Miranda PhD} is an equation of pasting diagrams in the $2$-category $\mathbf{Gray}$-$\mathbf{Cat}(\mathfrak{A}^{T}, \mathfrak{A}^{T})(F, F)$. The objects appearing in these pasting diagrams are semi-strictly decomposable trinatural transformations, so that in fact this is an equation of $3$-cells in the semi-strictly generated internal hom $\mathbf{Gray}$-category $[\mathfrak{A}^{T}, \mathfrak{A}^{T}]$. 
\end{example}
\section{The data of the closed structure}\label{Data of the closed structure}

\subsection{The internal hom, unit, identity assigners and constant map assigners}

\noindent The goal of this subsection is to treat all aspects of the closed structure that do not involve the internal whiskering operators. The internal whiskering operators are significantly more complicated and will be treated in Subsections \ref{Subsection underlying functoriality of internal whiskering operators} and \ref{Subsection Gray functoriality of internal whiskering operators}. We first show in Proposition \ref{coherent hom on Gray Cat}, to follow, that the semi-strictly generated internal hom defines a functor $[-, ?]: \mathbf{Gray}\text{-}\mathbf{Cat}^\text{op} \times \mathbf{Gray}\text{-}\mathbf{Cat} \rightarrow \mathbf{Gray}\text{-}\mathbf{Cat}$.

\begin{proposition}\label{coherent hom on Gray Cat}
	Let $A: \mathfrak{A} \rightarrow \mathfrak{A}'$ and $B: \mathfrak{B} \rightarrow \mathfrak{B}'$ be $\mathbf{Gray}$-functors and consider the semi-strictly generated hom of Definition \ref{semi-strictly generated hom of Gray categories} part (2).
	
	\begin{enumerate}
		\item There is a $\mathbf{Gray}$-functor $[A, \mathfrak{B}]: [\mathfrak{A}', \mathfrak{B}] \rightarrow  [\mathfrak{A}, \mathfrak{B}]$.
		\item There is a $\mathbf{Gray}$-functor $[\mathfrak{A}, B]: [\mathfrak{A}, \mathfrak{B}] \rightarrow  [\mathfrak{A}, \mathfrak{B}']$.
		\item There is an equality of $\mathbf{Gray}$-functors $[A, \mathfrak{B}]\circ [\mathfrak{A}', B] = [\mathfrak{A},B]\circ [A, \mathfrak{B}]$.
		\item There is a functor $[-, ?]: \mathbf{Gray}$-$\mathbf{Cat}^\text{op}\times \mathbf{Gray}$-$\mathbf{Cat} \rightarrow \mathbf{Gray}$-$\mathbf{Cat}$ defined via $\left(\mathfrak{A}, \mathfrak{B}\right) \mapsto [\mathfrak{A}, \mathfrak{B}]$ on objects and via the $\mathbf{Gray}$-functor of part (3) on morphisms. 
	\end{enumerate}
\end{proposition}

\begin{proof}
	It is clear that if $j: F \rightarrow G: \mathfrak{A} \rightarrow \mathfrak{B}$ is a semi-strict trinatural transformation between $\mathbf{Gray}$-functors then $j_A$ is also semi-strict. In this way the $\mathbf{Gray}$-functor $\mathbf{Tricat}\left(A, \mathfrak{B}\right): \mathbf{Tricat}\left(\mathfrak{A}', \mathfrak{B}\right) \rightarrow \mathbf{Tricat}\left(\mathfrak{A}, \mathfrak{B}\right)$ in Theorem A.6 of \cite{Buhne PhD} defined by restriction along $A: \mathfrak{A} \rightarrow \mathfrak{A}'$, restricts to a $\mathbf{Gray}$-functor $[A, \mathfrak{B}]: [\mathfrak{A}', \mathfrak{B}] \rightarrow  [\mathfrak{A}, \mathfrak{B}]$. This proves part (1). For part (2), by Theorem A.3 of \cite{Buhne PhD}, there is at least a trihomomorphism $\mathbf{Tricat}\left(\mathfrak{A}, B\right): \mathbf{Tricat}\left(\mathfrak{A}, \mathfrak{B}\right) \rightarrow \mathbf{Tricat}\left(\mathfrak{A}, \mathfrak{B}'\right)$, and since $\mathbf{Gray}$-functors are closed under composition, this restricts to a trihomomorphism $\mathbf{Tricat}_{s}\left(\mathfrak{A}, B\right): \mathbf{Tricat}_{s}\left(\mathfrak{A}, \mathfrak{B}\right) \rightarrow \mathbf{Tricat}_{s}\left(\mathfrak{A}, \mathfrak{B}'\right)$. The coherences for this trihomomorphism are described in Proposition A.5 and Lemmas A.5 and A.6 of \cite{Buhne PhD}. They are given entirely in terms of the coherences for $B$, and hence they will all be identities since $B$ is assumed to be a $\mathbf{Gray}$-functor. Finally, it is once again clear that $Bj$ will also be semi-strict if $j$ is, so that this $\mathbf{Gray}$-functor restricts to $[\mathfrak{A}, B]: [\mathfrak{A}, \mathfrak{B}] \rightarrow  [\mathfrak{A}, \mathfrak{B}']$ as required for part (2). Part (3) uses the analysis in Remark \ref{higher cells between Gray categories as functions on globular input data}. It extends associativity of composition for $\mathbf{Gray}$-functors to associativity for composition of $\left(3, k\right)$-transfors in the case where $A \in [\mathfrak{A}, \mathfrak{A}']$ and $B \in [\mathfrak{B}, \mathfrak{B}']$ are just $\mathbf{Gray}$-functors but $\Phi \in [\mathfrak{A}', \mathfrak{B}]$ is of arbitrary dimension $k$, and can immediately be seen to follow from associativity of composition of functions between sets depicted below when $k \neq 1$.
	
	$$\begin{tikzcd}
		\mathfrak{A}_{m} \arrow[rr, "A_{m}"] && \mathfrak{A'}_{m} \arrow[rr, "\Phi_{m}"] && \mathfrak{B}_{k+m} \arrow[rr, "B_{k+m}"] &&\mathfrak{B'}_{k+m}
	\end{tikzcd}$$

\noindent When $k= m=1$ we also need to check that the rest of the adjoint equivalence agrees, but this is also follows by an easy inspection. Part (4) follows from part (3).
\end{proof}

\begin{remark}
	We warn the reader that there is some ambiguity in our notation $[\mathfrak{A}, -]$ for the extension $\mathbf{Gray}$-functors. Without context, $[\mathfrak{A}, -]$ could also be interpreted as the \emph{functor} from $\mathbf{Gray}$-$\mathbf{Cat}$ to $\mathbf{Gray}$-$\mathbf{Cat}$ given by the restriction of the functor from Proposition \ref{coherent hom on Gray Cat} along the functor $\mathbf{1} \rightarrow \mathbf{Gray}$-$\mathbf{Cat}$ which selects $\mathfrak{A}$. Nonetheless we think that this is a convenient abuse of notation with which it will be easy to avoid confusion simply by also denoting the source and target of $[\mathfrak{A}, -]$. The reason for this is that the functor $[\mathfrak{A}, -]: \mathbf{Gray}$-$\mathbf{Cat} \rightarrow \mathbf{Gray}$-$\mathbf{Cat}$ and the $\mathbf{Gray}$-functor $[\mathfrak{A}, -]: [\mathfrak{B}, \mathfrak{C}] \rightarrow [[\mathfrak{A}, \mathfrak{B}], [\mathfrak{A}, \mathfrak{C}]]$ have only one type of input data in common, namely $\mathbf{Gray}$-functors of the form $S: \mathfrak{B} \rightarrow \mathfrak{C}$. But we will need to define the extension $\mathbf{Gray}$-functor $[\mathfrak{A}, -]: [\mathfrak{B}, \mathfrak{C}] \rightarrow [[\mathfrak{A}, \mathfrak{B}], [\mathfrak{A}, \mathfrak{C}]]$ in such a way that it agrees with the functor $[\mathfrak{A}, -]: \mathbf{Gray}$-$\mathbf{Cat} \rightarrow \mathbf{Gray}$-$\mathbf{Cat}$ on this common datum. That is to say, both the functor $[\mathfrak{A}, -] : \mathbf{Gray}$-$\mathbf{Cat} \rightarrow \mathbf{Gray}$-$\mathbf{Cat}$ and the $\mathbf{Gray}$-functor  $[\mathfrak{A}, -]: [\mathfrak{B}, \mathfrak{C}] \rightarrow [[\mathfrak{A}, \mathfrak{B}], [\mathfrak{A}, \mathfrak{C}]]$ will assign $S: \mathfrak{A} \rightarrow \mathfrak{B}$ to the $\mathbf{Gray}$-functor $[\mathfrak{A}, S]: [\mathfrak{A}, \mathfrak{B}] \rightarrow [\mathfrak{A}, \mathfrak{C}]$.
\end{remark}

\noindent After Lemma \ref{Gray functors out of 1}, to follow, we will be ready to describe the identity assigning extra-natural transformation and prove that its components satisfy axiom (1) of Definition \ref{definition of a closed structure on a category}.

\begin{lemma}\label{Gray functors out of 1}
	Let $\mathbf{1}$ denote the terminal $\mathbf{Gray}$-category. The function $\text{ev}$ from the hom-set  $\mathbf{Gray}$-$\mathbf{Cat}\left(\mathbf{1},  \mathfrak{A}\right)$ to the set of objects $\mathfrak{A}_{0}$ which evaluates a $\mathbf{Gray}$-functor $\mathbf{1} \rightarrow \mathfrak{A}$ at the unique object in $\mathbf{1}$ is a bijection.
\end{lemma}

\begin{proof}
	It is clear that strictness of a $\mathbf{Gray}$-functor $\mathbf{X}: \mathbf{1} \rightarrow \mathfrak{A}$ means that it determines precisely the data of an object $X := \mathbf{X}\left(1\right)$ in $\mathfrak{A}$. There are no non-trivial coherences which could generate any extra data.
\end{proof}

\begin{proposition}\label{identity assigners in closed structure}
	For any $\mathbf{Gray}$-category $\mathfrak{A}$, consider the $\mathbf{Gray}$-functor $i_\mathfrak{A}: \mathbf{1} \rightarrow [\mathfrak{A}, \mathfrak{A}]$ which selects the identity $\mathbf{Gray}$-functor at $\mathfrak{A}$. \begin{enumerate}
		\item The assignation $\mathfrak{A} \mapsto i_\mathfrak{A}$ varies extra-naturally in $\mathfrak{A}$. 
		\item The function which sends a $\mathbf{Gray}$-functor $F: \mathfrak{A} \rightarrow \mathfrak{B}$ to the $\mathbf{Gray}$-functor depicted below defines a bijection between hom-sets $\mathbf{Gray}$-$\mathbf{Cat}\left(\mathfrak{A}, \mathfrak{B}\right) \rightarrow \mathbf{Gray}$-$\mathbf{Cat}\left(\mathbf{1}, [\mathfrak{A}, \mathfrak{B}]\right)$.
		
		$$\begin{tikzcd}
			\mathbf{1}\arrow[rr, "i_\mathfrak{A}"] && {[}\mathfrak{A}{,}\mathfrak{A}{]} \arrow[rr, "{[}1_\mathfrak{A}{,} F{]}"] && {[}\mathfrak{A}{,}\mathfrak{B}{]}
		\end{tikzcd}$$ 
	\end{enumerate}
\end{proposition}

\begin{proof}
	Part (1) amounts to the fact that for any $\mathbf{Gray}$-functor $F: \mathfrak{A} \rightarrow \mathfrak{B}$, the equality of $\mathbf{Gray}$-functors $1_\mathfrak{B}.F = F.1_\mathfrak{A}$ holds. This follows immediately from the left and right unit laws for $\mathbf{Gray}$-$\mathbf{Cat}$ as a category. Part (2) follows from Lemma \ref{Gray functors out of 1} and the right unit law $F.1_\mathfrak{A} = F$ for the category $\mathbf{Gray}$-$\mathbf{Cat}$.
\end{proof}

\noindent Finally, we describe the constant map assigning natural isomorphism.

\begin{proposition}\label{constant assigners closed structure}
	Let $\mathfrak{A}$ be a $\mathbf{Gray}$-category and let $\text{ev}_\mathfrak{A}: [\mathbf{1}, \mathfrak{A}] \rightarrow \mathfrak{A}$ be the $\mathbf{Gray}$-functor which evaluates each piece of data in $[\mathbf{1}, \mathfrak{A}]$ on the unique object of $\mathbf{1}$. Then
	
	\begin{enumerate}
		\item $\text{ev}_\mathfrak{A}$ has an inverse $c_\mathfrak{A}$ in $\mathbf{Gray}$-$\mathbf{Cat}$.
		\item The assignment $\mathfrak{A} \mapsto \text{ev}_\mathfrak{A}$ is natural in $\mathfrak{A}$.
		\item The assignment $\mathfrak{A} \mapsto c_\mathfrak{A}$ is the component at $\mathfrak{A}$ of a natural isomorphism.
	\end{enumerate}
\end{proposition}

\begin{proof}
	For part (1), we need to show that evaluation gives an isomorphism on $n$-cells for $0\leq n \leq 3$. For $n = 0$, this is Lemma \ref{Gray functors out of 1}. For $n = 1$, given a morphism $f: X \rightarrow Y$ in $\mathfrak{A}$ there is a semi-strict trinatural transformation $c_\mathfrak{A}\left(f\right): \mathbf{X} \rightarrow \mathbf{Y}$ whose $1$-cell component is $f$ and all other components are identities. It is clear that $\left(\text{ev}_\mathfrak{A}\circ c_\mathfrak{A}\right)\left(f\right) = f$. To see that $\left( c_\mathfrak{A} \circ \text{ev}_\mathfrak{A}\right)\left(\mathbf{f}\right) = \mathbf{f}$ for any semi-strict $\mathbf{f}$, it suffices to show that any semi-strict trinatural transformation $\mathbf{f}: \mathbf{X} \Rightarrow \mathbf{Y}: \mathbf{1} \rightarrow \mathfrak{A}$ is in fact strict. But by semi-strictness, the unitor and compositor for $\mathbf{f}$ are identities. Using this property for the unitor, we can see that the $2$-cell component of $\mathbf{f}$ at the unique identity morphism in $\mathbf{1}$ is the identity. By local pseudonaturality, the $3$-cell component at the unique identity $2$-cell in $\mathbf{1}$ is also the identity.
	\\
	\\
	For $n = 2$, we can similarly define $c_\mathfrak{A}\left(\phi: f \rightarrow g\right)$ to be the trimodification whose $2$-cell component is given by $\phi$ and whose $3$-cell component is the identity. Then clearly $\left(\text{ev}_\mathfrak{A}\circ c_\mathfrak{A}\right)\left(\phi\right) = \phi$. To see that $\left( c_\mathfrak{A} \circ \text{ev}_\mathfrak{A}\right)\left(\Phi\right) = \Phi$ holds for any trimodification $\Phi: \mathbf{f} \Rightarrow \mathbf{g}$ between semi-strict trinatural transformations, it suffices to see that the $3$-cell component of $\Phi$ must be the identity. But this follows from the unit law for $\Phi$. Finally, for $n = 3$ it is clear that a perturbation $\Phi \Rrightarrow \Psi$ is precisely the data of a $3$-cell $\Phi_{1} \Rrightarrow \Psi_{1}$. This proves part (1). Part (2) is precisely the definition of how $\mathbf{Gray}$-functors compose or whisker with data in $[\mathfrak{A}, \mathfrak{B}]$. Part (3) follows from part (2).
\end{proof}

\subsection{Underlying functoriality of the internal whiskering operators}\label{Subsection underlying functoriality of internal whiskering operators}

\noindent The fragment of the closed structure we have defined so far in Propositions \ref{coherent hom on Gray Cat}, \ref{identity assigners in closed structure} and \ref{constant assigners closed structure} has been relatively simple. In contrast, defining the extension $\mathbf{Gray}$-functors $[\mathfrak{A}, -]: [\mathfrak{B}, \mathfrak{C}] \rightarrow [[\mathfrak{A}, \mathfrak{B}], [\mathfrak{A}, \mathfrak{C}]]$ involves substantial calculation. Many of the lemmas needed are stated in \cite{Buhne PhD}, but with proof details omitted due to the size of the diagrams involved. In \cite{Bourke Lobbia Skew Approach to Enrichment for Gray-categories}, analogous lemmas are streamlined by analysing substitutions of the multimorphisms represented by their skew structures. Our treatment explicitly describes the various interchangers which mediate interactions between trinatural transformations and trimodifications. In the proofs of Propositions \ref{local modification axiom for interchager of trinatural transformations}, \ref{interchanging trimodification with tritransformation} and \ref{interchanging tritransformation with trimodification} we illustrate such proofs in full detail. Although we will omit some details in other proofs, we will give a strategy via which these details can be reconstructed by the sufficiently patient reader. As the reader may notice and as we will discuss in Remark \ref{proof strategy pasting diagram chases}, these proof strategies themselves can be systematically discovered.

\begin{remark}\label{Buhne's calculations}
	Many of the propositions to follow (eg. \ref{local modification axiom for interchager of trinatural transformations}, \ref{interchanging tritransformations}, \ref{interchanging trimodification with tritransformation} and \ref{interchanging tritransformation with trimodification}) are of the following form. They take as input a $m$-cell $\Phi$ in the $\mathbf{Gray}$-category $\mathbf{Tricat}_{s}\left(\mathfrak{A}, \mathfrak{B}\right)$ and an $n$-cell $\Psi$ in the $\mathbf{Gray}$-category $\mathbf{Tricat}_{s}\left(\mathfrak{B}, \mathfrak{C}\right)$ such that $2 \leq m + n \leq 3$, and return as output a $\left(m + n\right)$-cell in the $\mathbf{Gray}$-category $\mathbf{Tricat}_{s}\left(\mathfrak{A}, \mathfrak{C}\right)$ which can be viewed as a generalised interchanger between $\Psi$ and $\Phi$. Similarly, Propositions \ref{unitor of whiskering with p} and \ref{compositor of whiskering with p} establish data up to which these generalised interchangers respect the underlying category structure of $\mathbf{Tricat}_{s}\left(\mathfrak{A}, \mathfrak{B}\right)$. Each step in each of these proofs uses a property of trinatural transformations, trimodifications or perturbations. These properties will be described with reference to Appendix \ref{subsection weak higher maps between Gray functors}.
	\\
	\\
	\noindent The data constructed in Propositions \ref{unitor of whiskering with p} and \ref{compositor of whiskering with p} will feature as unitors and compositors of trinatural transformations $\mathbf{Tricat}_{s}\left(\mathfrak{A}, p\right): \mathbf{Tricat}_{s}\left(\mathfrak{A}, S\right) \rightarrow \mathbf{Tricat}_{s}\left(\mathfrak{A}, T\right)$ in Proposition \ref{Tricat A, p trinatural}. The assignment of trinatural transformations $p \mapsto \mathbf{Tricat}_\text{s}\left(\mathfrak{A}, p\right)$ is seen as the assignment on morphisms of a functor in Proposition \ref{functoriality of A, p}. This is shown to restrict to semi-strictly decomposable trinatural transformations, to define the assignment on morphisms of the underlying functor of the required internal whiskering operators $[\mathfrak{A}, -]: [\mathfrak{B}, \mathfrak{C}] \rightarrow [[\mathfrak{A}, \mathfrak{B}], [\mathfrak{A}, \mathfrak{C}]]$ in Corollary \ref{underlying functor of internal whiskering operator in closed structure}.
\end{remark}

\noindent Consider the following pair of trinatural transformations between $\mathbf{Gray}$-functors. 

\begin{equation}\label{pair of trinatural transformations between Gray-functors}
	\begin{tikzcd}
		\mathfrak{A} \arrow[rr, bend right, "G"' name = G]\arrow[rr, bend left, "F" name = F] &&\mathfrak{B} \arrow[rr, bend left, "S" name = S]\arrow[rr, bend right, "T"' name = T] && \mathfrak{C}
		\arrow[from = F, to = G, Rightarrow, shorten = 5, "j"]
		\arrow[from = S, to = T, Rightarrow, shorten = 5, "p"]
	\end{tikzcd}
\end{equation}

\noindent In Proposition \ref{interchanging tritransformations} we will show that these data give rise to an adjoint equivalence in $\mathbf{Tricat}_{s}\left(\mathfrak{A}, \mathfrak{C}\right)\left(SF, TG\right)$ whose left adjoint we will call $p_{j}: p_G.Sj \Rightarrow Tj.p_F$. Before doing this, we show in Proposition \ref{local modification axiom for interchager of trinatural transformations} that the maps $f \mapsto \left(p_{j}\right)_{f}$ are well-defined as modifications between the pseudonatural transformations between hom-$2$-categories of the trinatural transformations \begin{tikzcd}
	SF \arrow[rr, "pF"] && TF \arrow[rr, "Tj"] && TG
\end{tikzcd} and \begin{tikzcd}
	SF \arrow[rr, "Sj"] && SG \arrow[rr, "pG"] && TG\end{tikzcd}. See Remark \ref{Gray Functor Gray category with weak higher cells} part (4 a) for a description of the composite trinatural transformations.
 
\begin{proposition}\label{local modification axiom for interchager of trinatural transformations}
	Consider the family of $2$-cells in $\mathfrak{C}$ depicted below left, which are indexed by objects $X$ of $\mathfrak{A}$. Then the assignation of $f \in \mathfrak{A}\left(X, Y\right)$ to the $3$-cell in $\mathfrak{C}$ depicted below right defines a $2$-cell in $\mathbf{Gray}\left(\mathfrak{A}\left(X, Y\right), \mathfrak{C}\left(SFX, TGY\right)\right)$.
	\\
	
	\begin{tikzcd}[font=\fontsize{9}{6}]
		SFX \arrow[rrrr, "Sj_{X}"]\arrow[dddd, "p_{FX}"'] &&{}\arrow[dddd, Rightarrow, shorten = 20, "p_{j_{X}}"]
		&& SGX\arrow[dddd, "p_{GX}"]
		\\
		\\
		\\
		\\
		TFX \arrow[rrrr, "Tj_{X}"'] &&{}&& TGX&{}
	\end{tikzcd}\begin{tikzcd}[column sep= 25, font=\fontsize{9}{6}]
		p_{GY}.Sj_{Y}.SFf \arrow[rr, "p_{j_{Y}}.1"]
		\arrow[dd, "1.Sj_{f}"']
		\arrow[rrdd, "p_{j_{Y}.Ff}" description]
		&& Tj_{Y}.p_{FY}.SFf \arrow[dd, "1.p_{Ff}"]
		\\
		&{}\arrow[r, Rightarrow, shorten = 10, "p_{Ff{,}j_{Y}}^{-1}" shift left = 5]&{}
		\\
		p_{GY}.SGf.Sj_{X} \arrow[rr, shorten = 10, Rightarrow, "p_{j_{f}}"]
		\arrow[rrdd, "p_{Gf.j_{X}}" description]
		\arrow[dd, "p_{Gf}.1"'] 
		&& Tj_{Y}.TFf.p_{FX}\arrow[dd, "Tj_{f}.1"]
		\\
		{}\arrow[r, Rightarrow, shorten =10, "p_{j_{X}{,}Gf}"', shift right = 5]&{}
		\\
		TGf.p_{GX}.Sj_{X} \arrow[rr, "1. p_{j_{X}}"'] && TGf.Tj_{X}.p_{FX}
	\end{tikzcd}
	
\end{proposition}

\begin{proof}
	The modification axiom is a condition for every $\phi: f \Rightarrow g$ in $\mathfrak{A}$. On such data in $\mathfrak{A}$, the trinatural transformation $j$ determines the data in $\mathfrak{B}$ depicted below.
	
	$$
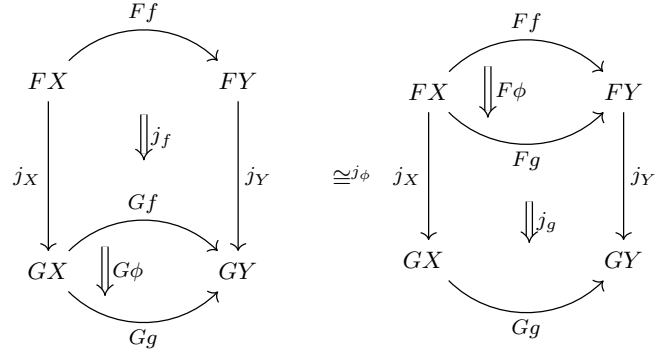
\begin{tikzcd}[column sep = 18, row sep = 20, font=\fontsize{9}{6}]
		FX
		\arrow[rr,bend left = 45, "Ff"]
		\arrow[dd,"j_{X}"']
		&
		{}
		\arrow[d,Rightarrow, shorten = 5,"j_{f}"]
		&
		FY
		\arrow[dd,"j_{Y}"]
		\\
		&{}&&\cong^{j_{\phi}}
		\\
		GX
		\arrow[rr,bend left = 45,"Gf" {name = C}]
		\arrow[rr, bend right = 45, "Gg"' {name = D}]
		&
		{}
		&
		GY
		\arrow[from =C, to =D, Rightarrow, shorten = 10, shift right = 5, "G\phi"]
	\end{tikzcd}	\begin{tikzcd}[column sep = 18, row sep = 20, font=\fontsize{9}{6}]
		FX
		\arrow[rr,bend left = 45, "Ff" {name = A}]
		\arrow[rr, bend right = 45, "Fg"' {name = B}]
		\arrow[dd,"j_{X}"']
		&
		&
		FY
		\arrow[dd,"j_{Y}"]
		\\
		&{}
		\arrow[d,Rightarrow, shorten = 5,"j_{g}"]
		\\
		GX\
		\arrow[rr, bend right = 45, "Gg"']
		&
		{}
		&
		GY &{}
		\arrow[from =A, to =B, Rightarrow, shorten = 10, shift right = 5, "F\phi"]
	\end{tikzcd} $$
	As we will now describe, the proof for the modification axiom uses all axioms for $p$ involving these data in $\mathfrak{B}$. Specifically, observe that the left hand side of the local modification axiom in Definition \ref{trimodification definition} is given by the diagram depicted below.
	
	$$
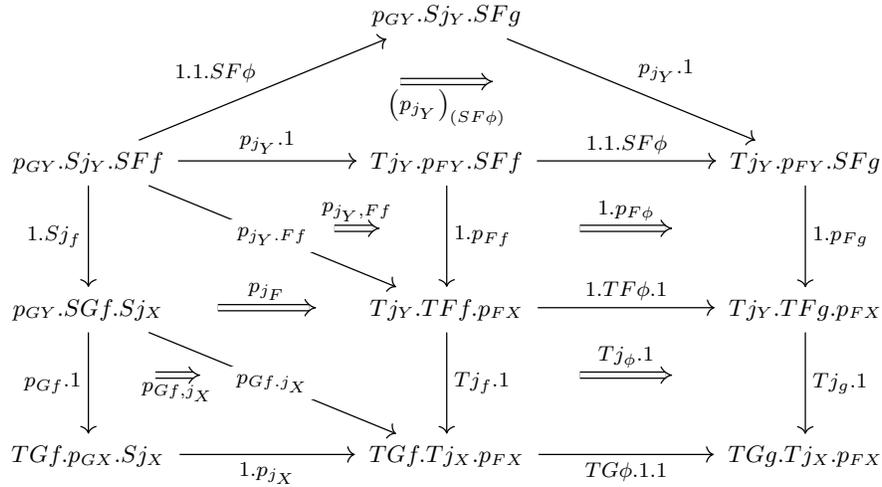
\begin{tikzcd}[font=\fontsize{9}{6}, row sep = 15]
		&&
		p_{GY}.Sj_{Y}.SFg
		\arrow[rrdd,"p_{j_Y}.1"]
		\\
		&
		{}
		\arrow[rr,Rightarrow,shorten=45,"{\left({p}_{j_Y}\right)}_{\left(SF\phi\right)}"']
		&&
		{}
		\\
		p_{GY}.Sj_{Y}.SFf
		\arrow[rruu,"1.1.SF\phi"]
		\arrow[rr,"p_{j_Y}.1"]
		\arrow[rrdd,"p_{j_Y.Ff}" description]
		\arrow[dd,"1.Sj_f"']
		&&
		Tj_{Y}.p_{FY}.SFf
		\arrow[rr,"1.1.SF\phi"]
		\arrow[dd,"1.p_{Ff}"]
		&&
		Tj_{Y}.p_{FY}.SFg
		\arrow[dd,"1.p_{Fg}"]
		\\
		&
		{}
		\arrow[r,Rightarrow,shorten=20,"p_{{j_Y}{,}{Ff}}"]
		&
		{}
		\arrow[rr,Rightarrow,shorten=45,"1.p_{F\phi}"]
		&&
		{}
		\\
		p_{GY}.SGf.Sj_{X}
		\arrow[rr,Rightarrow,shorten=15,"{p}_{j_F}"]
		\arrow[rrdd,"p_{{Gf}.{j_X}}" description]
		\arrow[dd,"p_{Gf}.1"']
		&&
		Tj_{Y}.TFf.p_{FX}
		\arrow[rr,"1.TF\phi.1"]
		\arrow[dd,"Tj_{f}.1"]
		&&
		Tj_{Y}.TFg.p_{FX}
		\arrow[dd,"Tj_{g}.1"]
		\\
		{}
		\arrow[r,Rightarrow,shorten=20,"p_{{Gf}{,}{j_X}}"']
		&
		{}
		&
		{}
		\arrow[rr,Rightarrow,shorten=45,"Tj_{\phi}.1"]
		&&
		{}
		\\
		TGf.p_{GX}.Sj_{X}
		\arrow[rr,"1.p_{j_X}"']
		&&
		TGf.Tj_{X}.p_{FX}
		\arrow[rr,"TG\phi.1.1"']
		&&
		TGg.Tj_{X}.p_{FX}
	\end{tikzcd}$$
	
	\noindent Now use $p$'s respect for the whiskering of $F\phi$ with $j_{Y}$.
	
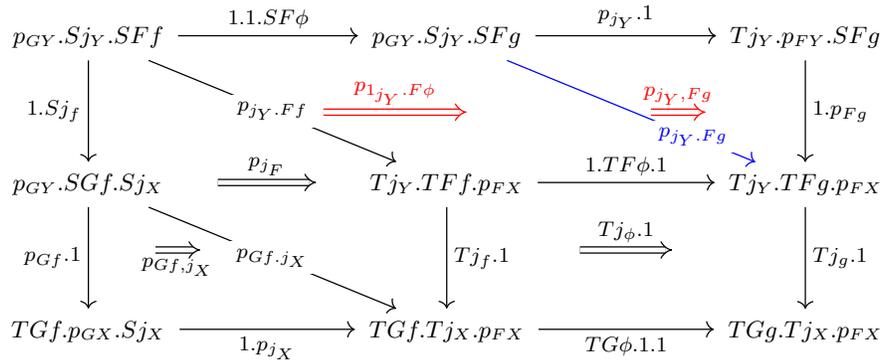
	$$\begin{tikzcd}[font=\fontsize{9}{6}, row sep = 15]
		p_{GY}.Sj_{Y}.SFf
		\arrow[rr,"1.1.SF\phi"]
		\arrow[rrr,Rightarrow,shorten=55,shift right=10, red,"p_{{1_{j_Y}}.{F\phi}}"]
		\arrow[rrdd,"p_{j_Y.Ff}" description]
		\arrow[dd,"1.Sj_f"']
		&&p_{GY}.Sj_{Y}.SFg
		\arrow[rr,"p_{j_Y}.1"]
		\arrow[rrdd,blue,"p_{{j_Y}.{Fg}}" description,near end]&
		{}
		\arrow[r,Rightarrow,red,shift right = 10,"p_{{j_Y}{,}{Fg}}",shorten=4]
		&
		Tj_{Y}.p_{FY}.SFg
		\arrow[dd,"1.p_{Fg}"]
		\\
		&
		{}
		&
		{}
		&&
		{}
		\\
		p_{GY}.SGf.Sj_{X}
		\arrow[rr,Rightarrow,shorten=15,"{p}_{j_F}"]
		\arrow[rrdd,"p_{{Gf}.{j_X}}" description]
		\arrow[dd,"p_{Gf}.1"']
		&&
		Tj_{Y}.TFf.p_{FX}
		\arrow[rr,"1.TF\phi.1"]
		\arrow[dd,"Tj_{f}.1"]
		&&
		Tj_{Y}.TFg.p_{FX}
		\arrow[dd,"Tj_{g}.1"]
		\\
		{}
		\arrow[r,Rightarrow,shorten=20,"p_{{Gf}{,}{j_X}}"']
		&
		{}
		&
		{}
		\arrow[rr,Rightarrow,shorten=45,"Tj_{\phi}.1"]
		&&
		{}
		\\
		TGf.p_{GX}.Sj_{X}
		\arrow[rr,"1.p_{j_X}"']
		&&
		TGf.Tj_{X}.p_{FX}
		\arrow[rr,"TG\phi.1.1"']
		&&
		TGg.Tj_{X}.p_{FX}
	\end{tikzcd}$$
	\noindent Now use local pseudonaturality of $p$ on the $3$-cell $j_{\phi}$.
	
	$$\begin{tikzcd}[font=\fontsize{9}{6}, row sep = 15]
		p_{GY}.Sj_{Y}.SFf
		\arrow[rr,"1.1.SF\phi"]
		\arrow[dd,"1.Sj_f"']
		\arrow[rr,Rightarrow,red,shorten=30,shift right = 10,"1.Sj_\phi"']
		&&
		p_{GY}.Sj_{Y}.SFg
		\arrow[rr,"p_{j_Y}.1"]
		\arrow[rrdd,"p_{{j_Y}.{Fg}}" description,near end]
		\arrow[dd,blue,"1.Sj_g"']
		&
		{}
		\arrow[r,Rightarrow,"p_{{j_Y}{,}{Fg}}"',shift right = 5, shorten=4]
		&
		Tj_{Y}.p_{FY}.SFg
		\arrow[dd,"1.p_{Fg}"]
		\\
		&
		{}
		&
		{}
		\arrow[r,Rightarrow,red,shorten=15,shift right = 4,"p_{j_g}"]
		&
		{}
		&
		{}
		\\
		p_{GY}.SGf.Sj_{X}
		\arrow[rr,blue,"1.SG\phi.1"]
		\arrow[rrdd,"p_{{Gf}.{j_X}}" description]
		\arrow[dd,"p_{Gf}.1"']
		&&
		\color{blue}p_{GY}.SGg.Sj_X
		\arrow[rrdd,blue,"p_{{Gg}.{j_X }}"]
		&&
		Tj_{Y}.TFg.p_{FX}
		\arrow[dd,"Tj_{g}.1"]
		\\
		{}
		\arrow[r,Rightarrow,shorten=20,"p_{{Gf}{,}{j_X}}"']
		&
		{}
		\arrow[rr,Rightarrow,red,shorten=45,"p_{{G\phi}.{1_{j_X}}}"']
		&
		{}
		&
		{}
		&
		{}
		\\
		TGf.p_{GX}.Sj_{X}
		\arrow[rr,"1.p_{j_X}"']
		&&
		TGf.Tj_{X}.p_{FX}
		\arrow[rr,"TG\phi.1.1"']
		&&
		TGg.Tj_{X}.p_{FX}
	\end{tikzcd}$$
	
	\noindent Finally, apply $p$'s respect for the whiskering of $G\phi$ with $j_{X}$, and observe that this results in the right hand side of the local modification axiom in Definition \ref{trimodification definition}, completing the proof.
	
	$$\begin{tikzcd}[font=\fontsize{9}{6}, row sep = 15]
		p_{GY}.Sj_{Y}.SFf
		\arrow[rr,"1.1.SF\phi"]
		\arrow[dd,"1.Sj_f"']
		\arrow[rr,Rightarrow,shorten=30,shift right = 10, "1.Sj_\phi"']
		&&
		p_{GY}.Sj_{Y}.SFg
		\arrow[rr,"p_{j_Y}.1"]
		\arrow[rrdd,"p_{{j_Y}.{Fg}}" description,near end]
		\arrow[dd,"1.Sj_g"']
		&
		{}
		\arrow[r,Rightarrow,shift right = 5, "p_{{j_Y}{,}{Fg}}"',shorten=4]
		&
		Tj_{Y}.p_{FY}.SFg
		\arrow[dd,"1.p_{Fg}"]
		\\
		&
		{}
		&
		{}
		\arrow[r,Rightarrow,shorten=15,shift right =2,"p_{j_g}"']
		&
		{}
		&
		{}
		\\
		p_{GY}.SGf.Sj_{X}
		\arrow[rr,"1.SG\phi.1"]
		\arrow[dd,"p_{Gf}.1"']
		&&
		p_{GY}.SGg.Sj_X
		\arrow[rrdd,"p_{{Gg}.{j_X }}"]
		\arrow[dd,blue,"p_{Gg}.1"']
		&&
		Tj_{Y}.TFg.p_{FX}
		\arrow[dd,"Tj_{g}.1"]
		\\
		{}
		\arrow[rr,Rightarrow,red,shorten=45,"p_{G\phi}.1"]
		&
		{}
		&
		{}
		\arrow[r,red,shorten=20,Rightarrow,"p_{{Gg}{,}{j_X}}"']
		&
		{}
		&
		{}
		\\
		TGf.p_{GX}.Sj_{X}
		\arrow[rrdd,"1.p_{j_X}"']
		\arrow[rr,blue,"TG{\phi}.1.1"]
		&&
		\color{blue}TGg.p_{GX}.Sj_{X}
		\arrow[rr,"1.p_{jX}",blue]
		&&
		TGg.Tj_{X}.p_{FX}
		\\
		&
		{}
		\arrow[rr,red,Rightarrow,shorten=30,"{\left(TG\phi\right)}_{\left(p_{j_X}\right)}"]
		&&
		{}
		\\
		&&
		TGf.Tj_{X}.p_{FX}
		\arrow[rruu,"TG\phi.1.1"']
	\end{tikzcd}$$
\end{proof}

\begin{proposition}\label{interchanging tritransformations}
	Let $j: F \rightarrow G$ and $p: S \rightarrow T$ be trinatural transformations as in Equation \ref{pair of trinatural transformations between Gray-functors}. 
	
	\begin{enumerate}
		\item There is a trimodification $p_{j}: p_G.Sj \Rightarrow Tj.p_F$, whose component at an object $X \in \mathfrak{A}$ is given by ${p}_{j_{X}}$ and whose $3$-cell component on a morphism $f: X \rightarrow Y$ is described in Proposition \ref{local modification axiom for interchager of trinatural transformations}.
		
		\item The trimodification $p_{j}$ is part of an adjoint equivalence in the $2$-category $\mathbf{Tricat}_{s}\left(\mathfrak{A}, \mathfrak{C}\right)\left(SF, TG\right)$.
		\item The trimodification $p_j$ is the identity if either $p$ is the identity or if $p$ is semi-strict and $j$ is the identity.
	\end{enumerate}
\end{proposition}

\begin{proof}
	Proposition \ref{local modification axiom for interchager of trinatural transformations} has established the local modification condition needed part (1). Indeed, parts (1) and (2) are part A.70 in Proposition A.12 of \cite{Buhne PhD}. However, their proofs omit details. We will describe how the steps in the proofs for the unit and composition laws can be discovered similarly to the proof of Proposition \ref{local modification axiom for interchager of trinatural transformations}.
	\\
	\\
	The unit law is a condition for each object $X \in \mathfrak{A}$. From this, $j$ determines the data in $\mathfrak{B}$ depicted below. 
	
	$$\begin{tikzcd}[font=\fontsize{9}{6}]
		FX \arrow[rr, "1_{FX}"]\arrow[dd, "j_{X}"'] &{}\arrow[dd, Rightarrow, shorten = 10, "j_{1_{X}}"]
		&FX \arrow[dd, "j_{X}"]
		&&FX \arrow[rr, "1_{FX}"]\arrow[dd, "j_{X}"'] &{}\arrow[dd, Rightarrow, shorten = 10, "1_{j_{X}}"]
		&FX\arrow[dd, "j_{X}"]
		\\
		&&&{\cong}^{j^{X}}
		\\
		GX \arrow[rr, "1_{GX}"'] &{}&
		GX
		&&GX\arrow[rr, "1_{GX}"']
		&{}&GX
	\end{tikzcd}$$
	\noindent The proof for the unit condition uses all axioms for $p$ determined by these data. This means that one uses \begin{itemize}
		\item The right unit axiom for $p$ on $j_{X}$,
		\item Local pseudonaturality of $p$ on $j^X$,
		\item The left unit axiom for $p$ on $j_{X}$.
	\end{itemize}
	
	\noindent Similarly, the composition axiom on \begin{tikzcd}
		X \arrow[r, "f"] & Y \arrow[r, "g"] &Z
	\end{tikzcd} uses all axioms for $p$ on the following data in $\mathfrak{B}$.
	
	$$\begin{tikzcd}[font=\fontsize{9}{6}]
		FX \arrow[rr, "Ff"]\arrow[dd, "j_{X}"'] &{}\arrow[dd, Rightarrow, shorten = 10, "j_{f}"]& FY \arrow[rr, "Fg"]\arrow[dd, "j_{Y}"] &{}\arrow[dd, Rightarrow, shorten = 10, "j_{g}"]& FZ \arrow[dd, "j_{Z}"]
		\\
		&&&&&\cong^{j_{g{,}f}}
		\\
		GX \arrow[rr, "GF"'] &{}& GY \arrow[rr, "Gg"'] &{}&GZ
	\end{tikzcd}\begin{tikzcd}[column sep = 35, row sep = 30]
		FX\arrow[dd, "j_{X}"'] \arrow[rr, "F\left(gf\right)"] 
		&{}\arrow[dd, Rightarrow, shorten = 10, "j_{gf}"]
		& FZ\arrow[dd, "j_{Z}"]
		\\
		\\
		GX \arrow[rr, "G\left(gf\right)"']&{}&GZ
	\end{tikzcd}$$
	
	\noindent This is to say, it uses
	
	\begin{itemize}
		\item The associativity axiom for $p$ on \begin{tikzcd}
			FX \arrow[rr, "Ff"] && FY \arrow[rr, "Fg"] && FZ \arrow[rr, "j_{Z}"] && GZ
		\end{tikzcd},
		\item $p$'s respect for the whiskering of $j_g$ with $Ff$,
		\item The associativity axiom for $p$ on \begin{tikzcd}
			FX \arrow[rr, "Ff"] && FY \arrow[rr, "j_{Y}"] && GY \arrow[rr, "Gg"] && GZ
		\end{tikzcd},
		\item $p$'s respect for the whiskering of $j_f$ with $Gg$
		\item The associativity axiom for $p$ on \begin{tikzcd}
			FX \arrow[rr, "j_{X}"] && GX \arrow[rr, "Gf"] && GY \arrow[rr, "Gg"] && GZ
		\end{tikzcd}
		\item Local pseudonaturality of $p$ on $j_{g{,}f}$.
	\end{itemize}
	\noindent In this way one can construct the proofs for part (1). In the special case where $p$ and $j$ are both semi-strict trinatural transformations, the proof is much simpler. The unit law holds immediately while the proof of the composition law only uses the whiskering laws mentioned. This simpler case is discussed in Remark \ref {proof strategy pasting diagram chases}.
	\\
	\\
	\noindent For part (2), one can similarly construct a trimodification $p_{j}^{*}$ with $2$-cell component at $X$ given by $p_{j_{X}}^{*}$. Then one can construct invertible perturbations $1 \rightarrow p_{j}^{*}p_{j}$ and $p_{j}p_{j}^{*} \rightarrow 1$, each of whose component at $X$ will be given by the unit and counit of $p_{j_{X}} \dashv p_{j_{X}}^{*}$ respectively. The perturbation axioms at $f: X \rightarrow Y$ follow from the modification axioms for $p^\eta$ and $p^\varepsilon$ at $j_{f}$ and $j_{f}^{*}$ respectively. The triangle identities for the adjoint equivalence follow from those for the adjunction $p_{j_{X}}^{*} \dashv p_{j_{X}}$. Part (3) is clear from the definition of $p_{j}$.
\end{proof}

\noindent Proposition \ref{interchanging tritransformations} describes how trinatural transformations interchange with one another. Propositions \ref{unitor of whiskering with p} and \ref{compositor of whiskering with p}, to follow, establish perturbations up to which this interchange respects the underlying category structure of $\mathbf{Tricat}_\text{s}\left(\mathfrak{A}, \mathfrak{B}\right)$.

\begin{proposition}\label{unitor of whiskering with p}
	Let $p: S \rightarrow T$ and $F$ be as in Equation \ref{pair of trinatural transformations between Gray-functors}, and consider the identity trinatural transformation $1_F$. Then there is a perturbation $1_{p_{F}} \Rrightarrow p_{1_{F}}$ whose component at $X$ is given by $p^{FX}$.
\end{proposition}

\begin{proof}
	The trinatural transformation $1_{F}: F \rightarrow F$ has the following $2$-cell component in $\mathfrak{B}$ on the morphism $f: X \rightarrow Y$ in $\mathfrak{A}$.
	
	$$\begin{tikzcd}[font=\fontsize{9}{6}]
		FX \arrow[rr, "Ff"]\arrow[dd, "1_{FX}"'] &{}\arrow[dd, Rightarrow, shorten = 10, "1_{Ff}"]& FY\arrow[dd, "1_{FY}"]
		\\
		\\
		FX \arrow[rr, "Ff"'] &{}& FY
	\end{tikzcd}$$
	
	\noindent The perturbation axiom at $f: X \rightarrow Y$ follows from the left and right unit laws for the trinatural transformation $p$ at the morphism $Ff$ in $\mathfrak{B}$.
\end{proof}

\begin{proposition}\label{compositor of whiskering with p}
	Let $j: F \rightarrow G$ and $p: S \rightarrow T$ be as in Equation \ref{pair of trinatural transformations between Gray-functors}, and let $p_j: p_{G}.Sj \Rightarrow Tj.p_{F}$ be the trimodification defined in Proposition \ref{interchanging tritransformations} part (1). Let $j': G \rightarrow H$ be another trinatural transformation and let $p_{j'}: p_{G}.Sj' \Rightarrow Tj'.p_{F}$ be the trimodification constructed similarly. Then there is an invertible perturbation as depicted below, whose component at $X \in \mathfrak{A}$ is given by the compositor $p_{j'_{X}, j_{X}}$.

	$$\begin{tikzcd}[font=\fontsize{9}{6}]
		&&Tj'.p_{G}.Sj
		\arrow[rrdd, "1.p_{j}"]
		\arrow[dd, shorten = 10, Rightarrow, "p_{j'{,}j}"]
		\\
		\\
		p_{H}.Sj'.Sj \arrow[rruu, "p_{j'}.1"]\arrow[rrrr, "p_{j'j}"']&&{}
		&& Tj'.Tj.p_{F}
	\end{tikzcd}$$
	
\end{proposition}

\begin{proof}
	The perturbation axiom is a condition on a morphism $f: X \rightarrow Y$ in $\mathfrak{A}$. This determines the data in $\mathfrak{B}$ depicted below.

	$$\begin{tikzcd}[font=\fontsize{9}{6}]
		FX \arrow[rr, "j_{X}"]\arrow[dd, "Ff"'] &{}\arrow[dd, Leftarrow, shorten = 10, "j_{f}"]& GX \arrow[rr, "j'_{X}"]\arrow[dd, "Gf"] &{}\arrow[dd, Leftarrow, shorten = 10, "j'_{f}"]& HX \arrow[dd, "Hf"]
		\\
		\\
		FY \arrow[rr, "j_{Y}"'] &{}& GY \arrow[rr, "j'_{Y}"'] &{}&HY
	\end{tikzcd}$$
	
	\noindent One uses all axioms for $p$ on this data in $\mathfrak{B}$. Note that the individual steps in this proof will be similar to those in the proof of the composition axiom for $p_{j}$ in Proposition \ref{interchanging tritransformations} except that it will not involve using local pseudonaturality of $p$.
\end{proof}

\begin{proposition}\label{interchanging trimodification with tritransformation}
	Let $j: F \rightarrow G$ and $p: S \rightarrow T$ be as in Equation \ref{pair of trinatural transformations between Gray-functors}, and let $\alpha: j \rightarrow k$ be a trimodification. \begin{enumerate}
		\item There is a perturbation as depicted below, whose component at $X$ is given by $p_{\alpha_{X}}$.

		$$\begin{tikzcd}[column sep = 18, row sep = 15, font=\fontsize{9}{6}]
			SF
			\arrow[rr,bend left = 45, "Sj"]
			\arrow[dd,"p_F"']
			&
			{}
			\arrow[d,Rightarrow, shorten = 5,"p_{j}"]
			&
			SG
			\arrow[dd,"p_G"]
			\\
			&{}&&\cong^{p_{\alpha}}
			\\
			TF
			\arrow[rr,bend left = 45,"Tj" {name = C}]
			\arrow[rr, bend right = 45, "Tk"' {name = D}]
			&
			{}
			&
			TG
			\arrow[from =C, to =D, Rightarrow, shorten = 10, shift right = 5, "T\alpha"]
		\end{tikzcd}	\begin{tikzcd}[column sep = 18, row sep = 15, font=\fontsize{9}{6}]
			SF
			\arrow[rr,bend left = 45, "Sj" {name = A}]
			\arrow[rr, bend right = 45, "Sk"' {name = B}]
			\arrow[dd,"p_F"']
			&
			&
			SG
			\arrow[dd,"p_G"]
			\\
			&{}
			\arrow[d,Rightarrow, shorten = 5,"p_{k}"]
			\\
			TF\
			\arrow[rr, bend right = 45, "Tk"']
			&
			{}
			&
			TG &{}
			\arrow[from =A, to =B, Rightarrow, shorten = 10, shift right = 5, "S\alpha"]
		\end{tikzcd} $$
		\item $p_\alpha$ is the identity if either $\alpha$ is strict of if $p$ is locally strict. In particular it is the identity if either $p$ is the identity or if $\alpha$ is the identity.
	\end{enumerate} 
\end{proposition}

\begin{proof}
	Part (1) is part A.74 in Proposition A.13 of \cite{Buhne PhD}. The idea of the proof of the perturbation axiom is the same as the ideas in the proof of Proposition \ref{interchanging tritransformations}. The perturbation axiom is a condition for a morphism $f: X \rightarrow Y$ in $\mathfrak{A}$. Under $\alpha$, this determines the data in $\mathfrak{B}$ depicted below. As we will describe, the proof uses all of the axioms for $p$ on this data.

	$$\begin{tikzcd}[column sep = 18, row sep = 15, font=\fontsize{9}{6}]
		FX
		\arrow[rr,bend left = 45, "k_{X}"]
		\arrow[dd,"Ff"']
		&
		{}
		\arrow[d,Leftarrow, shorten = 5,"k_{f}"]
		&
		GX
		\arrow[dd,"Gf"]
		\\
		&{}&&\cong^{\alpha_{f}}
		\\
		FY
		\arrow[rr,bend left = 45,"k_{Y}" {name = C}]
		\arrow[rr, bend right = 45, "j_{Y}"' {name = D}]
		&
		{}
		&
		GY
		\arrow[from =C, to =D, Leftarrow, shorten = 10, shift right = 5, "\alpha_{Y}"]
	\end{tikzcd}	\begin{tikzcd}[column sep = 18, row sep = 15, font=\fontsize{9}{6}]
		FX
		\arrow[rr,bend left = 45, "k_{X}" {name = A}]
		\arrow[rr, bend right = 45, "j_{X}"' {name = B}]
		\arrow[dd,"Ff"']
		&
		&
		GX
		\arrow[dd,"Gf"]
		\\
		&{}
		\arrow[d,Leftarrow, shorten = 5,"j_{f}"]
		\\
		FY\
		\arrow[rr, bend right = 45, "j_{Y}"']
		&
		{}
		&
		GY &{}
		\arrow[from =A, to =B, Leftarrow, shorten = 10, shift right = 5, "\alpha_{X}"]
	\end{tikzcd} $$
	\noindent Start with the right hand side of the perturbation axiom in Definition \ref{Definition perturbation}, which is given by the pasting depicted below.
	
	$$
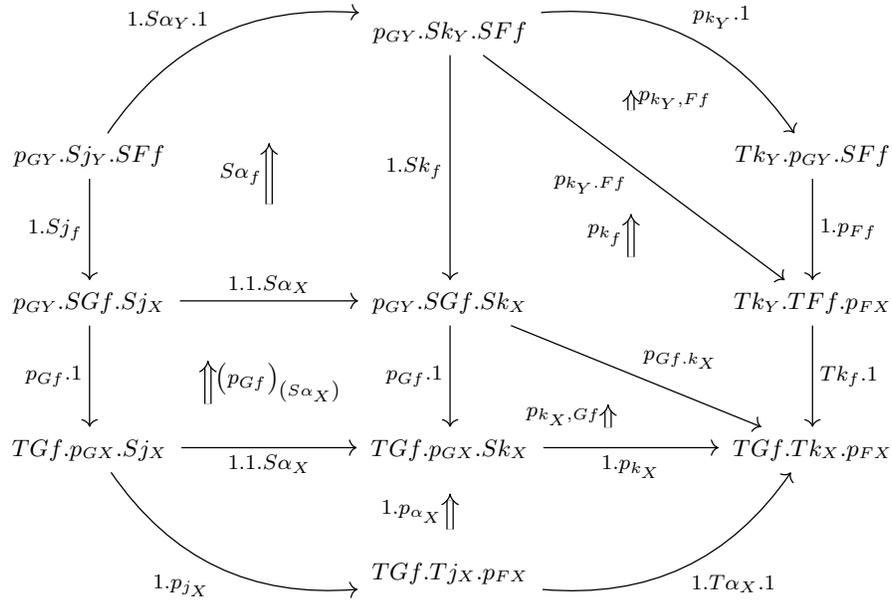
\begin{tikzcd}[row sep = 15, font=\fontsize{9}{6}]
		&
		{}
		&
		{p_{GY}}.{Sk_Y}.{SFf}
		\arrow[rrdd,bend left,"{p_{k_Y}}.1"]
		\arrow[rrdddd,"p_{{k_Y}.{Ff}}"']
		\arrow[dddd,"1.{Sk_f}"']
		&
		{}
		\\
		\\
		{p_{GY}}.{Sj_Y}.{SFf}
		\arrow[rruu,bend left,"1.{S{\alpha}_Y}.1"]
		\arrow[dd,"1.{Sj_f}"']
		&{}&
		&{}
		\arrow[uu,Rightarrow,shorten=15,"p_{{k_Y}{,}{Ff}}"']
		&
		{Tk_Y}.{p_{GY}}.{SFf}
		\arrow[dd,"1.{p_{Ff}}"]
		\\
		&{}
		\\
		{p_{GY}}.{SGf}.{Sj_X}
		\arrow[dd,"{p_{Gf}}.1"']
		\arrow[rr,"1.1.{S{\alpha}_X}"]
		&
		{}
		\arrow[uuuu,Rightarrow,shorten=35,"S{\alpha}_f"]
		&
		{p_{GY}}.{SGf}.{Sk_X}
		\arrow[rrdd,"p_{{Gf}.{k_X}}"]
		\arrow[dd,"{p_{Gf}}.1"']
		&{}
		\arrow[uu,Rightarrow,shorten=15,"p_{k_f}"]
		&
		{Tk_Y}.{TFf}.{p_{FX}}
		\arrow[dd,"{Tk_f}.1"]
		\\
		&&&
		{}
		\\
		{TGf}.{p_{GX}}.{Sj_X}
		\arrow[rrdd,bend right,"1.{p_{j_X}}"']
		\arrow[rr,"1.1.{S{\alpha}_X}"']
		&{}
		\arrow[uu,Rightarrow,shorten=15,"{\left(p_{Gf}\right)}_{\left(S{\alpha}_X\right)}"',shift left=8]
		&
		{TGf}.{p_{GX}}.{Sk_X}
		\arrow[rr,"1.{p_{k_X}}"']
		&
		{}
		\arrow[u,Rightarrow,shorten=6,"p_{{k_X}{,}{Gf}}",shift left=3]
		&
		{TGf}.{Tk_X}.{p_{FX}}
		\\
		\\
		&{}
		&
		{TGf}.{Tj_X}.{p_{FX}}
		\arrow[rruu,bend right,"1.{T{\alpha}_X}.1"']
		\arrow[uu,Rightarrow,shorten=8,"1.{p_{{\alpha}_X}}"]
		&{}
	\end{tikzcd}$$
	
	\noindent Use $p$'s respect for whiskering $\alpha_{X}$ with $Gf$ to reduce this to the pasting depicted below.
	
	$$
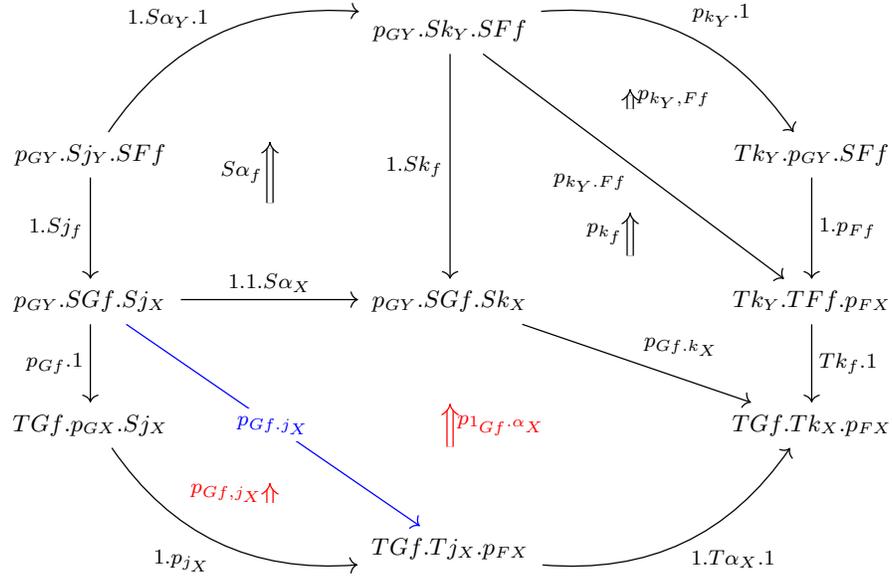
\begin{tikzcd}[row sep = 15, font=\fontsize{9}{6}]
		&
		{}
		&
		{p_{GY}}.{Sk_Y}.{SFf}
		\arrow[rrdd,bend left,"{p_{k_Y}}.1"]
		\arrow[rrdddd,"p_{{k_Y}.{Ff}}"']
		\arrow[dddd,"1.{Sk_f}"']
		&
		{}
		\\
		\\
		{p_{GY}}.{Sj_Y}.{SFf}
		\arrow[rruu,bend left,"1.{S{\alpha}_Y}.1"]
		\arrow[dd,"1.{Sj_f}"']
		&{}&
		&{}
		\arrow[uu,Rightarrow,shorten=15,"p_{{k_Y}{,}{Ff}}"']
		&
		{Tk_Y}.{p_{GY}}.{SFf}
		\arrow[dd,"1.{p_{Ff}}"]
		\\
		&{}
		\\
		{p_{GY}}.{SGf}.{Sj_X}
		\arrow[rrdddd,"p_{{Gf}.{j_X}}" description,blue]
		\arrow[dd,"{p_{Gf}}.1"']
		\arrow[rr,"1.1.{S{\alpha}_X}"]
		&
		{}
		\arrow[uuuu,Rightarrow,shorten=35,"S{\alpha}_f"]
		&
		{p_{GY}}.{SGf}.{Sk_X}
		\arrow[rrdd,"p_{{Gf}.{k_X}}"]
		&{}
		\arrow[uu,Rightarrow,shorten=15,"p_{k_f}"]
		&
		{Tk_Y}.{TFf}.{p_{FX}}
		\arrow[dd,"{Tk_f}.1"]
		\\
		\\
		{TGf}.{p_{GX}}.{Sj_X}
		\arrow[rrdd,bend right,"1.{p_{j_X}}"']
		&{}
		&&&
		{TGf}.{Tk_X}.{p_{FX}}
		\\
		\\
		&{}\arrow[uu, shorten = 15, Rightarrow, "p_{Gf{,}j_{X}}",red]
		&
		{TGf}.{Tj_X}.{p_{FX}}
		\arrow[rruu,bend right,"1.{T{\alpha}_X}.1"']
		\arrow[uuuu,Rightarrow,shorten=30,red,"p_{{1_{Gf}}.{{\alpha}_X}}"']
		&{}
	\end{tikzcd}$$
	
	\noindent Now use local pseudonaturality of $p$ on $\alpha_{f}$.
	
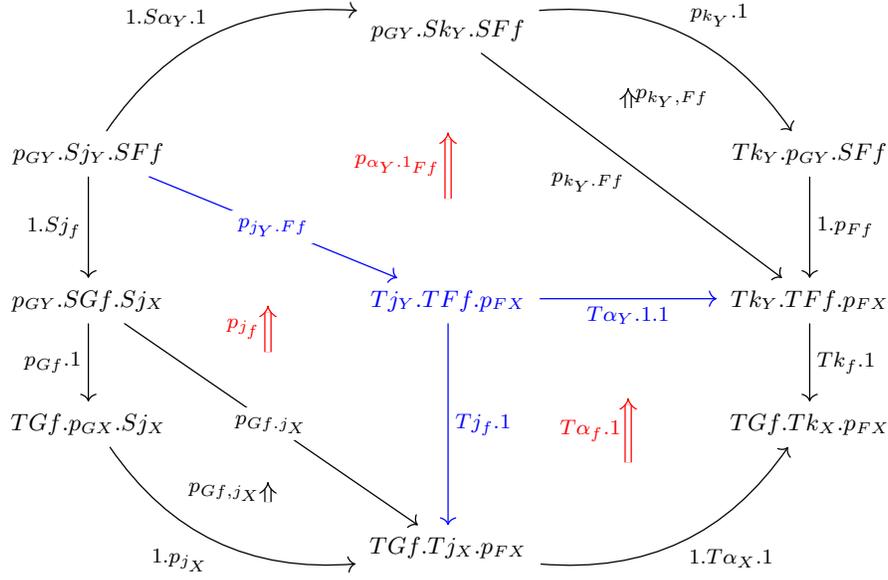
	$$\begin{tikzcd}[row sep = 15, font=\fontsize{9}{6}]
		&&
		{p_{GY}}.{Sk_Y}.{SFf}
		\arrow[rrdd,bend left,"{p_{k_Y}}.1"]
		\arrow[rrdddd,"p_{{k_Y}.{Ff}}"']
		&
		{}
		\\
		\\
		{p_{GY}}.{Sj_Y}.{SFf}
		\arrow[rruu,bend left,"1.{S{\alpha}_Y}.1"]
		\arrow[rrdd,"p_{{j_Y}.{Ff}}" description,blue]
		\arrow[dd,"1.{Sj_f}"']
		&{}&
		&{}
		\arrow[uu,Rightarrow,shorten=15,"p_{{k_Y}{,}{Ff}}"']
		&
		{Tk_Y}.{p_{GY}}.{SFf}
		\arrow[dd,"1.{p_{Ff}}"]
		\\
		&{}
		\\
		{p_{GY}}.{SGf}.{Sj_X}
		\arrow[rrdddd,"p_{{Gf}.{j_X}}" description]
		\arrow[dd,"{p_{Gf}}.1"']
		&&
		\color{blue}{Tj_Y}.{TFf}.{p_{FX}}
		\arrow[rr,"{T{\alpha}_Y}.1.1"',blue]
		\arrow[dddd,"{Tj_f}.1",blue]
		\arrow[uuuu,shorten=30,Rightarrow,red,"p_{{{\alpha}_Y}.{1_{Ff}}}"]
		&{}
		&
		{Tk_Y}.{TFf}.{p_{FX}}
		\arrow[dd,"{Tk_f}.1"]
		\\
		\\
		{TGf}.{p_{GX}}.{Sj_X}
		\arrow[rrdd,bend right,"1.{p_{j_X}}"']
		&{}\arrow[uuu, Rightarrow, shorten = 25, "p_{j_{f}}",red]
		&&&
		{TGf}.{Tk_X}.{p_{FX}}
		\\
		\\
		&{}\arrow[uu, shorten = 15, Rightarrow, "p_{Gf{,}j_{X}}"]
		&
		{TGf}.{Tj_X}.{p_{FX}}
		\arrow[rruu,bend right,"1.{T{\alpha}_X}.1"']
		&{}\arrow[uuuu, shorten = 30, Rightarrow, "T\alpha_{f}.1",red]
	\end{tikzcd}$$
	
	\noindent Finally, use $p$'s respect for whiskering $\alpha_{Y}$ with $Ff$ to reduce this to the pasting depicted below and observe that this is the required left hand side of the perturbation axiom in Definition \ref{Definition perturbation}, completing the proof for part (1).
	
	$$\begin{tikzcd}[row sep = 15, font=\fontsize{9}{6}]
		&&
		{p_{GY}}.{Sk_Y}.{SFf}
		\arrow[rrdd,bend left,"{p_{k_Y}}.1"]
		\\
		\\
		{p_{GY}}.{Sj_Y}.{SFf}
		\arrow[rruu,bend left,"1.{S{\alpha}_Y}.1"]
		\arrow[rr,blue,"{p_{j_Y}}.1"]
		\arrow[rrdd,"p_{{j_Y}.{Ff}}" description]
		\arrow[dd,"1.{Sj_f}"']
		&{}&
		\color{blue}{Tj_Y}.{p_{GY}}.{SFf}
		\arrow[rr,blue,"{T{\alpha}_Y}.1.1"]
		\arrow[dd,blue,"1.{p_{Ff}}"]
		\arrow[Rightarrow,uu,red,shorten=8, "p_{\alpha_{Y}}.1"]
		&{}&
		{Tk_Y}.{p_{GY}}.{SFf}
		\arrow[dd,"1.{p_{Ff}}"]
		\\
		&{}\arrow[u, shorten = 5, Rightarrow, red, "p_{j_{Y}{,}Ff}"']
		\\
		{p_{GY}}.{SGf}.{Sj_X}
		\arrow[rrdddd,"p_{{Gf}.{j_X}}" description]
		\arrow[dd,"{p_{Gf}}.1"']
		&&
		{Tj_Y}.{TFf}.{p_{FX}}
		\arrow[rr,"{T{\alpha}_Y}.1.1"']
		\arrow[dddd,"{Tj_f}.1"]
		&{}\arrow[uu, Rightarrow, shorten = 15, red, "{\left(T\alpha_{Y}\right)}_{\left(p_{Ff}\right)}"']&
		{Tk_Y}.{TFf}.{p_{FX}}
		\arrow[dd,"{Tk_f}.1"]
		\\
		\\
		{TGf}.{p_{GX}}.{Sj_X}
		\arrow[rrdd,bend right,"1.{p_{j_X}}"']
		&{}\arrow[uuu, Rightarrow, shorten = 25, "p_{j_{f}}"]
		&&&
		{TGf}.{Tk_X}.{p_{FX}}
		\\
		\\
		&{}\arrow[uu, shorten = 15, Rightarrow, "p_{Gf{,}j_{X}}"]
		&
		{TGf}.{Tj_X}.{p_{FX}}
		\arrow[rruu,bend right,"1.{T{\alpha}_X}.1"']
		&{}\arrow[uuuu, shorten = 30, Rightarrow, "T\alpha_{f}.1"]
	\end{tikzcd}$$
	
	\noindent Part (2) is clear from the definition of $p_\alpha$, since local pseudonaturality of $p$ forces its $3$-cell component at an identity $2$-cell to also be an identity, and if $p$ is locally strict says that its component $3$-cell at any $2$-cell is the identity, so in particular $p_{\alpha_{X}}$ will be the identity.
\end{proof}

\begin{remark}\label{proof strategy pasting diagram chases}
	The proofs of Propositions \ref{local modification axiom for interchager of trinatural transformations}, \ref{interchanging tritransformations}, \ref{unitor of whiskering with p}, \ref{compositor of whiskering with p} and \ref{interchanging trimodification with tritransformation} all follow a similar pattern. A similar proof will also be needed in Proposition \ref{interchanging tritransformation with trimodification} when we will need to describe a $3$-cell in $\mathbf{Tricat}_{s}\left(\mathfrak{A}, \mathfrak{C}\right)$ mediating interchange between a $1$-cell in $\mathbf{Tricat}_{s}\left(\mathfrak{A}, \mathfrak{B}\right)$ and a $2$-cell in $\mathbf{Tricat}_{s}\left(\mathfrak{B}, \mathfrak{C}\right)$. We describe the pattern in these proofs.
	\\
	\\
	\noindent For each axiom, a certain family of pairs of pasting diagrams in the hom-$2$-categories of $\mathfrak{C}$ are required to be equal. Each such family of pairs of pasting diagrams is indexed by some data in $\mathfrak{A}$. These data in $\mathfrak{A}$ determine some data in $\mathfrak{B}$. Finally, the highest dimensional $\left(3, k\right)$ transfor from $\mathfrak{B}$ to $\mathfrak{C}$ appearing in the statement of the Proposition has various axioms stipulated for these data in $\mathfrak{B}$. Those are the main axioms needed in the proof that the two pasting diagrams in $\mathfrak{C}$ are equal. Properties that we have omitted mention of but which are also implicit in our proofs include the following.
	\begin{itemize}
		\item $\mathbf{Gray}$-functoriality of $S$ and $T$,
		\item The composition and unit axioms for pseudonatural transformations $\left(f, \phi\right) \mapsto \left(p_{f}, p_{\phi}\right)$ that are part of trinatural transformations between $\mathbf{Gray}$-functors from $\mathfrak{B}$ to $\mathfrak{C}$,
		\item Middle-four interchange in the hom-$2$-categories of $\mathfrak{C}$.
	\end{itemize}	
	\noindent We illustrate the proof of the composition axiom for the trimodification $p_{j}: p_{G}.Sj \Rightarrow Tj.p_{F}$ on a composable pair \begin{tikzcd}
		X \arrow[r, "f"] & Y \arrow[r, "g"] & Z
	\end{tikzcd} in $\mathfrak{A}$. We work under the assumption that $p$ and $j$ are both semi-strict. A reader having trouble reconstructing the pasting diagrams that need to be proved to be equal in this Subsection may find it easier to adapt this proof to one without semi-strictness assumptions on $p$ and $j$ by inserting compositors in appropriate places.
	\\
	\\
	\noindent Various axioms mentioned in the proof of this property in Proposition \ref{interchanging tritransformations} are not needed in this simpler setting. These include the associativity axiom for $p: S \rightarrow T$, and naturality of $p$ on the compositors of $j$ on composable pairs \begin{tikzcd}
			X \arrow[r, "f"] & Y \arrow[r, "g"] & Z
		\end{tikzcd} in $\mathfrak{A}$. Instead, the two remaining steps in the proof are the two whiskering axioms for $p$, as described in the proof of Proposition \ref{interchanging tritransformations} with reference to the axioms for trinatural transformations described in Definition \ref{definition trinatural transformation}.
	
	\begin{diagram}\label{Diagram no interchange used in proof of composition axiom for interchanger of trinatural transformations}
		Proof of composition axiom for the trimodification $p_{j}: p_{G}.Sj \Rightarrow Tj.p_{F}$ when $p: S \rightarrow T$ and $j: F \rightarrow G$ are both semi-strict trinatural transformations.
		\\
		\\
		
		\noindent \begin{tikzcd}[column sep = 4, row sep = 20, font=\fontsize{6.5}{5}]
			&&&&p_{GZ}.Sj_{Z}.SFg.SFf
			\arrow[lldd, "1.Sj_{g}.1"']
			\arrow[rr, bend left = 40, "p_{j_{X}}.1.1"]
			&{}\arrow[dddd,red,  Rightarrow, shorten = 30, "p_{j_{g}}"]
			&Tj_{Z}.p_{FZ}.SFg.SFf
			\arrow[rrdd, "1.p_{Fg}.1"]
			\\
			\\
			&&p_{GZ}.SGg.Sj_{Y}.SFf
			\arrow[rrdd, "p_{Gg}.1.1"]
			\arrow[dddd, Rightarrow, blue, shorten = 30,shift right = 5,  "{\left(p_{Gg}\right)}_{\left(Sj_{f}\right)}"]
			\arrow[lldd, "1.1.Sj_{f}"']
			&&&&&&Tj_{Z}.TFg.p_{FY}.SFf
			\arrow[dddd, Rightarrow,red, shorten = 30,shift right = 5, "{\left(Tj_{g}\right)}_{\left(p_{Ff}\right)}"]
			\arrow[lldd, "Tj_{g}.1.1"']
			\arrow[rrdd, "1.1.p_{Ff}"]
			\\
			\\
			p_{GZ}.SGg.SGf.Sj_{X}
			\arrow[rrdd, "p_{Gg}.1.1"']
			&&&&TGg.p_{GY}.Sj_{Y}.SFf
			\arrow[lldd, "1.1.Sj_{f}"]
			\arrow[rr, "1.p_{j_{Y}}.1"]
			&{}
			\arrow[dddd, Rightarrow, blue, shorten = 30, "p_{j_{f}}"]
			&TGg.Tj_{Y}.p_{FY}.SFf
			\arrow[rrdd, "1.1.p_{Ff}"']
			&&&&Tj_{Z}.TFg.TFf.p_{FX}
			\arrow[lldd, "Tj_{g}.1.1"]
			\\
			\\
			&&TGg.p_{GY}.SGf.Sj_{X}
			\arrow[rrdd, "1.p_{Gf}.1"']
			&&&&&&TGg.Tj_{Y}.TFf.p_{FX}
			\arrow[lldd, "1.Tj_{f}.1"]
			\\
			\\
			&&&&TGg.TGf.p_{GX}.Sj_{X}
			\arrow[rr, bend right = 40, "1.1.p_{j_{X}}"']
			&{}&TGg.TGf.Tj_{X}.p_{FX}
		\end{tikzcd}

	\end{diagram}
	 In particular, considering the data depicted in Diagram \ref{Diagram no interchange used in proof of composition axiom for interchanger of trinatural transformations} above,
	
	\begin{itemize}
		\item the cells in \textcolor{red}{red} correspond to the $3$-cell component of $p$ on the whiskering depicted below left,
		\item the cells in \textcolor{blue}{blue} correspond to the $3$-cell component of $p$ on the whiskering depicted below right.
	\end{itemize} 
	
	$$\begin{tikzcd}[font=\fontsize{9}{6}]
		FX \arrow[rr, "Ff"] && FY \arrow[rr, bend right=40, "Gg. j_{Y}"' name = A]
		\arrow[rr, bend left = 40, "j_{Z}. Fg"name=B] 
		\arrow[from=B, to=A, Rightarrow, shorten = 5, "j_{g}"]&& GZ&{}
	\end{tikzcd}\begin{tikzcd}[font=\fontsize{9}{6}]
		FX \arrow[rr, bend right=40, "Gf. j_{X}"' name = A]
		\arrow[rr, bend left = 40, "j_{Y}. Ff"name=B] 
		\arrow[from=B, to=A, Rightarrow, shorten =5, "j_{f}"]&& GY \arrow[rr, "Gg"] && GZ
	\end{tikzcd}$$
\end{remark}

\noindent Thus far we have described how trinatural transformations $p: S \rightarrow T: \mathfrak{B} \rightarrow \mathfrak{C}$ interact with data in the $\mathbf{Gray}$-category $\mathbf{Tricat}_{s}\left(\mathfrak{A}, \mathfrak{B}\right)$ to produce data in the $\mathbf{Gray}$-category $\mathbf{Tricat}_\text{s}\left(\mathfrak{A}, \mathfrak{C}\right)$. We now collect this information in the form of a trinatural transformation in Proposition \ref{Tricat A, p trinatural}. The assignment of $p$ to such a trinatural transformation is then be shown to be behaviour on morphisms of a functor in Proposition \ref{functoriality of A, p}. Finally, this construction is restricted to semi-strictly generated trinatural transformations to give the underlying functors of the internal whiskering operators in our closed structure in Corollary \ref{underlying functor of internal whiskering operator in closed structure}.

\begin{proposition}\label{Tricat A, p trinatural}
	Let $p: S \rightarrow T: \mathfrak{B} \rightarrow \mathfrak{C}$ and $\mathfrak{A}$ be as in Equation \ref{pair of trinatural transformations between Gray-functors}.
	
	\begin{enumerate}
		\item There is a trinatural transformation $\mathbf{Tricat}_{s}\left(\mathfrak{A}, p\right): \mathbf{Tricat}_{s}\left(\mathfrak{A}, S\right) \Rightarrow \mathbf{Tricat}_{s}\left(\mathfrak{A}, T\right): \mathbf{Tricat}_{s}\left(\mathfrak{A}, \mathfrak{B}\right) \rightarrow \mathbf{Tricat}_{s}\left(\mathfrak{A}, \mathfrak{C}\right)$ whose
		
		\begin{itemize}
			\item $1$-cell component at an object $F \in \mathbf{Tricat}_\text{s}\left(\mathfrak{A}, \mathfrak{B}\right)$ is given by the trinatural transformation $p_F : SF \rightarrow TF$.
			\item component adjoint equivalence at a morphism $j: F \rightarrow G$ in the $\mathbf{Gray}$-category $\mathbf{Tricat}_\text{s}\left(\mathfrak{A}, \mathfrak{B}\right)$ is given by the trimodification adjoint equivalence $p_j \dashv p_{j}^{*}$ of Proposition \ref{interchanging tritransformations}.
			\item $3$-cell component at a $2$-cell $\alpha: j \Rightarrow k: F \rightarrow G$ in the $\mathbf{Gray}$-category $\mathbf{Tricat}_\text{s}\left(\mathfrak{A}, \mathfrak{B}\right)$ is given by the invertible perturbation $p_\alpha$ of Proposition \ref{interchanging trimodification with tritransformation}.
			\item Unitor on an object $F \in \mathbf{Tricat}_\text{s}\left(\mathfrak{A}, \mathfrak{B}\right)$ is given by the invertible perturbation $p^F$ of Proposition \ref{unitor of whiskering with p}.
			\item Compositor on a composable pair \begin{tikzcd}
				F \arrow[r, "j"] & G \arrow[r, "j'"] & H
			\end{tikzcd} in the $\mathbf{Gray}$-category $\mathbf{Tricat}_\text{s}\left(\mathfrak{A}, \mathfrak{B}\right)$ is given by the invertible perturbation $p_{j', j}$ of Proposition \ref{compositor of whiskering with p}.
		\end{itemize}
		\item If $p$ is semi-strict then so is $\mathbf{Tricat}_{s}\left(\mathfrak{A}, p\right)$
		\item If $p$ is semi-strict then $\mathbf{Tricat}_{s}\left(\mathfrak{A}, p\right)$ has semi-strict $1$-cell components.
		\item If $p$ is semi-strict then it restricts to a semi-strict trinatural transformation from the $\mathbf{Gray}$-functor $[\mathfrak{A}, S]: [\mathfrak{A}, \mathfrak{B}]\rightarrow [\mathfrak{A}, \mathfrak{C}]$ to the $\mathbf{Gray}$-functor $[\mathfrak{A}, T]: [\mathfrak{A}, \mathfrak{B}]\rightarrow [\mathfrak{A}, \mathfrak{C}]$.
	\end{enumerate} 
\end{proposition}

\begin{proof}
	For part (1), each axiom for $\mathbf{Tricat}_{s}\left(\mathfrak{A}, p\right)$ follows from the corresponding axiom for $p$ component-wise, as we will now describe.
	
	\begin{itemize}
		\item Local pseudonaturality conditions for the assignment $\left(j, \alpha\right) \mapsto \left(p_j, p_{\alpha}\right)$ follow from local pseudonaturality conditions for $p$ on each $\left(j_{X}, \alpha_{X}\right)$. The same is true for $\left(j, \alpha\right) \mapsto \left(p_{j}, p_{\alpha}\right)$.
		\item The whiskering conditions for $\mathbf{Tricat}_{s}\left(\mathfrak{A}, p\right)$ on $\alpha: j \Rightarrow k: F \rightarrow G$ and $j': G \rightarrow H$ follow from the whiskering condition for $p$ on each $\alpha_{X}: j_{X} \Rightarrow k_{X}: FX \rightarrow GX$ and $j'_{X}: GX \rightarrow HX$. The same is true for left whiskering.
		\item The left and right unit laws for  $\mathbf{Tricat}_{s}\left(\mathfrak{A}, p\right)$ on a morphism $j: F \rightarrow G$ hold by the same laws for $p$ on each component $j_{X}: FX \rightarrow GX$.
		\item The associativity law on a composable triple of trinatural transformations depicted below holds by the same law for $p$ on the composable triple of components in $\mathfrak{B}$, on each $X \in \mathfrak{A}$. $$\begin{tikzcd}
			F \arrow[r, "j"] & G \arrow[r, "k"] & H \arrow[r, "l"] & K
		\end{tikzcd}$$
	\end{itemize}
	Part (2) is clear from the constructions of the compositor and unitor of $\mathbf{Tricat}_{s}\left(\mathfrak{A}, p\right)$ in terms of those same data from $p$. Part (3) is just the fact that if $p$ is semi-strict then so is $p_F$, which has already been established as part of Proposition \ref{coherent hom on Gray Cat} part (1). Part (4) follows from part (3), since a semi-strict trinatural transformation is in particular semi-strictly decomposable.
\end{proof}

\noindent Next, we show in Proposition \ref{functoriality of A, p} that the construction we have defined so far involving arbitrary trinatural transformations is functorial between underlying categories. Once we know functoriality, results in Proposition \ref{Tricat A, p trinatural} parts (2), (3) and (4) will extend from semi-strict trinatural transformations to semi-strictly decomposable ones. This will then be used to prove underlying functoriality of the internal whiskering operators in the closed structure.

\begin{proposition}\label{functoriality of A, p}
	Consider the underlying categories $\mathbf{Tricat}_{s}\left(\mathfrak{B}, \mathfrak{C}\right)_{1}$ and $\mathbf{Tricat}_{s}\left(\mathfrak{K}, \mathfrak{L}\right)_{1}$ where $\mathfrak{K}$ is the $\mathbf{Gray}$-category $\mathbf{Tricat}_{s}\left(\mathfrak{A}, \mathfrak{B}\right)$, and $\mathfrak{L}$ is the $\mathbf{Gray}$-category $\mathbf{Tricat}_{s}\left(\mathfrak{A}, \mathfrak{C}\right)$. There is a functor $\mathbf{Tricat}_{s}\left(\mathfrak{B}, \mathfrak{C}\right)_{1}\rightarrow \mathbf{Tricat}_{s}\left(\mathfrak{K},\mathfrak{L}\right)_{1}$ which 
	
	\begin{itemize}
		\item on objects, sends a $\mathbf{Gray}$-functor $S: \mathfrak{B} \rightarrow \mathfrak{C} $ to the $\mathbf{Gray}$-functor $\mathbf{Tricat}_{s}\left(\mathfrak{A}, S\right): \mathbf{Tricat}_\text{s}\left(\mathfrak{A}, \mathfrak{B}\right) \rightarrow \mathbf{Tricat}_\text{s}\left(\mathfrak{A}, \mathfrak{C}\right)$,
		\item on morphisms, sends a trinatural transformation $p: S \rightarrow T$ to the trinatural transformation described in Proposition \ref{Tricat A, p trinatural}.
	\end{itemize}
\end{proposition}

\begin{proof}
	Preservation of identities follows from Proposition \ref{interchanging tritransformations} part (3), Proposition \ref{interchanging trimodification with tritransformation} part (2) and an inspection of Proposition \ref{unitor of whiskering with p} and Proposition \ref{compositor of whiskering with p}. To show preservation of composition we need to show that each component of $\mathbf{Tricat}_{s}\left(\mathfrak{A}, p\right) \circ \mathbf{Tricat}_{s}\left(\mathfrak{A}, r\right)$ and $\mathbf{Tricat}_{s}\left(\mathfrak{A}, p\circ r\right)$ coincide for any pair of composable trinatural transformations \begin{tikzcd}
		R \arrow[r, "r"] & S \arrow[r, "p"] &T
	\end{tikzcd}.
\\
\\
\noindent The following can be deduced by unwinding definitions, keeping Remark \ref{Gray Functor Gray category with weak higher cells} in mind.

	\begin{itemize}
		\item For both trinatural transformations $\mathbf{Tricat}_{s}\left(\mathfrak{A}, p\right) \circ \mathbf{Tricat}_{s}\left(\mathfrak{A}, r\right)$ and $\mathbf{Tricat}_{s}\left(\mathfrak{A}, p\circ r\right)$, the $1$-cell component at a $\mathbf{Gray}$-functor $F: \mathfrak{A} \rightarrow \mathfrak{B}$ is given by the trinatural transformation \begin{tikzcd}
			RF \arrow[rr, "r_F"] && SF \arrow[rr, "p_F"] &&TF 
		\end{tikzcd}
		\item For both trinatural transformations, the $2$-cell component at a trinatural transformation $j: F \rightarrow G$ is given by the trimodification whose $2$-cell component at $X$ is displayed below. 
		
		$$\begin{tikzcd}
			p_{GX}.r_{GX}.Rj_{X} \arrow[r, "1.r_{j_{X}}"] & p_{GX}.Sj_{X}.r_{FX} \arrow[r, "p_{j_{X}} .1"] & Tj_{X}.p_{FX}.r_{FX}
		\end{tikzcd}$$
	
		\noindent We describe why for a morphism $f: X \rightarrow Y$ in $\mathfrak{A}$, both trimodifications have the same $3$-cell component, namely the one depicted in Diagram \ref{functoriality pasting diagram}. 
		\begin{itemize}
			\item For the trimodification $\mathbf{Tricat}_{s}\left(\mathfrak{A}, p\circ r\right)_{j}$, we must first compute the $3$-cell component of the trinatural transformation $p \circ r$ at the $2$-cell $j_{f}$. This is given by the pasting of $p_{j_{f}}.1$ and $1.r_{j_{f}}$ along $1.Sj_{X}.1$. Then we must attach the compositors of $p \circ r$ at the pair of composable pairs \begin{tikzcd}
				FX \arrow[r, "j_{X}"] & GX \arrow[r, "Gf"] & GY
			\end{tikzcd} and \begin{tikzcd}
				FX \arrow[r, "Ff"] & FY \arrow[r, "j_{Y}"] & GY
			\end{tikzcd} in $\mathfrak{B}$. Using the description of the compositor of the composite trinatural transformation $p \circ r$ given in Remark \ref{Gray Functor Gray category with weak higher cells}, we see that this is the pasting depicted in Diagram \ref{functoriality pasting diagram}.
			\item For the trimodification ${\left(\mathbf{Tricat}_{s}\left(\mathfrak{A}, p\right)\circ \mathbf{Tricat}_{s}\left(\mathfrak{A}, r\right)\right)}_{j}$, the $3$-cell component at $f$ will be given by that of the following trimodification.
			
			$$\begin{tikzcd}
				RF \arrow[rr, "r_{F}"]
				\arrow[dd, "Rj"'] &{}\arrow[dd, Leftarrow, shorten = 10, "r_{j}"]
				& SF \arrow[rr, "p_{F}"]\arrow[dd, "Sj"] &{}\arrow[dd, Leftarrow, shorten = 10, "p_{j}"]
				& TF\arrow[dd, "Tj"]
				\\
				\\
				RG \arrow[rr, "r_{G}"'] &{}& SG \arrow[rr, "p_{G}"'] &{}& TG
			\end{tikzcd}$$
			
			\noindent This trimodification is the composite of the whiskering of $r_{j}$ with $p_{G}$ with the whiskering of $p_{j}$ with $r_{F}$. But \begin{itemize}
				\item The trimodification $r_{j}$ has $3$-cell component on $f$ given by pasting together $r_{j_{f}}$ with the compositors $r_{Gf, j_{X}}$ and $r_{j_{Y}, Ff}^{-1}$. Hence by part (4) of Remark \ref{Gray Functor Gray category with weak higher cells}, its whiskering with the trinatural transformation $p_{G}$ will have $3$-cell component on $f$ given by this $3$-cell pasted with the interchanger ${\left(p_{Gf}\right)}_{\left(r_{j_{X}}\right)}$.
				\item Similarly, the trimodification $p_{j}$ has $3$-cell component on $f$ given by pasting together $p_{j_{f}}$ with the compositors $p_{Gf, j_{X}}$ and $p_{j_{Y}, Ff}^{-1}$. Hence by Remark \ref{Gray Functor Gray category with weak higher cells} part (4), its whiskering with the trinatural transformation $r_{F}$ will have $3$-cell component on $f$ given by this $3$-cell pasted with the interchanger ${\left(p_{j_{Y}}\right)}_{\left(r_{Ff}\right)}^{-1}$.
			\end{itemize}
			The desired $3$-cell will hence be the pasting of these, as depicted in Diagram \ref{functoriality pasting diagram}.
		\end{itemize}
		Although we have described how to construct Diagram \ref{functoriality pasting diagram} in two different ways, these are related by middle four interchange in a genuine $2$-category. Indeed, the $3$-cell components of $\mathbf{Tricat}_{s}\left(\mathfrak{A}, p\circ r\right)_{j}$ and ${\left(\mathbf{Tricat}_{s}\left(\mathfrak{A}, p\right) \circ \mathbf{Tricat}_{s}\left(\mathfrak{A}, r\right)\right)}_{j}$ on $f$ are equal by the pasting lemma \cite{Power 2 cat pasting} in the $2$-category $\mathfrak{B}\left(RFX, TGY\right)$.
	\end{itemize}

	\begin{diagram}\label{functoriality pasting diagram}
		$3$-cell component of the trimodifications ${\left(\mathbf{Tricat}_{s}\left(\mathfrak{A}, p\right) \circ \mathbf{Tricat}_{s}\left(\mathfrak{A}, r\right)\right)}_{j}$ and ${\mathbf{Tricat}_{s}\left(\mathfrak{A}, pr\right)}_{j}$ at the morphism $f \in \mathfrak{A}\left(X, Y\right)$.
		\\
		\\
		
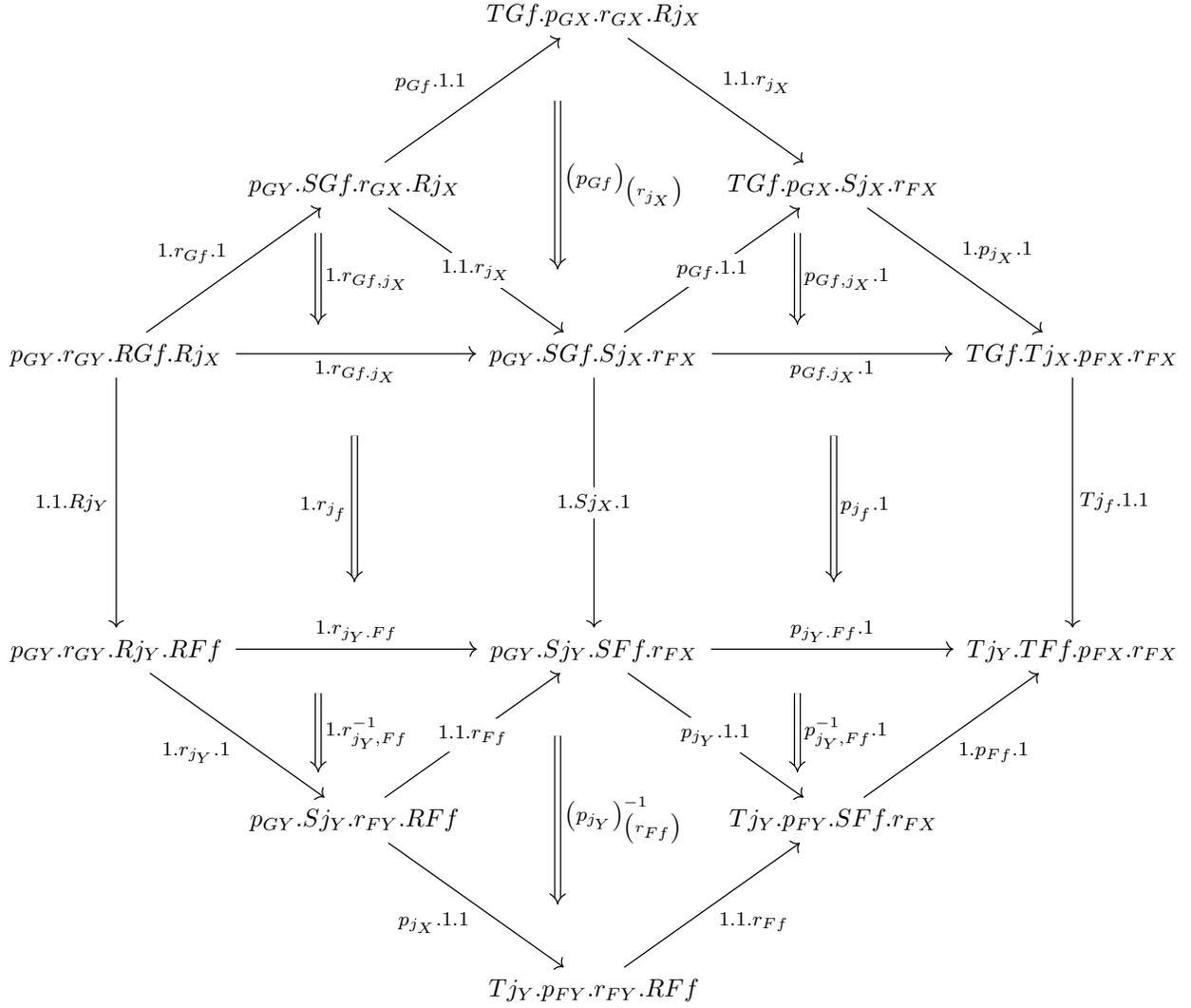
\noindent\begin{tikzcd}[column sep= 0.1, row sep = 25, font=\fontsize{10}{7}]
			&&&&TGf.p_{GX}.r_{GX}.Rj_{X}
			\arrow[dddd, Rightarrow, shorten = 25, "{\left(p_{Gf}\right)}_{\left(r_{j_{X}}\right)}", shift right = 5]
			\arrow[rrdd, "1.1.{r}_{j_{X}}"]
			\\
			\\
			&&p_{GY}.SGf.r_{GX}.Rj_{X}
			\arrow[dd, Rightarrow, shorten = 10, "1.r_{Gf{,}j_{X}}", shift right = 5]
			\arrow[rruu, "p_{Gf}.1.1"]
			\arrow[rrdd, "1.1.r_{j_{X}}" description]
			&&&&TGf.p_{GX}.Sj_{X}.r_{FX}
			\arrow[dd, Rightarrow, shorten = 10, "p_{Gf{,}j_{X}}.1", shift right = 5]
			\arrow[rrdd, "1.p_{j_{X}}.1"]
			\\
			\\
			p_{GY}.r_{GY}.RGf.Rj_{X}
			\arrow[dddd, "1.1.Rj_{Y}"']
			\arrow[rrrr, "1.r_{Gf.j_{X}}"']
			\arrow[rruu, "1.r_{Gf}.1"]
			&&{}\arrow[dddd, "1.r_{j_{f}}"', shorten = 25, Rightarrow]
			&&p_{GY}.SGf.Sj_{X}.r_{FX}
			\arrow[dddd, "1.Sj_{X}.1" description]
			\arrow[uurr, "p_{Gf}.1.1" description]
			\arrow[rrrr, "p_{Gf.j_{X}}.1"']
			&&{}\arrow[dddd, "p_{j_{f}}.1", shorten = 25, Rightarrow]
			&&TGf.Tj_{X}.p_{FX}.r_{FX}
			\arrow[dddd, "Tj_{f}.1.1"]
			\\
			\\
			\\
			\\
			p_{GY}.r_{GY}.Rj_{Y}.RFf
			\arrow[rrrr, "1.r_{j_{Y}.Ff}"]
			\arrow[rrdd, "1.r_{j_{Y}}.1"']
			&&{}\arrow[dd, Rightarrow, shorten = 10, shift right = 5, "1.r_{j_{Y}{,}Ff}^{-1}"]
			&&p_{GY}.Sj_{Y}.SFf.r_{FX}
			\arrow[ddrr, "p_{j_{Y}}.1.1" description]
			\arrow[rrrr, "p_{j_{Y}.Ff}.1"]
			\arrow[dddd, Rightarrow, shorten = 25, "{\left(p_{j_{Y}}\right)}_{\left(r_{Ff}\right)}^{-1}", shift right = 5]
			&&{}\arrow[dd, Rightarrow, shorten =10, "p_{j_{Y}{,}Ff}^{-1}.1", shift right = 5]
			&&Tj_{Y}.TFf.p_{FX}.r_{FX}
			\\
			\\
			&&p_{GY}.Sj_{Y}.r_{FY}.RFf
			\arrow[rrdd, "p_{j_{X}}.1.1"']
			\arrow[rruu, "1.1.r_{Ff}" description]
			&&&&Tj_{Y}.p_{FY}.SFf.r_{FX}
			\arrow[rruu, "1.p_{Ff}.1"']
			\\
			\\
			&&&&Tj_{Y}.p_{FY}.r_{FY}.RFf
			\arrow[rruu, "1.1.r_{Ff}"']
		\end{tikzcd}
	\end{diagram} 
	\begin{itemize}
		
		\item For both trinatural transformations, the $3$-cell component at a trimodification $\alpha: j \rightarrow k$ is given by the perturbation
		\\
		\begin{tikzcd}
			p_{G}.r_{G}.Rj
			\arrow[rr, "1.r_{j}"] \arrow[dd, "1.1.R\alpha"']
			&{}\arrow[dd, Rightarrow, shorten = 10, "1.r_{\alpha}"]
			&
			p_{G}.Sj.r_{F}
			\arrow[rr, "p_{j} .1"]
			\arrow[dd, "1.S\alpha.1"] &{}\arrow[dd, Rightarrow, shorten = 10, "p_{\alpha}.1"]
			& Tj.p_{F}.r_{F}
			\arrow[dd, "T\alpha.1.1"]
			\\
			\\
			p_G.r_G.Rk
			\arrow[rr, "1.r_{k}"']
			&{}&
			p_{G}.Sk.r_{F}
			\arrow[rr, "p_{k}.1"']
			&{}& 
			Tk.p_{F}.r_{F}
		\end{tikzcd}
		\item For each trinatural transformation, the unitor at a $\mathbf{Gray}$-functor $F$ is given by the horizontal composite depicted below.
		
		$$\begin{tikzcd}
			RF
			\arrow[rr, bend left, "1_{r_{F}}" name = A]
			\arrow[rr, bend right, "r_{1_{F}}"' name = B]
			&& SF
			\arrow[rr,bend left, "1_{p_{F}}" name = C]
			\arrow[rr, bend right, "p_{1_{F}}"' name = D]
			&& TF
			\arrow[from = A, to = B, Rightarrow, shorten = 5, "r^F"]
			\arrow[from = C, to = D, Rightarrow, shorten =5, "p^F"]
		\end{tikzcd}$$
		\item Each of the two trinatural transformations have a compositor determined by any composable pair of trinatural transformations \begin{tikzcd}
			F \arrow[r, "j"] & G \arrow[r, "k"] &H
		\end{tikzcd}. They are both given by the following pasting diagram of perturbations.
	\end{itemize}
	$$\begin{tikzcd}[column sep= 15, font=\fontsize{10}{7}]
		&&&&Tk.p_{G}.r_{G}.Rj
		\arrow[dddd, Rightarrow, shorten = 25, "{\left(p_{k}\right)}_{\left(r_{j}\right)}"]
		\arrow[rrdd, "1.1.{r}_{j}"]
		\\
		\\
		&&p_{H}.Sk.r_{G}.Rj
		\arrow[dd, Rightarrow, shorten = 10, "1.r_{k{,}j}"]
		\arrow[rruu, "p_{k}.1.1"]
		\arrow[rrdd, "1.1.r_{j}"]
		&&&&Tk.p_{G}.Sj.r_{F}
		\arrow[dd, Rightarrow, shorten = 10, "p_{k{,}j}.1"]
		\arrow[rrdd, "1.p_{j}.1"]
		\\
		\\
		p_{H}.r_{H}.Rk.Rj
		\arrow[rrrr, "1.r_{kj}"']
		\arrow[rruu, "1.r_{k}.1"]
		&&{}
		&&p_{H}.Sk.Sj.r_{F}
		\arrow[uurr, "p_{k}.1.1"]
		\arrow[rrrr, "p_{kj}.1"']
		&&{}
		&&Tk.Tj.p_{F}.r_{F}
	\end{tikzcd}$$
	\begin{itemize}
		\item Each of the two trinatural transformations has a component adjoint equivalence. The argument that their units and counits respectively coincide is similar to the argument above for compositors.
	\end{itemize}
\end{proof}

\noindent We are finally in a position to restrict the functor to one between semi-strictly generated homs of $\mathbf{Gray}$-categories, and establish the underlying functors $[\mathfrak{A}, -]_{1}: [\mathfrak{B}, \mathfrak{C}]_{1} \rightarrow [[\mathfrak{A}, \mathfrak{B}], [\mathfrak{A}, \mathfrak{C}]]_{1}$ of the internal whiskering operators in the closed structure.

\begin{corollary}\label{strengthening results about semi strictness to results about coherence}
	Let $p: S \Rightarrow T: \mathfrak{B} \rightarrow \mathfrak{C}$ be semi-strictly decomposable. Then the trinatural transformation $\mathbf{Tricat}_{s}\left(\mathfrak{A}, p\right)$ of Proposition \ref{Tricat A, p trinatural}
	\begin{enumerate}
		\item has semi-strictly decomposable $1$-cell components.
		\item is itself semi-strictly decomposable.
		\item restricts to a semi-strictly decomposable trinatural transformation  $[\mathfrak{A}, p]: [\mathfrak{A}, S] \rightarrow [\mathfrak{A}, T]: [\mathfrak{A}, \mathfrak{B}] \rightarrow [\mathfrak{A}, \mathfrak{C}]$. 
	\end{enumerate}
\end{corollary}

\begin{proof}
	For part (1), notice that the component of $\mathbf{Tricat}_{s}\left(\mathfrak{A}, p\right): \mathbf{Tricat}_{s}\left(\mathfrak{A}, S\right) \rightarrow  \mathbf{Tricat}_{s}\left(\mathfrak{A}, T\right)$ on $F$ will coincide with the image of $p: S \rightarrow T$ under the $\mathbf{Gray}$-functor $[F, \mathfrak{C}]: [\mathfrak{B}, \mathfrak{C}] \rightarrow [\mathfrak{A}, \mathfrak{C}]$ from Proposition \ref{coherent hom on Gray Cat} part (1). It is part of that result that this trinatural transformation will be semi-strictly decomposable. For part (2), let $p = p_{n}\circ ...\circ p_{1}$ be a factorisation of $p$ into semi-strict trinatural transformations and let $ 1\leq i \leq n$. Each $\mathbf{Tricat}_{s}\left(\mathfrak{A}, p_{i}\right)$ has semi-strict components by Proposition \ref{Tricat A, p trinatural} part (3). Hence by functoriality of $\mathbf{Tricat}_{s}\left(\mathfrak{A}, -\right)$ from Proposition \ref{functoriality of A, p}, this means that $\mathbf{Tricat}_{s}\left(\mathfrak{A}, p\right)$ also be semi-strictly decomposable. Part (3) similarly follows from Proposition \ref{Tricat A, p trinatural} part (4) using functoriality.
\end{proof}

\begin{corollary}\label{underlying functor of internal whiskering operator in closed structure}
	The functor between underlying categories considered in Proposition \ref{functoriality of A, p} restricts to a functor between underlying categories $[\mathfrak{A}, -]_{1}: [\mathfrak{B}, \mathfrak{C}]_{1} \rightarrow [[\mathfrak{A}, \mathfrak{B}], [\mathfrak{A}, \mathfrak{C}]]_{1}$.
\end{corollary}

\begin{proof}
	This follows from Proposition \ref{functoriality of A, p} and Corollary \ref{strengthening results about semi strictness to results about coherence}.
\end{proof}

\subsection{$\mathbf{Gray}$-functoriality of the internal whiskering operators}\label{Subsection Gray functoriality of internal whiskering operators}

\noindent We now turn to giving the functors $[\mathfrak{A}, -]_{1}: [\mathfrak{B}, \mathfrak{C}]_{1} \rightarrow [[\mathfrak{A}, \mathfrak{B}], [\mathfrak{A}, \mathfrak{C}]]_{1}$ of Corollary \ref{underlying functor of internal whiskering operator in closed structure} an enrichment over $\mathbf{Gray}$. We first describe a $3$-cell $\sigma_{j}$ in $\mathbf{Tricat}_{s}\left(\mathfrak{A}, \mathfrak{C}\right)$ which mediates interchange between a $1$-cell $j: F \rightarrow G$ in $\mathbf{Tricat}_{s}\left(\mathfrak{A}, \mathfrak{B}\right)$ and a $2$-cell $\sigma: p \Rightarrow q$ in $\mathbf{Tricat}_{s}\left(\mathfrak{B}, \mathfrak{C}\right)$. This description is given in Proposition \ref{interchanging tritransformation with trimodification}. The perturbation $\sigma_{j}$ will feature as the $3$-cell component of the trimodification $\mathbf{Tricat}_{s}\left(\mathfrak{A}, \sigma\right): \mathbf{Tricat}_{s}\left(\mathfrak{A}, p\right) \Rightarrow \mathbf{Tricat}_{s}\left(\mathfrak{A}, q\right)$ on the $1$-cell $j: F \rightarrow G$, in Proposition \ref{Tricat A sigma trimodification}.

\begin{proposition}\label{interchanging tritransformation with trimodification}
	Let $p: S \Rightarrow T: \mathfrak{B} \rightarrow \mathfrak{C}$ and $j: F \Rightarrow G: \mathfrak{A} \rightarrow \mathfrak{B}$ be as in Equation \ref{pair of trinatural transformations between Gray-functors} and let $\sigma: p \Rrightarrow q: S \Rightarrow T$ be a trimodification. 
	
	\begin{enumerate}
		\item There is a invertible perturbation $\sigma_{j}$ as depicted below left, whose component on $X \in \mathfrak{A}$ is given by the $3$-cell in $\mathfrak{C}$ depicted below right.
		\item The perturbation $\sigma_{j}$ is the identity if $j$ is a pseudo-icon or if $\sigma$ is strict. In particular it is the identity if either $\sigma$ or $j$ is the identity.
	\end{enumerate}

	\begin{tikzcd}[column sep = 18, row sep = 20]
		SF
		\arrow[rr,bend left = 45, "q_{F}"]
		\arrow[dd,"Sj"']
		&
		{}
		\arrow[d,Leftarrow, shorten = 5,"q_{j}"]
		&
		TF
		\arrow[dd,"Tj"]
		\\
		&{}&&\cong^{\sigma_{j}}
		\\
		SG
		\arrow[rr,bend left = 45,"q_{G}" {name = C}]
		\arrow[rr, bend right = 45, "p_{G}"' {name = D}]
		&
		{}
		&
		GY
		\arrow[from =D, to =C, Rightarrow, shorten = 10, shift right = 5, "\sigma_{G}"]
	\end{tikzcd}	\begin{tikzcd}[column sep = 18, row sep = 20]
		SF
		\arrow[rr,bend left = 45, "q_{F}" {name = A}]
		\arrow[rr, bend right = 45, "p_{F}"' {name = B}]
		\arrow[dd,"Sj"']
		&
		&
		TF
		\arrow[dd,"Tj"]
		\\
		&{}
		\arrow[d,Leftarrow, shorten = 5,"p_{j}"]
		\\
		SG\
		\arrow[rr, bend right = 45, "p_{G}"']
		&
		{}
		&
		TG &{}
		\arrow[from =A, to =B, Leftarrow, shorten = 10, shift right = 5, "\sigma_{F}"]
	\end{tikzcd} \begin{tikzcd}
		p_{GX}.Sj_{X} \arrow[rr, "p_{j_{X}}"] \arrow[dd, "\sigma_{GX}.1"']
		&{}\arrow[dd, Rightarrow, shorten = 10, "\sigma_{j_{X}}"]
		&Tj_{X}.p_{FX}\arrow[dd, "1.\sigma_{FX}"]
		\\
		\\
		q_{GX}.Sj_{X}\arrow[rr, "q_{j_{X}}"'] &{}& Tj_{X}.q_{FX}
	\end{tikzcd}
\end{proposition}

\begin{proof}
	This is A.71 of Proposition A.12 in \cite{Buhne PhD}. The proof follows the pattern described in Remark \ref{proof strategy pasting diagram chases}. The perturbation axiom is a condition on a morphism $f: X \rightarrow Y$ in $\mathfrak{A}$. Now $j$ maps $f$ to the $2$-cell $j_f: j_{Y}.Ff \Rightarrow Gf.j_{X}$ in $\mathfrak{B}$. As we detail below, the proof hence uses all axioms for the trimodification $\sigma: p \Rightarrow q$ on these data in $\mathfrak{B}$.
	\\
	\\
	\noindent Start with the the right hand side of the perturbation axiom in Definition \ref{Definition perturbation}, which is given by the pasting depicted below.
	
	$$\begin{tikzcd}[column sep = 40, row sep = 15, font=\fontsize{9}{6}]
		{p_{GY}}.{Sj_Y}.{SFf}
		\arrow[rr,"{{\sigma}_{GY}}.1.1"]
		\arrow[dd,"1.{Sj_f}"']
		&{}&
		{q_{GY}}.{Sj_Y}.{SFf}
		\arrow[rr,"q_{j_Y}.1"]
		\arrow[rrdd, "q_{j_{Y}.Ff}" {description, name=A}]
		\arrow[dd, "1.Sj_{f}" description]
		&{}&
		{Tj_Y}.{q_{FY}}.{SFf}
		\arrow[dd,"1.{q_{Ff}}"]
		\\
		&{}
		\\
		{p_{GY}}.{SGf}.{Sj_X}
		\arrow[rr, "\sigma_{GY}.1.1"]
		\arrow[dd,"{p_{Gf}}.1"']
		&{}\arrow[uu, Rightarrow, shorten = 15, "{\left(\sigma_{GY}\right)}_{\left(Sj_{f}\right)}"]&
		{Tj_Y}.{TFf}.{p_{FX}}
		\arrow[dd, "q_{Gf}.1" description]
		\arrow[ddrr,"q_{Gf.j_{X}}"name=B]
		\arrow[from=B, to=A, Rightarrow, shorten = 10, "q_{j_{f}}"]
		&{}&
		{Tj_Y}.{TFf}.{q_{FX}}
		\arrow[dd,"{Tj_f}.1"]
		\arrow[uu, shorten = 15, shift left = 5, Rightarrow, "q_{j_{Y}{,}Ff}", near end]
		\\
		&{}
		&&{}\arrow[d, Leftarrow, shorten = 5, "q_{Gf{,}j_{X}}", shift right = 10]
		\\
		{TGf}.{p_{GX}}.{Sj_X}
		\arrow[rr, "1.\sigma_{GX}.1"]
		\arrow[rrdd,"1.{p_{j_X}}"']
		&{}\arrow[uu, Rightarrow, shorten = 15,"\sigma_{Gf}.1"]
		&TGf.q_{GX}.Sj_{X}\arrow[rr, "1.q_{j_{X}}"']
		\arrow[dd, Leftarrow, shorten = 10, "1.\sigma_{j_{X}}"]
		&{}&
		{TGf}.{Tj_X}.{q_{FX}}
		\\
		\\
		&&
		{TGf}.{Tj_X}.{p_{FX}}
		\arrow[rruu,"1.1.{{\sigma}_{FX}}"']
	\end{tikzcd}$$
	
	\noindent Apply the composition axiom for $\sigma$ on the composable pair \begin{tikzcd}
		FX \arrow[r, "Ff"] & FY \arrow[r, "j_{Y}"] & GY
	\end{tikzcd}
	
	$$\begin{tikzcd}[column sep = 40,row sep=15, font=\fontsize{9}{6}]
		{p_{GY}}.{Sj_Y}.{SFf}
		\arrow[rr,"{{\sigma}_{GY}}.1.1"]
		\arrow[dd,"1.{Sj_f}"']
		&{}&
		{q_{GY}}.{Sj_Y}.{SFf}
		\arrow[rr,"q_{j_Y}.1"]
		\arrow[rrdd, "q_{j_{Y},Ff}" {description,name=B}]
		\arrow[dd, "1.Sj_{f}" description]
		&{}&
		{Tj_Y}.{q_{FY}}.{SFf}
		\arrow[dd,"1.{q_{Ff}}"]
		\\
		&{}
		\\
		{p_{GY}}.{SGf}.{Sj_X}
		\arrow[rr, "\sigma_{GY}.1.1"]
		\arrow[rrdd,"p_{{Gf}.{j_X}}" description]
		\arrow[dd,"{p_{Gf}}.1"']
		&{}\arrow[uu, Rightarrow, shorten = 15, "{\left(\sigma_{GY}\right)}_{\left(Sj_{f}\right)}"]&
		{Tj_Y}.{TFf}.{p_{FX}}
		\arrow[ddrr,"q_{Gf.j_{X}}"'{name=A}]
		\arrow[from=A, to=B, Rightarrow, shorten = 20, "q_{j_{f}}"]
		&{}&
		{Tj_Y}.{TFf}.{q_{FX}}
		\arrow[dd,"{Tj_f}.1"]
		\arrow[uu, shorten = 15, shift left = 5, Rightarrow, "q_{j_{Y}{,}Ff}", near end]
		\\
		&{}
		&&
		\\
		{TGf}.{p_{GX}}.{Sj_X}
		\arrow[uu, shorten = 10, shift right = 5, Rightarrow, "p_{{Gf}{,}{j_{X}}}"', near start, red]
		\arrow[rr,"1.{p_{j_X}}"']
		&&{TGf}.{Tj_X}.{p_{FX}}
		\arrow[uu, shorten = 10, Rightarrow, "\sigma_{Gf.j_{X}}", red]
		\arrow[rr,"1.1.{{\sigma}_{FX}}"']
		&{}&
		{TGf}.{Tj_X}.{q_{FX}}
	\end{tikzcd}$$
	
	\noindent Now apply the local modification axiom for $\sigma$ on the $2$-cell $j_{f}$.
	
	$$\begin{tikzcd}[column sep = 40, row sep = 15, font=\fontsize{9}{6}]
		{p_{GY}}.{Sj_Y}.{SFf}
		\arrow[rr,"{{\sigma}_{GY}}.1.1"]
		\arrow[rrdd,"p_{{j_Y}.{Ff}}" description, blue]
		\arrow[dd,"1.{Sj_f}"']
		&{}&
		{q_{GY}}.{Sj_Y}.{SFf}
		\arrow[rr,"q_{j_Y}.1"]
		\arrow[rrdd, "q_{j_{Y}.Ff}" description]
		&{}&
		{Tj_Y}.{q_{FY}}.{SFf}
		\arrow[dd,"1.{q_{Ff}}"]
		\\
		&{}
		\\
		{p_{GY}}.{SGf}.{Sj_X}
		\arrow[rrdd,"p_{{Gf}.{j_X}}" description]
		\arrow[dd,"{p_{Gf}}.1"']
		&&
		\color{blue}{Tj_Y}.{TFf}.{p_{FX}}
		\arrow[rr,"1.1.{{\sigma}_{FX}}", blue]
		\arrow[dd,"{Tj_f}.1"', blue]
		\arrow[uu, Rightarrow, shorten = 10, red, "\sigma_{j_{Y}.Ff}"]
		&{}&
		{Tj_Y}.{TFf}.{q_{FX}}
		\arrow[dd,"{Tj_f}.1"]
		\arrow[uu, shorten = 15, shift left = 5, Rightarrow, "q_{j_{Y}{,}Ff}", near end]
		\\
		&{}\arrow[uu, Rightarrow, shorten = 15, "p_{j_{f}}", red]
		&&
		\\
		{TGf}.{p_{GX}}.{Sj_X}
		\arrow[uu, shorten = 10, shift right = 5, Rightarrow, "p_{{Gf}{,}{j_{X}}}"', near start]
		\arrow[rr,"1.{p_{j_X}}"']
		&&{TGf}.{Tj_X}.{p_{FX}}
		\arrow[rr,"1.1.{{\sigma}_{FX}}"']
		&{}\arrow[uu, Rightarrow, shorten = 15, "{\left(Tj_{f}\right)}_{\left(\sigma_{FX}\right)}", red]&
		{TGf}.{Tj_X}.{q_{FX}}
	\end{tikzcd}$$
	
	\noindent Finally, apply the composition axiom for $\sigma$ on the composable pair \begin{tikzcd}
		FX \arrow[r, "j_{X}"] & GX \arrow[r, "Gf"] & GY
	\end{tikzcd}, and observe that this results in the left hand side of the perturbation axiom in Definition \ref{Definition perturbation}, completing the proof for part (1).
	
	$$\begin{tikzcd}[column sep = 40, row sep = 15, font=\fontsize{9}{6}]
		&&
		{q_{GY}}.{Sj_Y}.{SFf}
		\arrow[rrd,"q_{j_Y}.1"]
		\\
		{p_{GY}}.{Sj_Y}.{SFf}
		\arrow[rru,"{{\sigma}_{GY}}.1.1"]
		\arrow[rr,"p_{j_Y}.1" description,blue]
		\arrow[rrdd,"p_{{j_Y}.{Ff}}" description]
		\arrow[dd,"1.{Sj_f}"']
		&{}&
		\color{blue}{Tj_Y}.{p_{FY}}.{SFf}
		\arrow[rr,blue,"1.{{\sigma}_{FY}}.1" description]
		\arrow[dd,"1.{p_{Ff}}",blue]
		\arrow[u,Rightarrow,shorten=3,red,"{{\sigma}_{j_Y}}.1"]
		&{}&
		{Tj_Y}.{q_{FY}}.{SFf}
		\arrow[dd,"1.{q_{Ff}}"]
		\\
		&{}\arrow[u, Rightarrow, shorten = 4, "p_{j_{Y}{,}Ff}"', red]
		\\
		{p_{GY}}.{SGf}.{Sj_X}
		\arrow[rrdd,"p_{{Gf}.{j_X}}" description]
		\arrow[dd,"{p_{Gf}}.1"']
		&&
		{Tj_Y}.{TFf}.{p_{FX}}
		\arrow[rr,"1.1.{{\sigma}_{FX}}"]
		\arrow[dd,"{Tj_f}.1"']
		&{}\arrow[uu, Rightarrow, shorten = 15, red, "1.\sigma_{Ff}"]&
		{Tj_Y}.{TFf}.{q_{FX}}
		\arrow[dd,"{Tj_f}.1"]
		\\
		&{}\arrow[uu, Rightarrow, shorten = 15, "p_{j_{f}}"]
		&&
		\\
		{TGf}.{p_{GX}}.{Sj_X}
		\arrow[uu, shorten = 10, shift right = 5, Rightarrow, "p_{{Gf}{,}{j_{X}}}"', near start]
		\arrow[rr,"1.{p_{j_X}}"']
		&&{TGf}.{Tj_X}.{p_{FX}}
		\arrow[rr,"1.1.{{\sigma}_{FX}}"']
		&{}\arrow[uu, Rightarrow, shorten = 15, "{\left(Tj_{f}\right)}_{\left(\sigma_{FX}\right)}"]&
		{TGf}.{Tj_X}.{q_{FX}}	
	\end{tikzcd}$$
	
	\noindent For part (2), $\sigma_{j}$ is the identity if $j_{X}$ is by the unit law for $\sigma$, and it is clearly the identity if $\sigma$ has identity $2$-cell components.
\end{proof}

\begin{proposition}\label{Tricat A sigma trimodification}
	Let $\sigma: p \Rrightarrow q: S \Rightarrow T:  \mathfrak{B} \rightarrow \mathfrak{C}$ be a trimodification. \begin{enumerate}
		\item There is a trimodification $\mathbf{Tricat}_{s}\left(\mathfrak{A}, \sigma\right): \mathbf{Tricat}_{s}\left(\mathfrak{A}, p\right) \Rightarrow \mathbf{Tricat}_{s}\left(\mathfrak{A}, q\right)$ whose component on $F: \mathfrak{A} \rightarrow \mathfrak{B}$ is the trimodification $\sigma_{F}$ defined in Proposition A.8 of \cite{Buhne PhD}, and whose component on a trinatural transformation $j: F \rightarrow G$ is the perturbation $\sigma_{j}$ of Proposition \ref{interchanging tritransformation with trimodification} part (2).
		\item When $p$ and $q$ are semi-strictly decomposable, the trimodification $\mathbf{Tricat}_{s}\left(\mathfrak{A}, \sigma\right): \mathbf{Tricat}_{s}\left(\mathfrak{A}, p\right) \Rightarrow \mathbf{Tricat}_{s}\left(\mathfrak{A}, q\right)$ restricts to a trimodification $[\mathfrak{A}, \sigma]: [\mathfrak{A}, p] \Rightarrow [\mathfrak{A}, q]$, where $[\mathfrak{A}, p]$ and $[\mathfrak{A}, q]$ are defined as in Corollary \ref{strengthening results about semi strictness to results about coherence} part (3).
	\end{enumerate}
\end{proposition}

\begin{proof}
	For part (1), the trimodification axioms for data in $\mathbf{Tricat}_{s}\left(\mathfrak{A}, \mathfrak{B}\right)$ follow directly from the corresponding axioms for $\sigma$ on their components in $\mathfrak{B}$. The restriction to semi-strictly decomposable trinatural transformations for part (2) is an immediate specialisation of part (1).
\end{proof}

\begin{proposition}\label{whiskering perturbations is functorial}
	Let $\sigma: p \Rightarrow q$ be as in Proposition \ref{interchanging tritransformation with trimodification} and let $\Omega: \sigma \Rrightarrow \tau$ be a perturbation.
	
	\begin{enumerate}
		\item There is a perturbation $\mathbf{Tricat}_{s}\left(\mathfrak{A}, \Omega\right): \mathbf{Tricat}_{s}\left(\mathfrak{A}, \sigma\right) \Rrightarrow \mathbf{Tricat}_{s}\left(\mathfrak{A}, \tau\right)$ whose component at $F$ is the perturbation $\Omega_{F}$ of Proposition A.9 of \cite{Buhne PhD}. 
		\item The perturbation  $\mathbf{Tricat}_{s}\left(\mathfrak{A}, \Omega\right): \mathbf{Tricat}_{s}\left(\mathfrak{A}, \sigma\right) \Rrightarrow \mathbf{Tricat}_{s}\left(\mathfrak{A}, \tau\right)$ restricts to a perturbation $[\mathfrak{A}, \sigma] \Rrightarrow [\mathfrak{A}, \tau]$.
		\item The assignment $\left(\Omega: \sigma \Rrightarrow \tau\right) \mapsto \left([\mathfrak{A}, \Omega]: [\mathfrak{A}, \sigma] \Rrightarrow [\mathfrak{A}, \tau]\right)$ defines a functor from the category $[\mathfrak{B}, \mathfrak{C}]\left(S, T\right)\left(p, q\right)$ to the hom-category of the $\mathbf{Gray}$-category $ [[\mathfrak{A}, \mathfrak{B}], [\mathfrak{A}, \mathfrak{C}]]$ determined by the morphisms $[\mathfrak{A}, p]$ and $[\mathfrak{A}, q]$.
	\end{enumerate} 
\end{proposition}

\begin{proof}
	For part (1), the perturbation axiom at $j: F \rightarrow G$ is precisely the perturbation axiom for $\Omega$ at each $j_{X}$. The restriction in part (2) is immediate. Part (3) follows directly from the definition of composition of perturbations.
\end{proof}

\begin{remark}\label{un-curried internal whiskering operator in terms of functions}
	Thus far we have completely described the underlying morphism of $3$-globular sets of the internal whiskering operator $[\mathfrak{A}, -]: [\mathfrak{B}, \mathfrak{C}] \rightarrow [[\mathfrak{A}, \mathfrak{B}], [\mathfrak{A}, \mathfrak{C}]]$. We have also shown functoriality at the level of underlying categories (Corollary \ref{underlying functor of internal whiskering operator in closed structure}) and functoriality at the level of hom-categories (Proposition \ref{whiskering perturbations is functorial} part (3)). Before continuing to prove the remaining aspects of $\mathbf{Gray}$-functoriality, we describe how the morphism of $3$-globular sets `un-curries' to various functions between sets in the spirit of Remark \ref{higher cells between Gray categories as functions on globular input data}. The aim of organising these data in this way is to clarify what further calculations are needed in proving $\mathbf{Gray}$-functoriality, and well-definedness of the resulting closed structure.
	\\
	\\
	\noindent By Remark \ref{higher cells between Gray categories as functions on globular input data}, for $n \leq 3$ an $n$-cell $\Phi$ in $[\mathfrak{B}, \mathfrak{C}]$ is sent by $[\mathfrak{A}, -]$ to an $n$-cell in $[[\mathfrak{A}, \mathfrak{B}], [\mathfrak{A}, \mathfrak{C}]]$. We may assume that $\Phi$ is semi-strict if $n=1$, since functoriality from Corollary \ref{underlying functor of internal whiskering operator in closed structure} means that we can extend a comparison of semi-strict trinatural transformations to a comparison of semi-strictly decomposable ones. Now $\Phi$ itself is determined by functions $\Phi_{m}: \mathfrak{B}_{m} \rightarrow \mathfrak{C}_{n+m}$, unless $m=n=1$ in which case the factorisation $\Phi_{1}': \mathfrak{B}_{1} \rightarrow \mathfrak{C}_{\text{eq}}$ is needed to determine $\Phi$. Similarly, the $n$-cell $[\mathfrak{A}, \Phi] \in [[\mathfrak{A}, \mathfrak{B}], [\mathfrak{A}, \mathfrak{C}]]$ is determined by functions $[\mathfrak{A}, \Phi]_{i}: [\mathfrak{A}, \mathfrak{B}]_{i}\rightarrow [\mathfrak{A}, \mathfrak{C}]_{n+i}$ for $i + n \leq 3$, and $[\mathfrak{A}, \Phi]_{1}': [\mathfrak{A}, \mathfrak{B}]_{1} \rightarrow [\mathfrak{A}, \mathfrak{C}]_{eq}$ if $i = n = 1$. But an input $l$-cell $\Psi \in [\mathfrak{A}, \mathfrak{B}]$ is similarly determined by functions $\Psi_{k}: \mathfrak{A}_{k} \rightarrow \mathfrak{B}_{k+l}$ and $\Psi_{k}': \mathfrak{A}_{1} \rightarrow \mathfrak{B}_{\text{eq}}$ if $k = l = 1$. Combining these observations, the morphism of $3$-globular sets $[\mathfrak{A}, -]: [\mathfrak{B}, \mathfrak{C}] \rightarrow [[\mathfrak{A}, \mathfrak{B}], [\mathfrak{A}, \mathfrak{C}]]$ is determined by the composites 
	\begin{tikzcd}
		\mathfrak{A}_{k} \arrow[rr, "\Psi_{k}"] && \mathfrak{B}_{k+l} \arrow[rr, "\Phi_{k+l}"] && \mathfrak{C}_{k+l+n}
	\end{tikzcd} for all $k+l+n \leq 3$, unless $k=l=1$ or $l=n = 1$. In these cases, the functions depicted below are needed to determine the $\mathbf{Gray}$-functor $[\mathfrak{A}, -]: [\mathfrak{B}, \mathfrak{C}] \rightarrow [[\mathfrak{A}, \mathfrak{B}], [\mathfrak{A}, \mathfrak{C}]]$.
	\\
	\\
	\begin{tikzcd}
		\mathfrak{A}_{1} 
		\arrow[rrdd, "\Psi_{1}"']
		\arrow[rr, "\Psi_{1}"]
		&& \mathfrak{B}_{eq}\arrow[dd]
		\arrow[rr, "\Phi_\text{eq}"] && \mathfrak{C}_{eq}
		\arrow[dd]
		\\
		\\
		&& \mathfrak{B}_{2} 
		\arrow[rr, "\Phi_{2}"'] && \mathfrak{C}_{2}&{}
	\end{tikzcd}\begin{tikzcd}
		\mathfrak{A}_{1} 
		\arrow[rr, "\Psi_{1}"] && \mathfrak{B}_{1}
		\arrow[rr, "\Phi_\text{1}'"] && \mathfrak{C}_{eq}
	\end{tikzcd}
\end{remark}

\begin{proposition}\label{local 2 functoriality of whiskering}
	The functors of Proposition \ref{whiskering perturbations is functorial} are the functors between hom-categories of a $2$-functor $[\mathfrak{B}, \mathfrak{C}]\left(S, T\right) \rightarrow [\mathfrak{L}, \mathfrak{K}]\left([\mathfrak{A}, S], [\mathfrak{A}, T]\right)$, where $\mathfrak{L} := [\mathfrak{A}, \mathfrak{B}]$ and $\mathfrak{K}:= [\mathfrak{A}, \mathfrak{C}]$. The action on objects of this $2$-functor is given by $p \mapsto [\mathfrak{A}, p]$. 
\end{proposition}

\begin{proof}
	We first check functoriality between underlying categories of the required $2$-functor. Preservation of identities follows from Proposition \ref{interchanging tritransformation with trimodification} part (2). For preservation of composition, we need to show that for any \begin{tikzcd}
		p \arrow[r, "\sigma"] & q \arrow[r, "\tau"] & r
	\end{tikzcd}, the trimodifications $[\mathfrak{A}, \tau \circ \sigma]$ and $[\mathfrak{A}, \tau] \circ [\mathfrak{A}, \sigma]$ have the same $2$-cell component for every $F: \mathfrak{A} \rightarrow \mathfrak{B}$, and the same $3$-cell component for every $j: F \rightarrow G$ in $[\mathfrak{A}, \mathfrak{B}]$. But observe that their $2$-cell components will both be the trimodification depicted below left and their $3$-cell components will both be the pasting of perturbations depicted below centre. This proves functoriality between underlying categories. Finally, $2$-functoriality follows directly from the definition of horizontal composition of perturbations. It says that on a pair of horizontally composable perturbations depicted below right, the components of both $[\mathfrak{A}, \Omega]*[\mathfrak{A}, \Sigma]$ and $[\mathfrak{A}, \Omega * \Sigma]$ at a $\mathbf{Gray}$-functor $F: \mathfrak{A} \rightarrow \mathfrak{B}$ are given by the same perturbation. This asks for the equality ${\left(\Omega * \Sigma\right)}_{F} = \Omega_{F}*\Sigma_{F}$ to hold. The fact that it does hold is clear since the horizontal composite of perturbations is defined as the component-wise horizontal composite in the hom-$2$-categories of $\mathfrak{B}$.
	
	$$\begin{tikzcd}
		p_{F} \arrow[r, "\sigma_{F}"] & q_{F} \arrow[r, "\tau_{F}"] & r_{F} &{}
	\end{tikzcd}  \begin{tikzcd}[column sep = 15]
		{p}_{F}
		\arrow[dd, "{p}_{j}"]
		\arrow[rr, "{\sigma}_{F}"] &
		{} \arrow[dd, Rightarrow, shorten = 15, "{\sigma}_{j}"]
		& {q}_{F} \arrow[dd, "q_{j}"]
		\arrow[rr, "{\tau}_{F}"] &
		{}\arrow[dd, Rightarrow, shorten = 15, "{\tau}_{j}"]
		& {r}_{F} \arrow[dd, "{r}_{j}"]
		\\
		\\
		{p}_{G} 
		\arrow[rr, "\sigma_{G}"'] &{}& 
		{q}_{G} 
		\arrow[rr, "{\tau}_{G}"'] 
		&{}& {r}_{G}&{}
	\end{tikzcd}
	\begin{tikzcd}
		p
		\arrow[rr, bend left, "\sigma" name = A]
		\arrow[rr, bend right, "\sigma '"' name = B]
		&& q
		\arrow[rr,bend left, "\tau" name = C]
		\arrow[rr, bend right, "\tau '"' name = D]
		&& r
		\arrow[from = A, to = B, Rightarrow, shorten = 5, "\Sigma"]
		\arrow[from = C, to = D, Rightarrow, shorten =5, "\Omega"]
	\end{tikzcd}$$
	
\end{proof}

\noindent We are now ready to establish the last and most complicated fragment of the closed structure on $\mathbf{Gray}$-$\mathbf{Cat}$, namely $\mathbf{Gray}$-functoriality of $[\mathfrak{B}, \mathfrak{C}] \rightarrow [[\mathfrak{A}, \mathfrak{B}], [\mathfrak{A}, \mathfrak{C}]]$ and its various (extra)naturality properties. The proof of $\mathbf{Gray}$-functoriality follows from unwinding definitions, while the proof of the (extra)naturality properties uses the familiar phenomenon of `associativity of substitution'. This amounts to using associativity of the functions of the form $\Psi_{l}: \mathfrak{A}_{k} \rightarrow \mathfrak{B}_{k+l}$, as described in Remark \ref{higher cells between Gray categories as functions on globular input data}.

\begin{theorem}\label{internal whiskering operators closed structure on Gray cat}
	\hspace{1mm}
	\begin{enumerate}
		\item 
		The functor $[\mathfrak{A}, -]_{1}: [\mathfrak{B}, \mathfrak{C}]_{1} \rightarrow [[\mathfrak{A}, \mathfrak{B}], [\mathfrak{A}, \mathfrak{C}]]_{1}$ of Corollary \ref{underlying functor of internal whiskering operator in closed structure} enriches to define a $\mathbf{Gray}$-functor whose actions between hom-$2$-categories are the $2$-functors of Proposition \ref{local 2 functoriality of whiskering}.
		\item The assignment of $\mathfrak{B}$ to the $\mathbf{Gray}$-functor $[\mathfrak{A}, -]: [\mathfrak{B}, \mathfrak{C}] \rightarrow [[\mathfrak{A}, \mathfrak{B}], [\mathfrak{A}, \mathfrak{C}]]$ is natural as $\mathfrak{B}$ varies along a $\mathbf{Gray}$-functor $B: \mathfrak{B}' \rightarrow \mathfrak{B}$.
		\item The assignment of $\mathfrak{C}$ to the $\mathbf{Gray}$-functor $[\mathfrak{A}, -]: [\mathfrak{B}, \mathfrak{C}] \rightarrow [[\mathfrak{A}, \mathfrak{B}], [\mathfrak{A}, \mathfrak{C}]]$ is natural as $\mathfrak{C}$ varies along a $\mathbf{Gray}$-functor $C: \mathfrak{C} \rightarrow \mathfrak{C}'$ .
		\item The assignment of $\mathfrak{A}$ to the $\mathbf{Gray}$-functor $[\mathfrak{A}, -]: [\mathfrak{B}, \mathfrak{C}] \rightarrow [[\mathfrak{A}, \mathfrak{B}], [\mathfrak{A}, \mathfrak{C}]]$ is extra-natural as $\mathfrak{A}$ varies along a $\mathbf{Gray}$-functor $A: \mathfrak{A} \rightarrow \mathfrak{A}'$.
	\end{enumerate}
\end{theorem}

\begin{proof}
	For part (1), we need to check respect for composition. We describe how all aspects follow from `associativity of substitution'. Consider the following interchangeable pair of trimodifications in $[\mathfrak{B}, \mathfrak{C}]$: \begin{tikzcd}
		R \arrow[rr, bend right, "r"' name = G]\arrow[rr, bend left, "s" name = F] &&S \arrow[rr, bend left, "p" name = S]\arrow[rr, bend right, "q"' name = T] && T
		\arrow[from = F, to = G, Rightarrow, shorten = 5, "\sigma"]
		\arrow[from = S, to = T, Rightarrow, shorten = 5, "\tau"]
	\end{tikzcd}.
	
	\begin{itemize}
		\item For the whiskering of $\tau$ with $r$, both trimodifications $[\mathfrak{A}, \tau.r]$ and $[\mathfrak{A}, \tau].[\mathfrak{A}, r]$ have component $2$-cell on $F: \mathfrak{A} \rightarrow \mathfrak{B}$ given by the trimodification  depicted below left, and component $3$-cell on $j: F \rightarrow G$ given by the perturbation depicted below right.
	\end{itemize}
	$$\begin{tikzcd}
		RF \arrow[rr, "r_{F}"]&&SF \arrow[rr, bend left, "pF" name = S]\arrow[rr, bend right, "qF"' name = T] && T
		\arrow[from = S, to = T, Rightarrow, shorten = 5, "\tau_{F}"] &{}
	\end{tikzcd} \begin{tikzcd}[column sep = 12]
		p_{G}.r_{G}.Rj
		\arrow[dd,"\tau_{G}.1.1"']
		\arrow[rr,"1.r_j"]
		&
		{}
		\arrow[dd,Rightarrow,shorten=15,"{\left(\tau_{G}\right)}_{\left(r_{j}\right)}", shift right = 5]
		&
		p_{G}.Sj.r_{F} \arrow[rr, "p_{j}.1"]
		\arrow[dd,"\tau_{G}.1.1"]
		&{}\arrow[dd, Rightarrow, shorten = 15, "\tau_{j}.1"]
		& Tj.p_{F}.r_{F}\arrow[dd, "1.\tau_{F}.1"]
		\\
		\\
		{q}_{G}{r}_{G}.Rj
		\arrow[rr,"1.{r}_{j}"']
		&
		{}
		&
		{q}_{G}.Sj.{r}_{F} \arrow[rr, "{q}_{j}.1"']
		&{}
		& Tj.{q}_{F}.{r}_{F}
	\end{tikzcd}$$
	\begin{itemize}
		\item If $\Omega: \tau \Rrightarrow \tau'$ is a perturbation, then for the whiskering of $\Omega$ with $r$ we find that the perturbations $[\mathfrak{A}, \Omega.r]$ and $[\mathfrak{A}, \Omega].[\mathfrak{A}, p]$ have $3$-cell component on $F: \mathfrak{A} \rightarrow \mathfrak{B}$ given by the perturbation $\Omega_{F}.r_{F}$; the whiskering of the perturbation $\Omega_{F}$ with $r_{F}$. 
		\item Preservation of the whiskering of $\sigma$ with $p$ and of the whiskering of a perturbation $\Sigma: \sigma \Rrightarrow \sigma'$ with $p$ follows by a similar analysis to that in the previous two dot points.
		\item For interchanger of $\tau$ with $\sigma$, we need to show that the interchanger of the $2$-cells $[\mathfrak{A}, \sigma]$ and $[\mathfrak{A}, \tau]$ in the $\mathbf{Gray}$-category $[[\mathfrak{A}, \mathfrak{B}], [\mathfrak{A}, \mathfrak{C}]]$ is the same perturbation as the one to which the $\mathbf{Gray}$-functor $[\mathfrak{A}, -]$ maps the interchanger $\tau_\sigma \in [\mathfrak{B}, \mathfrak{C}]$. The first of these will have component on $F \in [\mathfrak{A}, \mathfrak{B}]$ given by the interchanger of $\tau_{F}$ and $\sigma_{F}$, and this will in turn be the perturbation whose component on $X \in \mathfrak{A}$ will be given by the interchanger ${\left(\tau_{FX}\right)}_{\left(\sigma_{FX}\right)}$ in $\mathfrak{C}$. The second of these will have component on $F\in [\mathfrak{A}, \mathfrak{B}]$ given by the perturbation whose component on $X$ will be the component of $\tau_{\sigma}$ on $FX$, and by part (4 e) of Remark \ref{Gray Functor Gray category with weak higher cells}, this is seen to also be the interchanger ${\left(\tau_{FX}\right)}_{\left(\sigma_{FX}\right)}$ in $\mathfrak{C}$.
	\end{itemize}
	
	\noindent The proofs for parts (2), (3) and (4) all also use associativity of substitution. Fix $0 \leq n \leq 3$ and let $\Phi$ be an $n$-cell in $[\mathfrak{B}, \mathfrak{C}]$. Part (2) says that the $n$-cell $[\mathfrak{A}, \Phi_{B}]$ in $[[\mathfrak{A}, \mathfrak{B}'], [\mathfrak{A}, \mathfrak{C}]]$ has the same components on data in $[\mathfrak{A}, \mathfrak{B}']$ as the restriction of the $n$ -cell $[\mathfrak{A}, \Phi]$ along the $\mathbf{Gray}$-functor $[\mathfrak{A}, B]: [\mathfrak{A}, \mathfrak{B}'] \rightarrow [\mathfrak{A}, \mathfrak{B}]$. Part (3) says that the $n$-cell  $[\mathfrak{A}, C\Phi]$ in $[[\mathfrak{A}, \mathfrak{B}], [\mathfrak{A}, \mathfrak{C}']]$ has the same components on data in $[\mathfrak{A}, \mathfrak{B}]$ as the extension of $[\mathfrak{A}, \Phi]$ along the $\mathbf{Gray}$-functor $[\mathfrak{A}, C]: [\mathfrak{A}, \mathfrak{C}] \rightarrow [\mathfrak{A}, \mathfrak{C}']$. Finally, part (4) says something very similar to part (2), namely that the $n$-cell $[\mathfrak{A}, \Phi_A]$ has the same components on data in $[\mathfrak{A}, \mathfrak{B}]$ as the $n$-cell $[\mathfrak{A}', \Phi]$ restricted along $[A, \Phi]$.
\end{proof}

\section{The axioms of the closed structure}\label{Axioms for closed structure}

\subsection{Axioms other than associativity}

\begin{proposition}\label{Non-associativity axioms for closed structure on Gray cat}
	Let $\mathfrak{A}$, $\mathfrak{B}$ be $\mathbf{Gray}$-categories and recall the families of $\mathbf{Gray}$-functors $i_\mathfrak{A}: \mathbf{1} \rightarrow [\mathfrak{A}, \mathfrak{A}]$ of Proposition \ref{identity assigners in closed structure}, $c_\mathfrak{A}: \mathfrak{A} \rightarrow  [\mathbf{1}, \mathfrak{A}]$ of Proposition \ref{constant assigners closed structure} and $[\mathfrak{A}, -]: [\mathfrak{B}, \mathfrak{C}] \rightarrow [[\mathfrak{A}, \mathfrak{B}], [\mathfrak{A}, \mathfrak{C}]]$ of Theorem \ref{internal whiskering operators closed structure on Gray cat}. Then axioms (1), (2), (3) and (4) in Definition \ref{definition of a closed structure on a category} hold for these data.
\end{proposition}

\begin{proof}
	Axiom (1) was established as Proposition \ref{identity assigners in closed structure} part (2). We explain how the commutativity of each of the diagrams required for axioms (2), (3) and (4) follow once the definitions are unwound. Axiom (2) just says that if $\Phi$ is any $n$-cell for $0 \leq n \leq 3$ in $[\mathfrak{A}, \mathfrak{B}]$ then it does not change when whiskered with the identity $\mathbf{Gray}$-functor $1_\mathfrak{B}$. Axiom (3) says that the component of $[\mathfrak{A}, \Phi] \in [[\mathfrak{A}, \mathfrak{A}], [\mathfrak{A}, \mathfrak{B}]]$ on the identity $\mathbf{Gray}$-functor $1_\mathfrak{A}$ is just $\Phi$ itself. Axiom (4) is a strengthening of Proposition \ref{constant assigners closed structure} part (3). It says that the family of $\mathbf{Gray}$-functors $c_\mathfrak{A}: \mathfrak{A} \rightarrow [\mathbf{1}, \mathfrak{A}]$ is not only natural with respect to $\mathbf{Gray}$-functors, but also with respect to semi-strictly decomposable trinatural transformations, trimodifications and perturbations. But this strengthening is also clear from the definition of $c$.
\end{proof}

\subsection{The associativity axiom}\label{Subsection the associativity axiom}

\noindent The next few propositions establish axiom (5) of Definition \ref{definition of a closed structure on a category} one dimension at a time. This involves
\begin{itemize}
	\item considering every triple $\left(\Theta, \Psi, \Phi\right)$ consisting of an $n$-cell $\Theta$ in $[\mathfrak{C}, \mathfrak{D}]$, an $m$-cell $\Psi$ in $[\mathfrak{B}, \mathfrak{C}]$ and an $l$-cell $\Phi$ in $[\mathfrak{A}, \mathfrak{B}]$ satisfying $l + m + n \leq 3$
	\item establishing that the two ways of combining this data to get an $\left(l + m + n\right)$-cell in $[\mathfrak{A}, \mathfrak{D}]$ give the same result.
\end{itemize} 

\noindent By Remark \ref{higher cells between Gray categories as functions on globular input data}, most of this will routinely follow from associativity of the functions depicted below for $k+l+m+n \leq 3$. The one exception to this is when $k=0$ and $l = m = n = 1$, which will use semi-strict decomposability of $\Theta$.

\begin{equation}\label{composable triple of functions from A k to D k+l+m+n}
	\begin{tikzcd}
		\mathfrak{A}_{k} \arrow[rr, "\Phi_{k}"] && \mathfrak{B}_{k+l} \arrow[rr, "\Psi_{k+l}"] && \mathfrak{C}_{k+l+m} \arrow[rr, "\Theta_{k+l+m}"] &&\mathfrak{D}_{k+l+m+n}
	\end{tikzcd}
\end{equation}

\noindent Doing this `one-dimension at a time' means that we will be taking the following route.

\begin{itemize}
	\item Proposition \ref{associativity condition underlying sets} will prove the associativity condition on underlying sets, and hence address the cases $n = 0$, $l + m \leq 3$.
	\item Proposition \ref{assoc condition underlying categories} will strengthen Proposition \ref{associativity condition underlying sets} to prove the associativity condition at the level of underlying categories, and hence address the cases $n = 1$, $l+m \leq 2$.
	\item  Proposition \ref{associativity condition for closed structure underlying sesquicategories} will strengthen Proposition \ref{assoc condition underlying categories} to prove the associativity condition at the level of underlying sesquicategories (i.e. also on $2$-cells), and hence address the cases $n = 2$, $l+m \leq 1$.
	\item Finally, the proof of Theorem \ref{Gray cat closed structure final theorem} will just need to check the case $n = 3$, $l= m=0$.
\end{itemize}

\begin{proposition}\label{associativity condition underlying sets}
	The associativity condition holds at the level of underlying sets.
\end{proposition}

\begin{proof}
	We describe how everything follows easily from associativity of composition of the functions in Equation \ref{composable triple of functions from A k to D k+l+m+n}. Let $V : \mathfrak{C} \rightarrow \mathfrak{D}$ be a $\mathbf{Gray}$-functor, let $\bar{V}: [\mathfrak{B}, \mathfrak{C}] \rightarrow [[\mathfrak{A}, \mathfrak{B}], [\mathfrak{A}, \mathfrak{D}]]$ be its image under the clockwise path in Definition \ref{definition of a closed structure on a category} axiom (5), and let $\underline{V}: [\mathfrak{B}, \mathfrak{C}] \rightarrow [[\mathfrak{A}, \mathfrak{B}], [\mathfrak{A}, \mathfrak{D}]]$ be its image under the anti-clockwise path. We need to show that $\bar{V} = \underline{V}$ on all data in $[\mathfrak{B}, \mathfrak{C}]$. This says that for an $m$-cell $\Psi \in [\mathfrak{B}, \mathfrak{C}]$ and an $l$-cell $\Phi \in [\mathfrak{A}, \mathfrak{B}]$ with $m + l \leq 3$ so that the $\left(l +m\right)$-cell $\Psi_\Phi \in [\mathfrak{A}, \mathfrak{C}]$ is defined, the equality $V \left(\Psi_\Phi \right) = {\left(V \Psi\right)}_{\Phi}$ holds. But this is clear by definition of whiskering with $V$.
\end{proof}

\begin{proposition}\label{assoc condition underlying categories}
	Consider arbitrary trinatural transformations as displayed below. $$\begin{tikzcd}
		\mathfrak{A} \arrow[rr, bend right, "G"' name = G]\arrow[rr, bend left, "F" name = F] &&\mathfrak{B} \arrow[rr, bend left, "S" name = S]\arrow[rr, bend right, "T"' name = T] && \mathfrak{C}\arrow[rr, bend left, "U" name = U]\arrow[rr, bend right, "V"' name = V]  && \mathfrak{D}
		\arrow[from = F, to = G, Rightarrow, shorten = 5, "j"]
		\arrow[from = S, to = T, Rightarrow, shorten = 5, "p"]
		\arrow[from = U, to = V, Rightarrow, shorten = 5, "w"]
	\end{tikzcd}$$
	\begin{enumerate}
		\item Let $w_{p}$ be the trimodification which mediates interchange between $w$ and $p$ as defined in Proposition \ref{interchanging tritransformations} and let $p_{j}$ be the trimodification defined similarly from $p$ and $j$. Let ${\left(w_{p}\right)}_{j}$ be the perturbation which mediates interchange between the trimodification $w_{p}$ and the trinatural transformation $j$ as per Proposition \ref{interchanging tritransformation with trimodification}, and let $w_{\left(p_{j}\right)}$ be the perturbation which mediates interchange between the trinatural transformation $w$ and the trimodification $p_{j}$ as per Proposition \ref{interchanging trimodification with tritransformation}.
		
		\begin{enumerate}
			\item The perturbation ${\left(w_{p}\right)}_{j}$ is given by the pasting of perturbations depicted below, where $w_{p_{G}, S_{j}}$ and $w_{Tj, p_{F}}$ are the perturbations defined as in Proposition \ref{compositor of whiskering with p}.
			$$\begin{tikzcd}[column sep= 40, font=\fontsize{9}{6}]
				w_{TG}.Up_{G}.USj \arrow[rr, "w_{p_{G}}.1"]
				\arrow[dd, "1.Up_{j}"']
				\arrow[rrdd, "w_{p_{G}.Sj}" description]
				&& Vp_{G}.w_{SG}.USj \arrow[dd, "1.w_{Sj}"]
				\\
				&{}\arrow[r, Rightarrow, shorten = 15, "w_{p_G{,}Sj}^{-1}" shift left = 5]&{}
				\\
				w_{TG}.UTj.Up_{F} \arrow[rr, shorten = 15, Rightarrow, "w_{\left(p_j\right)}"]
				\arrow[rrdd, "w_{Tj.p_{F}}" description]
				\arrow[dd, "w_{Tj}.1"'] 
				&& Vp_{G}.VSj.w_{SF}\arrow[dd, "Vp_{j}.1"]
				\\
				{}\arrow[r, Rightarrow, shorten =15, "w_{Tj{,}p_F}"', shift right = 5]&{}
				\\
				VTj.w_{TF}.Up_{F} \arrow[rr, "1. w_{p_F}"'] && VTj.Vp_{F}.w_{SF}
			\end{tikzcd}$$
			\item The equation ${\left(w_{p}\right)}_{j} = {w}_{\left(p_{j}\right)}$ holds under either of the following hypotheses.
			
			\begin{enumerate}
				\item $w$ is semi-strict, or
				\item $j$ is a pseudo-icon equivalence and $p$ is unital. In this case both ${\left(w_{p}\right)}_{j}$ and  ${w}_{\left(p_{j}\right)}$ are the identity. 
			\end{enumerate}
		\end{enumerate}
		\item 
		The associativity condition holds at the level of underlying categories.
	\end{enumerate}
\end{proposition}

\begin{proof}
	Part (1 a) is an exercise in unwinding definitions. For each $X \in \mathfrak{A}$, the trimodification $p_{j}$ has $2$-cell component on $X$ given by $p_{j_{X}}$. Hence the perturbation $w_{\left(p_{j}\right)}$ has component on $X$ given by the $3$-cell component of $w$ on the $2$-cell $p_{j_{X}}$. On the other hand, the trimodification $w_p$ has $3$-cell component on $j_{X}$ given as depicted. Part (1 b i) follows from part (1 a), since semi-strictness means that in particular the compositors $w_{Tj, p_{F}}$ and $w_{p_{G}, Sj}$ are identities. Meanwhile for part (1 b ii), on the one hand $j$ being a pseudo-icon equivalence means that $j_{X}$ is the identity, so the unit law for the trimodification $w_{p}$ says that $\left(w_{p}\right)_{j_{X}}$ is the identity. On the other hand, $p$ being unital means that if $j_{X}$ is the identity then $p_{j_{X}}$ is also the identity, so the unit aspect of local pseudonaturality of $w$ says that $w_{\left(p_{j_{X}}\right)}$ will also be the identity. Finally, the left and right unit laws for $w$ and identity preservation of $S$ and $T$ say that since $j_{X}$ is the identity, the compositor components $w_{Tj_{X}, p_{FX}}$ and $w_{p_{GX}, Sj_{X}}^{-1}$ are also the identity.
	\\
	\\
	\noindent We now turn to part (2). By Proposition \ref{associativity condition underlying sets}, it remains only to show that the associativity condition holds for morphisms. By functoriality established in Corollary \ref{underlying functor of internal whiskering operator in closed structure}, since every semi-strictly decomposable trinatural transformation is a composite of semi-strict ones, it suffices to show that the associativity condition holds for a semi-strict trinatural transformation $w: U \rightarrow V$. Hence suppose $w$ is semi-strict, let $\bar{w}$ denote the image of $w$ traversing the pentagon in Definition \ref{definition of a closed structure on a category} axiom (5) clockwise, and let $\underline{w}$ be its image traversing the diagram anti-clockwise. First note that $\underline{w}$ and $\bar{w}$ are both also semi-strict, and that their components on an object $S \in [\mathfrak{B}, \mathfrak{C}]$ are both the same trinatural transformation. Next, the components of $\bar{w}$ and $\underline{w}$ on a morphism $p: S \rightarrow T$ are given by the trimodifications which agree on objects and whose $3$-cell components on a morphism $j: F \rightarrow G$ are the perturbations ${w}_{\left(p_{j}\right)}$ and ${\left(w_{p}\right)}_{j}$. But by Part (1, b), these are equal since $w$ is assumed to be semi-strict. Finally, the components of both trinatural transformations on a $2$-cell $\sigma: p \Rightarrow q$ are given by the perturbation whose component on $F \in [\mathfrak{A}, \mathfrak{B}]$ is the perturbation $w_{{\sigma}_{F}}$.
\end{proof}

\begin{proposition}\label{associativity condition for closed structure underlying sesquicategories}
	Let $\chi: w \Rightarrow x: U \rightarrow V$ be a $2$-cell in $[\mathfrak{C}, \mathfrak{D}]$, let $\bar{\chi}$ denote its image under the clockwise path in Definition \ref{definition of a closed structure on a category} axiom 5, and let $\underline{\chi}$ denote its image under the anti-clockwise path. Then $\bar{\chi} = \underline{\chi}$. 
\end{proposition}

\begin{proof}
	This asks for each of the following equations to hold for $p: S \rightarrow T$ in $[\mathfrak{B}, \mathfrak{C}]$ and $j: F \rightarrow G$ in $[\mathfrak{A}, \mathfrak{B}]$. The fact that they do hold follows from associativity of composition of the functions in Equation \ref{composable triple of functions from A k to D k+l+m+n}, with $l+m \leq 1$ and $n = 1$.
	
	$$\begin{tikzcd}
		\chi_{\left(SF\right)} = {\left(\chi_{S}\right)}_{F}
		&&\chi_{\left(Sj\right)} = {\left(\chi_{S}\right)}_{j} &&\chi_{\left(p_F\right)} = {\left(\chi_{p}\right)}_{F}
	\end{tikzcd}$$ 
\end{proof}

\begin{theorem}\label{Gray cat closed structure final theorem}
	There is a closed structure on the category $\mathcal{E} = \mathbf{Gray}$-$\mathbf{Cat}$ whose \begin{enumerate}
		\item Internal hom \begin{tikzcd}\mathcal{E}^\text{op}\times \mathcal{E} \arrow[r, "{[}-{,}?{]}"] & \mathcal{E}
		\end{tikzcd} is as described in Proposition \ref{coherent hom on Gray Cat}.
		\item Unit object is the terminal $\mathbf{Gray}$-category $\mathbf{1}$,
		\item Identity assigners are described in Proposition \ref{identity assigners in closed structure},
		\item Constant map assigners are described in Proposition \ref{constant assigners closed structure},
		\item Internal whiskering operators \begin{tikzcd}
			{[}\mathfrak{B}{,}\mathfrak{C}{]} \arrow[rr, "{[}\mathfrak{A}{,}-{]}"] && {[}{[}\mathfrak{A}{,}\mathfrak{B}{]}{,}{[}\mathfrak{A}{,}\mathfrak{C}{]}{]}
		\end{tikzcd} are described in Theorem \ref{internal whiskering operators closed structure on Gray cat}.
	\end{enumerate}
\end{theorem}

\begin{proof}
	All that remains is to show that the associativity condition holds on $3$-cells $\Omega: \chi \Rrightarrow \zeta$ in $[\mathfrak{C}, \mathfrak{D}]$. But the equation $\Omega_{\left(SF\right)} = {\left(\Omega_{S}\right)}_{F}$ can again be seen to hold by associativity of composition of the functions in Equation \ref{composable triple of functions from A k to D k+l+m+n}, with $l = m = 0$ and $n = 3$.
\end{proof}

\begin{remark}\label{undesirable aspects of enrichment}
	\noindent We defer detailed analysis of enrichment and models of four-dimensional categories to a forthcoming paper \cite{Miranda weak interchange 4-categories}\footnote{Also see Chapter 3 of \cite{Miranda PhD}}. However, we record here an undesirable feature that such enriched categories appear to have. In particular, let $\mathscr{C}$ be such a category enriched over $\left(\mathbf{Gray}\text{-}\mathbf{Cat}, [-, ?]\right)$, let $X \in \mathscr{C}$ be an object and let $\phi: f \rightarrow g$ be a morphism in a hom-$\mathbf{Gray}$-category $\mathscr{C}\left(Y, Z\right)$. As part of its structure, $\mathscr{C}$ has a $\mathbf{Gray}$-functor $\mathscr{C}\left(X, -\right): \mathscr{C}\left(Y, Z\right) \rightarrow [\mathscr{C}\left(X, Y\right), \mathscr{C}\left(X, Z\right)]$ under which we can consider the image of $\phi$. This must be a semi-strictly generated trinatural transformation $\mathscr{C}\left(X, \phi\right): \mathscr{C}\left(X, f\right) \rightarrow \mathscr{C}\left(X, g\right)$. It does not appear to be possible to encode the requirement that $\mathscr{C}\left(X, \phi\right)$ should be semi-strictly decomposable in terms of equations on data in $\mathscr{C}$. In particular, it appears as though $\left(\mathbf{Gray}\text{-}\mathbf{Cat}, [-, ?]\right)\text{-}\mathbf{Cat}$ may not be a locally presentable category. In a forthcoming paper \cite{Miranda weak interchange 4-categories}, we will describe a different base whose enriched categories will be categories enriched over $\left(\mathbf{Gray}\text{-}\mathbf{Cat}, [-, ?]\right)$ with better algebraic properties.
\end{remark}

\begin{remark}\label{Remark relation to Crans}
	Although the semi-strictly generated closed structure of Theorem \ref{Gray cat closed structure final theorem} does not have a corresponding monoidal structure, it does have a restricted relationship to a pseudo variant of the Crans monoidal structure defined in \cite{crans tensor of gray categories}. Specifically, we consider the variant in which 
	
	\begin{itemize}
		\item the generating $2$-cells $\left(f, g\right): \left(g, X'\right).\left(f, Y\right) \Rightarrow \left(f, Y'\right).\left(g, X\right)$ in $\mathfrak{A} \otimes \mathfrak{B}$ determined by morphisms $f: X \rightarrow X'$ in $\mathfrak{A}$ and $g: Y \rightarrow Y'$ in $\mathfrak{B}$, described in the third dot point of Section 3.1 of \cite{crans tensor of gray categories}, are part of adjoint equivalences.
		\item the generating $3$-cells described in the fourth and fifth dot points of Section 3.1 of \cite{crans tensor of gray categories}, and displayed as $2$-cells in certain hom-$2$-categories of $\mathfrak{A} \otimes \mathfrak{B}$ below left and below right respectively, are invertible.
	\end{itemize}

	$$\begin{tikzcd}[column sep = 12]
		\left(X'{,}g\right).\left(f{,}Y\right)
		\arrow[rrrr,"1.\left(\phi{,}Y\right)"]
		\arrow[dd,"\left(f{,} g\right)"']
		&&
		{}
		\arrow[dd,Rightarrow,shorten=15,"\left(\phi{,} g\right)"]
		&&
		\left(X'{,}g\right).\left(f'{,}Y\right)
		\arrow[dd,"\left(f'{,}g\right)"]
		\\
		\\
		\left(f{,}Y'\right).\left(X{,}g\right)
		\arrow[rrrr,"\left(\phi{,}Y'\right).1"']
		&&
		{}
		&&
		\left(f'{,}Y'\right).\left(X{,}g\right) &{}
	\end{tikzcd}\begin{tikzcd}[column sep = 12]
		\left(X'{,}g\right).\left(f{,}Y\right)
		\arrow[rrrr,"\left(X{,}\psi\right).1"]
		\arrow[dd,"\left(f{,} g\right)"']
		&&
		{}
		\arrow[dd,Rightarrow,shorten=15,"\left(f{,}\psi\right)"]
		&&
		\left(X'{,}g'\right).\left(f{,}Y\right)
		\arrow[dd,"\left(f{,}g'\right)"]
		\\
		\\
		\left(f{,}Y'\right).\left(X{,}g\right)
		\arrow[rrrr,"1.\left(X{,}\psi\right)"']
		&&
		{}
		&&
		\left(f{,}Y'\right).\left(X{,}g'\right) &{}
	\end{tikzcd}$$

	\noindent A detailed presentation of this variant of the Crans tensor product is described in Definition 3.1.1 of \cite{Miranda PhD}.
	\\
	\\
	\noindent Recall also from Definition 5.1.2 of \cite{Miranda strictifying operational coherences} that a $\mathbf{Gray}$-category is said to have \emph{only trivial $1$-cell composition} if for any composable pair of morphisms \begin{tikzcd}
	X \arrow[r, "f"] & Y \arrow[r, "g"] & Z
	\end{tikzcd} either $f$ or $g$ is an identity. Let $\mathfrak{A}$ be a $\mathbf{Gray}$-category with only trivial $1$-cell composition. Then the functor $[\mathfrak{A}, -]: \mathbf{Gray}$-$\mathbf{Cat} \rightarrow \mathbf{Gray}$-$\mathbf{Cat}$ is right adjoint to the functor $\mathfrak{A} \otimes -$, with counit component at $\mathfrak{B}$ given by a $\mathbf{Gray}$-functor $\varepsilon_{\mathfrak{B}}^{\mathfrak{A}}: \mathfrak{A} \otimes [\mathfrak{A}, \mathfrak{B}] \to \mathfrak{B}$ which maps generating data $(\alpha, \Phi)$ to the evaluation of $\Phi$ at $\alpha$.
\\
\\
\noindent As such, the $\mathbf{Gray}$-functors $\varepsilon_{\mathfrak{B}}^{\mathfrak{A}}$ exhibit adjunctions $\mathfrak{A} \otimes - \dashv [\mathfrak{A}, -]$ whenever $\mathfrak{A}$ has trivial $1$-cell composition. Moreover, it can be checked that the multicategories represented by the pseudo-Crans monoidal structure $\otimes$ and the semi-strictly generated closed structure $[-, ?]$ are isomorphic, when restricted to such $\mathbf{Gray}$-categories. Most of the details in this correspondence are provided in Proposition 3.1.3 and Corollary 3.1.4 of \cite{Miranda PhD}.
\end{remark}

\section{Relation to the Bourke-Lobbia skew structure}\label{subsection interaction with the skew structure}

\noindent In \cite{Bourke Lobbia Skew Approach to Enrichment for Gray-categories}, several closed, skew monoidal structures are described on $\mathbf{Gray}$-$\mathbf{Cat}$. In this paper we have also described a closed structure on $\mathbf{Gray}$-$\mathbf{Cat}$. Our internal homs restrict the morphisms in the sharp homs $\mathbf{Ps}\left(\mathfrak{A}, \mathfrak{B}\right)_\text{s}$ of \cite{Bourke Lobbia Skew Approach to Enrichment for Gray-categories}, from unital trinatural transformations to just those which admit decompositions into semi-strict trinatural transformations.

\begin{remark}\label{comparing techniques to Bourke Lobbia}
	Our work in Sections \ref{Data of the closed structure} and \ref{Axioms for closed structure} almost serve to re-prove the closed fragment of the closed, skew monoidal structure $\mathbf{Ps}\left(-, ?\right)_{\text{s}}$ on $\mathbf{Gray}$-$\mathbf{Cat}$ defined in \cite{Bourke Lobbia Skew Approach to Enrichment for Gray-categories} in a `first-principles' way, avoiding the use of skew multicategories. The main gap in this approach is the proof of the associativity axiom for the case $k=0$, $l=m=n=1$. Our proof relies on semi-strict decomposability of the trinatural transformation $w: U \rightarrow V$. Although the equation $(w_{p})_{j} = w_{(p_{j})}$ only type-checks when the compositors $w_{p_{G}, Sj}$ and $w_{Tj, p_{F}}$ are identities, axiom (5) of Definition \ref{definition of a closed structure on a category} also holds for the unital trinatural transformations considered in \cite{Bourke Lobbia Skew Approach to Enrichment for Gray-categories}.
\end{remark}

\begin{proposition}\label{proposition relationship to bourke lobbia}
	Consider the inclusion $\mathbf{Gray}$-functors $\chi_{\mathfrak{B}, \mathfrak{C}}: [\mathfrak{B}, \mathfrak{C}] \rightarrow \mathbf{Ps}\left(\mathfrak{B}, \mathfrak{C}\right)_{s}$.
	\begin{enumerate}
		\item This family of functors equips the identity functor on $\mathbf{Gray}$-$\mathbf{Cat}$ with the structure of a normal closed functor in the sense of Definition \ref{Definition normal skew closed functor}.
		\item If the underlying sesquicategory of $\mathfrak{B}$ is free on a $2$-computad, \footnote{equivalently, if $\mathfrak{B}$ is a cofibrant $\mathbf{Gray}$-category in the sense of \cite{Quillen Gray Cat}, or if $\mathfrak{B}$ can be described via presentations in which the only relations are between $3$-cells.} then the $\mathbf{Gray}$-functor $\chi_{\mathfrak{B}, \mathfrak{C}}$ is a biequivalence internal to the tricategory of tricategories defined in \cite{Low dimensional structures formed by tricategories} and recalled in Proposition 4.3.4 of \cite{Miranda strictifying operational coherences}.
	\end{enumerate} 
\end{proposition}

\begin{proof}
	For part (1), it is straightforward to verify the axioms pertaining to naturality of $\chi$, normality, and compatibility with identity and constant map assigners. The internal whiskering operators for $[-, ?]$ are constructed via a restriction of $\mathbf{Gray}$-functors $\mathbf{Tricat}_\text{s}\left(\mathfrak{B}, \mathfrak{C}\right) \rightarrow \mathbf{Tricat}_\text{s}\left(\left(\mathfrak{A}, \mathfrak{B}\right),\left(\mathfrak{A}, \mathfrak{C}\right)\right)$, but they are also restrictions of the internal whiskering operators for $\mathbf{Ps}_\text{s}\left(-, ?\right)$. It follows that the $\mathbf{Gray}$-functors $\chi$ are compatible with these internal whiskering operators. Alternatively, compatibility may be directly checked for all combinations $\left(\Theta, \Psi, \Phi\right)$ of $n$-cells $\Theta \in [\mathfrak{B}, \mathfrak{C}]$, $m$-cells in $\Psi \in [\mathfrak{A}, \mathfrak{B}]$ and $l$-cells in $\Phi \in \mathfrak{A}$ satisfying $l + m + n \leq 3$. For $n=1$, it suffices to check this on morphisms in $[\mathfrak{B}, \mathcal{C}]$ which are semi-strict trinatural transformations; the required commutativity on semi-strictly decomposable trinatural transformations follows by functoriality at the level of underlying categories of either side of the diagram in axiom (2). This is similar to the proofs in Subsection \ref{Subsection the associativity axiom}, so we omit details. As in the proof of Proposition \ref{assoc condition underlying categories}, the proof of the case $l = m = n = 1$ proceeds by first assuming that $\Theta$ is semi-strict, and then using functoriality from Proposition \ref{functoriality of A, p} to extend this to semi-strictly decomposable $\Theta$.
	\\
	\\
	\noindent For part (2), the $\mathbf{Gray}$-functor is clearly the identity on objects and between hom-categories. It therefore suffices to show that it is essentially surjective on morphisms. This amounts to any unital trinatural transformation between $\mathbf{Gray}$-functors out of a cofibrant $\mathbf{Gray}$-category being isomorphic to a semi-strict trinatural transformation via a trimodification with identity $2$-cell components. But this follows from Proposition 4.3.2 of \cite{Miranda strictifying operational coherences}.
\end{proof}

\section{Centres of braided $\mathbf{Gray}$-monoids}\label{Section Centres of braided Gray monoids}

\noindent We relate $\mathbf{Gray}$-categories of the form $\mathbf{Ps}\left(\mathfrak{K}, \mathfrak{K}\right)_\text{s}$ where $\mathfrak{K}$ has one object, to generalised centres described in \cite{Crans Generalised Centers}. The reader is encouraged to refer to the explicit description of the $\mathbf{Gray}$-category given in Remark \ref{Gray Functor Gray category with weak higher cells}. We describe how generalised interchangers of Section \ref{Data of the closed structure} equip the full sub-$\mathbf{Gray}$-monoid on $1_{\mathfrak{K}}$ with a braid structure. The strictification results of \cite{Miranda strictifying operational coherences} are then used to improve upon the strictification results for braided monoidal bicategories \cite{Gurski Loop Spaces}.
\\
\\
\noindent Let $\mathfrak{K}$ be a $\mathbf{Gray}$-category with one object, and let $\left(\mathcal{K}, \bigoplus, I\right)$ denote its hom-$2$-category, equipped with the structure of a $\mathbf{Gray}$-monoid with tensor $\bigoplus$ corresponding to the $\mathbf{Gray}$-categorical composition and unit $I$ corresponding to the identity morphism. Consider the $2$-category $\mathbf{Ps}\left(\mathfrak{K}, \mathfrak{K}\right)_\text{s}\left(1_\mathfrak{K}, 1_\mathfrak{K}\right)$ of unital trinatural transformations, trimodifications and perturbations, equipped with the $\mathbf{Gray}$-monoid structure corresponding to composition of trinatural transformations. As we will explain, this $\mathbf{Gray}$-monoid is precisely the generalised centre $\mathscr{Z}\left(\mathcal{K}\right)$ described in \cite{Crans Generalised Centers}. Moreover, the generalised interchangers between trimodifications and trinatural transformations that we have analysed in Section \ref{Data of the closed structure} exactly correspond to the braid structure described in \cite{Crans Generalised Centers}.
	\\
	\\
	\noindent A unital trinatural transformation $\left(X, \beta_{-}^{X}, r_{-, ?}^{X}\right): 1_\mathfrak{K} \rightarrow 1_\mathfrak{K}$ consists of the following data.
	
	\begin{itemize}
		\item Since there is a unique object of $\mathfrak{K}$, its component consists of a $1$-cell in $\mathfrak{K}$ which is given by an object $X \in \mathcal{K}$.
		\item For every object $Y \in \mathcal{K}$, there is an adjoint equivalence in $\mathcal{K}$ with left adjoint $\beta_{Y}^{X}: X \bigoplus Y \rightarrow Y \bigoplus X$.
		\item For every morphism $g: Y \rightarrow Y'$ in $\mathcal{K}$, there is an invertible $2$-cell $\beta_{g}^{X}$ as depicted below left.
		\item Since every pair of objects in $\mathcal{K}$ are composable as morphisms in $\mathfrak{K}$, every $Y, Z \in \mathcal{K}$ has a compositor. These compositors are invertible $2$-cells as depicted below right.
		
		$$\begin{tikzcd}[font=\fontsize{9}{6}]
			X \bigoplus Y \arrow[dd, "\beta_{Y}^{X}"']\arrow[rr, "X\bigoplus g"] &{}\arrow[dd, Rightarrow, shorten = 10, "\beta_{g}^{X}"]
			&
			X\bigoplus Y'
			\arrow[dd, "\beta_{Y'}^{X}"]
			\\
			\\
			Y\bigoplus X \arrow[rr, "g\bigoplus X"']
			&{}& Y'\bigoplus X &{}
		\end{tikzcd} 	\begin{tikzcd}[ font=\fontsize{9}{6}]
			X \bigoplus Y \bigoplus Z \arrow[rr, "\beta_{Y\bigoplus Z}^{X}"]\arrow[rdd, "\beta_{Y}^{X}\bigoplus Z"'] &{}& Y \bigoplus Z \bigoplus X
			\\
			\\
			&Y \bigoplus X \bigoplus Z \arrow[ruu, "Y \bigoplus \beta_{Z}^{X}"']\arrow[uu, Rightarrow, shorten = 10, "r_{Y, Z}^{X}"]
		\end{tikzcd}$$
		
	\end{itemize}
	\noindent These data must then satisfy the axioms for a trinatural transformation. We compare these conditions to those discussed in 3.11 of \cite{Crans Generalised Centers}, where $\beta_{-}^{X}$ is denoted $R_{X, -}$ and $r_{-, ?}^{X}$ is denoted $\tilde{R}_{X|-, ?}$. Pseudonaturality of $\left(Y, g\right) \mapsto \left(\beta_{Y}^{X}, \beta_{g}^{X}\right)$ matches with the pseudonaturality condition from \cite{Crans Generalised Centers}, and the modification condition asked for in \cite{Crans Generalised Centers} is precisely what we called the left and right whiskering laws in Definition \ref{definition trinatural transformation}. Finally, condition (2.4) of \cite{Crans Generalised Centers} is what we call the associativity condition, while the condition $R_{X, I} = 1_{X}$, which in our notation is $\beta_{I}^{X} = 1_{X}$, says precisely that the trinatural transformation is unital. Finally, equations (2.10) and (2.11) of \cite{Crans Generalised Centers} are precisely the left and right unit laws for the trinatural transformation $\left(X, \beta_{-}^{X}, r_{-, ?}^{X}\right)$.
	\\
	\\
	\noindent A trimodification $\left(f, \beta_{-}^{f}\right): \left(X, \beta_{-}^{X}, r_{-, ?}^{X}\right) \Rightarrow \left(Y, \beta_{-}^{Y}, r_{-, ?}^{Y}\right)$ consists of the data of a morphism $f: X \rightarrow Y$, and for every $Z \in \mathcal{K}$ invertible $2$-cell $\beta_{Z}^{f}$ as depicted below.
	
	$$\begin{tikzcd}[ font=\fontsize{9}{6}]
		X \bigoplus Z \arrow[dd, "\beta_{Z}^{X}"']\arrow[rr, "f\bigoplus Z"] 
		&{}\arrow[dd, Leftarrow, shorten = 10, "\beta_{Z}^{f}"]
		&Y \bigoplus Z\arrow[dd, "\beta_{Z}^{Y}"]
		\\
		\\
		Z \bigoplus X \arrow[rr, "Z \bigoplus f"'] &{}& Z \bigoplus Y
	\end{tikzcd}$$
	
	\noindent These data are then subject to three conditions. What we call the local modification axiom corresponds to the condition specified in \cite{Crans Generalised Centers} that $Z \mapsto \beta_{Z}^{f}$ is a modification. What we call the unit condition for a trimodification is written as $R_{f, I} = \text{id}_{f}$ in \cite{Crans Generalised Centers}, while what we call the composition condition for a trimodification corresponds to what is referred to as `naturality of $\tilde{r}_{X | -, ?}$ in $f: X \rightarrow Y$' in \cite{Crans Generalised Centers}. Finally, a perturbation $\phi: \left(f, \beta_{-}^{f}\right) \Rrightarrow \left(g, \beta_{-}^{g}\right)$ consists of a $2$-cell $\phi: f \Rightarrow g$ in $\mathcal{K}$ subject to a certain commutativity condition which is reasonable to call `$2$-naturality of $\beta_{-}^{f}$ in $\phi$', as it is in \cite{Crans Generalised Centers}.
	\\
	\\
	\noindent The $2$-category structure on $\mathbf{Ps}\left(\mathfrak{K}, \mathfrak{K}\right)_\text{s}\left(1_\mathfrak{K}, 1_\mathfrak{K}\right)$ is described explicitly in Remark \ref{Gray Functor Gray category with weak higher cells}. From this description it is easily seen to coincide with the one described in \cite{Crans Generalised Centers}. Similarly, the $\mathbf{Gray}$-monoid structure described on the centre of $\mathcal{K}$ in \cite{Crans Generalised Centers} also corresponds to composition of trinatural transformations, whiskerings between trinatural transformations and trimodifications, and interchange between trimodifications, in the $\mathbf{Gray}$-category $\mathbf{Ps}\left(\mathfrak{K}, \mathfrak{K}\right)_\text{s}$, as described in Remark \ref{Gray Functor Gray category with weak higher cells}.
	\\
	\\
	\noindent We now analyse the braid structure. Looking at the braidings pseudonatural equivalences $\beta_{-, ?}$ and Reidemeister invertible modifications $r_{-, ?}^{!}$ and $l_{-, ?}^{!}$ described in Section 3.1.3 of \cite{Crans Generalised Centers}, we find that these structures are all captured by the internal whiskering operators $\mathbf{Ps}\left(\mathfrak{K}, -\right)_\text{s}: \mathbf{Ps}\left(\mathbf{Ps}\left(\mathfrak{K}, \mathfrak{K}\right)_\text{s}, \mathbf{Ps}\left(\mathfrak{K}, \mathfrak{K}\right)_\text{s}\right)_\text{s}$. Details follow.
	
	\begin{itemize}
		\item The braiding of objects corresponds exactly to the trimodification mediating interchange between two trinatural transformations as described in Proposition \ref{interchanging tritransformations}.
		\item The braiding of an object with a morphism is precisely the perturbation mediating interchange between a trinatural transformation with a trimodification as in Proposition \ref{interchanging tritransformation with trimodification}.
		\item The braiding of a morphism with an object is precisely the perturbation mediating interchange between a trimodification with a trinatural transformation, as described in Proposition \ref{interchanging trimodification with tritransformation}.
		\item The right Reidemeister described in \cite{Crans Generalised Centers} is precisely the compositor perturbation of Proposition \ref{compositor of whiskering with p}.
		\item The fact that the left Reidemeister described in \cite{Crans Generalised Centers} is the identity corresponds to the fact that $\mathbf{Ps}\left(\mathfrak{K}, -\right)_\text{s}$ is functorial between underlying categories. This is shown in Proposition \ref{functoriality of A, p}.
	\end{itemize}

\begin{proposition}\label{Proposition centre as 3 k transfors}
	Let $\mathcal{K}$ be a $\mathbf{Gray}$-monoid. Then the semi-strict braided monoidal $2$-category $\mathscr{Z}\left(\mathcal{K}\right)$ of \cite{Crans Generalised Centers} coincides with the analogous structure on $\mathbf{Ps}\left(\Sigma \mathcal{K}, \Sigma \mathcal{K}\right)_{s}\left(1, 1\right)$ defined via composition and interchangers of trinatural transformations, trimodifications and perturbations, as described above.
\end{proposition}

	\noindent As shown in Theorem 27 of \cite{Gurski Loop Spaces} if $\mathcal{K}$ is a braided monoidal bicategory seen as a one object tricategory, then the triequivalences $\mathbf{Gr}\left(\Sigma\mathcal{K}\right) \rightarrow \Sigma\mathcal{K}$ and $\Sigma\mathcal{K} \rightarrow \mathbf{Gr}\left(\Sigma\mathcal{K}\right)$ will moreover have structure making them compatible with the braidings. In this way, any braided monoidal bicategory $\mathcal{K}$ is braided monoidal biequivalent to a semi-strict braided monoidal $2$-category $\overline{\mathcal{K}}$ in the sense of \cite{Crans Generalised Centers}. By Theorem 18 of \cite{Baez Neuchl Braided Monoidal 2-categories}, as corrected in Section 3.2 of \cite{Crans Generalised Centers}, every semi-strict braided monoidal $2$-category is in turn braided monoidal biequivalent to one in which moreover each left Reidemeister modification $l_{-, ?}^{X}$ is an identity, via a braided monoidal biequivalence depicted as $\zeta$ below.
	
	$$\begin{tikzcd}
		\mathcal{K} \arrow[rr, "\sim"] 
		&& \overline{\mathcal{K}} \arrow[rr, "\zeta"]
		&& \mathscr{Z}\left(\overline{\mathcal{K}}\right)
	\end{tikzcd}$$

\begin{theorem}\label{coherence for braided Gray monoids}
	Every braided monoidal bicategory $\mathcal{K}$ in the sense of \cite{Gurski Loop Spaces} is braided monoidally biequivalent to a semi-strict braided monoidal $2$-category $\mathcal{K}'$ in the sense of \cite{Crans Generalised Centers}, in which moreover 
	
	\begin{itemize}
		\item every left Reidemeister modification $l^X_{-, ?}$ is the identity,
		\item every object $X$ is isomorphic to some object $X'$ for which the right Reidemeister components $r_{Y, Z}^{X'}$ are all identities.
		\item every object $X$ admits a decomposition $X = X_{1}\bigoplus ... \bigoplus X_{n}$ with each factor $X_{i}$ being such that the right Reidemeister modification's components $r^{X_{i}}_{Y, Z}$ are all identities.
	\end{itemize}
\end{theorem}

\begin{proof}
	\noindent By the analysis in Subsection 4.1 of \cite{Miranda strictifying operational coherences}, the choice of $\overline{\mathcal{K}}$ given in \cite{Gurski Loop Spaces} is such that its suspension $\overline{\mathfrak{K}}:= \Sigma \overline{\mathcal{K}}$ is cofibrant as a one object $\mathbf{Gray}$-category. We claim that the required $\mathcal{K}'$ can be taken to be the full sub-$2$-category of $\mathscr{Z}\left(\overline{\mathcal{K}}\right)$ on those objects which, when considered as trinatural transformations $1_{\overline{\mathcal{K}}} \rightarrow 1_{\overline{\mathcal{K}}}$, are semi-strictly decomposable. The braided monoidal biequivalence between the original $\mathcal{K}$ and $\mathscr{Z}\left(\overline{\mathcal{K}}\right)$ follows from the results of \cite{Baez Neuchl Braided Monoidal 2-categories}, \cite{Crans Generalised Centers} and \cite{Gurski Loop Spaces}, as described above. That the inclusion of the full-sub-$2$-category of $\mathscr{Z}\left(\overline{\mathcal{K}}\right)$ on semi-strictly generated trinatural transformations is a monoidal biequivalence follows from Theorem 5.3.4 part (1) of \cite{Miranda strictifying operational coherences}. But this $2$-functor is moreover strictly compatible with the braid structure. This corresponds to the aspect of Proposition \ref{proposition relationship to bourke lobbia} pertaining to compatibility with the internal whiskering operators, which encode the braid structure.
	\\
	\\
	\noindent By Proposition \ref{Proposition centre as 3 k transfors}, the strictification results for trinatural transformations out of cofibrant $\mathbf{Gray}$-categories in \cite{Miranda strictifying operational coherences} can be applied to prove parts (2) and (3). In particular, part (2) follows from Proposition 4.3.2 of \cite{Miranda strictifying operational coherences} since objects $X \in \mathscr{Z}\left(\overline{\mathcal{K}}\right)$ are certain trinatural transformations out of a cofibrant $\mathbf{Gray}$-category. Finally, the condition specified in part (3) corresponds to semi-strict decomposability of trinatural transformations, via Proposition \ref{Proposition centre as 3 k transfors}. By construction, the objects of $\mathcal{K}'$ are indeed semi-strictly decomposable trinatural transformations.
\end{proof}

\noindent In this Section we have applied the explicit understanding of the substitution structure of trinatural transformations and trimodifications developed in Section \ref{Data of the closed structure} to recognise these structures as braid structures on centres of $\mathbf{Gray}$-monoids. We used this relationship in Theorem \ref{coherence for braided Gray monoids} to improve Gurski and Crans' strictification results for braided monoidal bicategories. This will be further improved in the forthcoming \cite{Miranda weak interchange 4-categories}.

\section{Conclusion}

\noindent Although semi-strict trinatural transformations are not closed under composition, once they are closed under composition they form the morphisms of internal homs $[\mathfrak{A}, \mathfrak{B}]$ of $\mathbf{Gray}$-categories $\mathfrak{A}$ and $\mathfrak{B}$. These internal homs indeed form part of a closed structure on $\mathbf{Gray}$-$\mathbf{Cat}$, although this closed structure is not monoidal and it has poor properties in terms of enrichment. Understanding generalised interchangers between trinatural transformations and trimodifications has enabled a new perspective on generalised centres, which has in turn slightly improved strictification results for braided monoidal bicategories. We have also described a close relationship between the semi-strictly generated closed structure and the sharp closed structure of \cite{Bourke Lobbia Skew Approach to Enrichment for Gray-categories}, and in Remark \ref{Remark relation to Crans} we observed that there is a partial relationship between the semi-strictly generated closed structure and the symmetric monoidal structure of \cite{crans tensor of gray categories}. In a forthcoming paper \cite{Miranda weak interchange 4-categories}\footnote{Also see Chapter 3 of \cite{Miranda PhD}} we consider another enrichment context with better properties, which will be closely related to $\left(\mathbf{Gray}\text{-}\mathbf{Cat}, [-, ?]\right)$ and to the sharp closed structure of \cite{Bourke Lobbia Skew Approach to Enrichment for Gray-categories}.

\section{Appendices}\label{Appendices}

\subsection{Definitions of $\left(3, k\right)$-transfors between $\mathbf{Gray}$-functors}\label{subsection weak higher maps between Gray functors}

\noindent We recall the structure of trinatural transformations, trimodifications and perturbations between $\mathbf{Gray}$-functors.

\begin{definition}\label{definition trinatural transformation}
	Let $F, G: \mathfrak{A} \rightarrow \mathfrak{B}$ be $\mathbf{Gray}$-functors. A \emph{trinatural transformation} $p: F \Rightarrow G$ consists of \begin{itemize}
		\item For every $X \in \mathfrak{A}$, an arrow $p_{X}: FX \rightarrow GX$ in $\mathfrak{B}$,
		\item For every arrow $f: X \rightarrow Y$ in $\mathfrak{A}$, an adjoint equivalence $p_{f} \dashv p_{f}^{*}$ in $\mathfrak{B}\left(FX, GY\right)$ with left adjoint $p_{f}: p_{Y}.Ff \Rightarrow Gf.p_{X}$.
		\item For every $2$-cell $\phi: f \Rightarrow g: X \rightarrow Y$ in $\mathfrak{A}$, an invertible $3$-cell $p_\phi$ called the \emph{local pseudonaturality constraint}, as depicted below left,
		\item For every object $X \in \mathfrak{A}$, an invertible $3$-cell $p^X$ called the \emph{unitor}, as depicted below centre,
		\item For every \begin{tikzcd}
			X \arrow[r, "f"] & Y \arrow[r, "g"] & Z
		\end{tikzcd} in $\mathfrak{A}$, an invertible $3$-cell $p_{g, f}$ called the \emph{compositor}, as depicted below right.
	\end{itemize}
	
	$$\begin{tikzcd}[column sep = 12, font=\fontsize{9}{6}]
		p_Y.Ff
		\arrow[rr,"1.F\phi"]
		\arrow[dd,"p_f"']
		&
		{}
		\arrow[dd,Rightarrow,shorten=15,"p_\phi"]
		&
		p_Y.Fg
		\arrow[dd,"p_g"]
		\\
		\\
		Gf.p_X
		\arrow[rr,"G\phi.1"']
		&
		{}
		&
		Gg.p_X &{}
	\end{tikzcd}\begin{tikzcd}[column sep = 4, font=\fontsize{9}{6}]
		p_X.F1_X
		\arrow[rr,equal]
		\arrow[dd,"p_{1_X}"']
		&
		{}
		\arrow[dd,Rightarrow,shorten=15,"p^X"]
		&
		p_X.1_{FX}
		\arrow[dd,"1_{p_X}"]
		\\
		\\
		G1_X.p_X
		\arrow[rr, equal]
		&
		{}
		&
		{1_{GX}}.p_X&{}
	\end{tikzcd}\begin{tikzcd}[column sep = 1, font=\fontsize{9}{6}]
		p_Z.Fg.Ff
		\arrow[rr,"p_{gf}"]
		\arrow[rdd,"p_g.1"']
		&
		{}
		&
		Gg.Gf.p_X
		\\
		\\
		&
		Gg.p_Y.Ff
		\arrow[ruu,"1.p_f"']
		\arrow[uu,Rightarrow,shorten=10,"p_{g{,}f}"']
	\end{tikzcd}$$
	
	\noindent These data are subject to the following axioms.
	
	\begin{itemize}
		\item The assignation $\left(f, \phi\right) \mapsto \left(p_{f}, p_{\phi}\right)$ is a pseudonatural transformation. For $\left(f, \phi\right) \mapsto \left(p_{f}, p_{\phi}\right)$ this says that $p_{1_{f}} = 1_{p_{f}}$ and that the following equations hold for every $\Omega: \phi \Rrightarrow \phi'$ and every \begin{tikzcd}
			f \arrow[r, Rightarrow, "\phi"] & g \arrow[r, Rightarrow, "\psi"] & h
		\end{tikzcd} in $\mathfrak{A}$.

		$$\begin{tikzcd}[column sep = 18, row sep = 20, font=\fontsize{9}{6}]
			p_Y.Ff
			\arrow[rr,bend left = 45, "1.F\phi"]
			\arrow[dd,"p_f"']
			&
			{}
			\arrow[d,Rightarrow, shorten = 5,"p_\phi"]
			&
			p_Y.Fg
			\arrow[dd,"p_g"]
			\\
			&{}&&=
			\\
			Gf.p_X
			\arrow[rr,bend left = 45,"G\phi.1" {name = C}]
			\arrow[rr, bend right = 45, "G\phi'.1"' {name = D}]
			&
			{}
			&
			Gg.p_X
			\arrow[from =C, to =D, Rightarrow, shorten = 10, shift right = 5, "G\Omega .1"]
		\end{tikzcd}	\begin{tikzcd}[column sep = 18, row sep = 20, font=\fontsize{9}{6}]
			p_Y.Ff
			\arrow[rr,bend left = 45, "1.F\phi" {name = A}]
			\arrow[rr, bend right = 45, "1.F\phi'"' {name = B}]
			\arrow[dd,"p_f"']
			&
			&
			p_Y.Fg
			\arrow[dd,"p_g"]
			\\
			&{}
			\arrow[d,Rightarrow, shorten = 5,"p_\phi"]
			\\
			Gf.p_X\
			\arrow[rr, bend right = 45, "G\phi'.1"']
			&
			{}
			&
			Gg.p_X &{}
			\arrow[from =A, to =B, Rightarrow, shorten = 10, shift right = 5, "1.F\Omega"]
		\end{tikzcd}$$ 
		
		$$\begin{tikzcd}[column sep = 18, font=\fontsize{9}{6}]
			p_Y.Ff
			\arrow[rr,"1.F\phi"]
			\arrow[dd,"p_f"']
			&
			{}
			\arrow[dd,Rightarrow,shorten=15,"p_\phi"]
			&
			p_Y.Fg
			\arrow[dd,"p_g"] \arrow[rr, "1.F\psi"]&{}\arrow[dd, Rightarrow, shorten = 15, "p_{\psi}"]& p_{Y}.Fh \arrow[dd, "p_{h}"]
			\\
			&&&&&=
			\\
			Gf.p_X
			\arrow[rr,"G\phi.1"']
			&
			{}
			&
			Gg.p_X\arrow[rr, "G\psi.1"'] &{}&Gh.p_{X}
		\end{tikzcd}
		\begin{tikzcd}[column sep = 18, row sep = 24, font=\fontsize{9}{6}]
			p_Y.Ff
			\arrow[rr,"1.F\left(\psi\phi\right)"]
			\arrow[dd,"p_f"']
			&
			{}
			\arrow[dd,Rightarrow,shorten=15,"p_{\psi\phi}"]
			&
			p_Y.Fh
			\arrow[dd,"p_h"]
			\\
			\\
			Gf.p_X
			\arrow[rr,"G\left(\psi.\phi\right).1"']
			&
			{}
			&Gh.p_{X}
		\end{tikzcd}$$
		
		\item The assignation $\left(g, f\right) \mapsto p_{g, f}$ given in the compositor is a modification. This says that \begin{itemize}
			\item For every $\phi: f \Rightarrow f': X \rightarrow Y$  and $g: Y \rightarrow Z$ the first equation depicted below holds,
			\item For every $f: X \rightarrow Y$ and $\psi: g \Rightarrow g': Y \Rightarrow Z$ the second equation depicted below holds.
		\end{itemize} 
		We will respectively refer to these conditions as the \emph{right whiskering law} and \emph{left whiskering law} for $p$. They are specified by $2$-cells in $\mathfrak{A}$ being whiskered by $1$-cells on the left or on the right respectively.
	\end{itemize}

	$$\begin{tikzcd}[row sep = 20, column sep = 15, font=\fontsize{9}{6}]
		&Gg.p_{Y}.Ff\arrow[rd, "1.p_{f}"]
		\arrow[d, Rightarrow, shorten = 5, "p_{g{,}f}"]
		\\
		p_{Z}.Fg.Ff
		\arrow[ddd, "1.1_{Fg}F\phi" description]
		\arrow[ru, "p_{g}.1"]
		\arrow[rr, "p_{gf}"']
		&{}\arrow[ddd, Rightarrow, shorten = 18, "p_{g\phi}"]
		& Gg.Gf.p_{X}\arrow[ddd, "1.G\phi.1" description]
		\\
		&&&=
		\\
		\\
		p_{Z}.Fg.Ff' \arrow[rr, "p_{gf'}"']
		&{}
		& Gg.Gf'.p_{X}
	\end{tikzcd}\begin{tikzcd}[column sep = 15, font=\fontsize{9}{6}]
		&Gg.p_{Y}.Ff\arrow[rd, "1.p_{f}"]
		\arrow[ddd, shorten = 10, "1.1.F\phi" description]
		\\
		p_{Z}.Fg.Ff
		\arrow[dd, Rightarrow, shorten = 10, "{\left(p_{g}\right)}_{\left(F\phi\right)}", shift left = 10]
		\arrow[ddd, "1.1.F\phi"']
		\arrow[ru, "p_{g}.1"]
		&& Gg.Gf.p_{X}
		\arrow[dd, Rightarrow, shorten = 10, "1.p_{\phi}"', shift right = 10]
		\arrow[ddd, "1.G\phi.1"]
		\\
		\\
		{}&Gg.p_{Y}.Ff'
		\arrow[rd, "1.p_{f'}"]
		\arrow[d, Rightarrow, shorten = 5, "p_{g{,}f'}"]
		&{}
		\\
		p_{Z}.Fg.Ff'
		\arrow[ru, "p_{g}.1"]
		\arrow[rr, "p_{gf'}"']
		&{}
		& Gg.Gf'.p_{X}
	\end{tikzcd}$$
	\\
	$$\begin{tikzcd}[row sep = 20, column sep = 15, font=\fontsize{9}{6}]
		&Gg.p_{Y}.Ff\arrow[rd, "1.p_{f}"]
		\arrow[d, Rightarrow, shorten = 5, "p_{g{,}f}"]
		\\
		p_{Z}.Fg.Ff
		\arrow[ddd, "1.F\psi.1_{Ff}" description]
		\arrow[ru, "p_{g}.1"]
		\arrow[rr, "p_{gf}"']
		&{}\arrow[ddd, Rightarrow, shorten = 18, "p_{\psi f}"]
		& Gg.Gf.p_{X}\arrow[ddd, "G\psi.1.1" description]
		\\
		&&&=
		\\
		\\
		p_{Z}.Fg'.Ff \arrow[rr, "p_{g'f}"']
		&{}
		& Gg'.Gf.p_{X}
	\end{tikzcd}\begin{tikzcd}[column sep = 15, font=\fontsize{9}{6}]
		&Gg.p_{Y}.Ff\arrow[rd, "1.p_{f}"]
		\arrow[ddd, shorten = 10, "G\psi.1.1" description]
		\\
		p_{Z}.Fg.Ff
		\arrow[dd, Rightarrow, shorten = 10, "p_{\psi}.1", shift left = 10]
		\arrow[ddd, "1_{p_{Z}}.F\psi.1"']
		\arrow[ru, "p_{g}.1"]
		&& Gg.Gf.p_{X}
		\arrow[dd, Rightarrow, shorten = 10, "{\left(G\psi\right)}_{\left(p_{f}\right)}"', shift right = 10]
		\arrow[ddd, "G\psi.1.1"]
		\\
		\\
		{}&Gg'.p_{Y}.Ff
		\arrow[rd, "1.p_{f}"]
		\arrow[d, Rightarrow, shorten = 5, "p_{g'{,}f}"]
		&{}
		\\
		p_{Z}.Fg'.Ff
		\arrow[ru, "p_{g'}.1"]
		\arrow[rr, "p_{g'f}"']
		&{}
		& Gg'.Gf.p_{X}
	\end{tikzcd}$$
	\begin{itemize}
		\item For every composable triple \begin{tikzcd}
			W \arrow[r, "e"] & X \arrow[r, "f"] &Y \arrow[r, "g"] & Z 
		\end{tikzcd} in $\mathfrak{A}$, the following equation, called the \emph{associativity coherence} holds in the hom-$2$-category $\mathfrak{B}\left(W, Z\right)$. 
		\\
		\\
		\begin{tikzcd}[font=\fontsize{9}{6}]
			&
			Gg.p_Y.Ff.Fe
			\arrow[rr,"1.p_f.1"]
			\arrow[rrrdd,"1.p_{fe}"']
			\arrow[dd,Rightarrow,shorten=15,"p_{{g}{,}{fe}}"]
			&&
			Gg.Gf.p_X.Fe
			\arrow[rdd,"1.1.p_e"]
			\arrow[d,Rightarrow,shorten=4,"1.p_{{f}{,}{e}}"']
			\\
			&&&
			{}
			\\
			p_Z.Fg.Ff.Fe
			\arrow[ruu,"p_g.1.1"]
			\arrow[rrrr,"p_{gfe}"']
			&
			{}
			&&&
			Gg.Gf.Ge.p_W
		\end{tikzcd}=
		\\
		\\
		\\
		\begin{tikzcd}[font=\fontsize{9}{6}]
			&
			Gg.p_Y.Ff.Fe
			\arrow[rr,"1.p_f.1"]
			\arrow[d,Rightarrow,shorten=4,"p_{{g}{,}{f}}.1"]
			&&
			Gg.Gf.p_X.Fe
			\arrow[rdd,"1.1.p_e"]
			\arrow[dd,Rightarrow,shorten=15,"p_{{gf}{,}{e}}"']
			\\
			&
			{}
			\\
			p_Z.Fg.Ff.Fe
			\arrow[ruu,"p_g.1.1"]
			\arrow[rrrr,"p_{gfe}"']
			\arrow[rrruu,"p_{gf}.1"']
			&&&
			{}
			&
			Gg.Gf.Ge.p_W
		\end{tikzcd}
		\item For every $f: X \rightarrow Y$, the following pastings of $3$-cells in $\mathfrak{B}$ are both equal to the identity on $p_{f}$. These equations are respectively called the \emph{left and right unit laws}.
	\end{itemize}
	\noindent 
	\begin{tikzcd}[column sep = 15, font=\fontsize{9}{6}]
		&
		Gf.p_X.F1_X
		\arrow[rdd,bend right,"{1_{Gf}}.{p_{1_X}}"'{name=B}]
		\arrow[rdd,bend left=40,"{1_{Gf}}.{1_{p_X}}"{name=A}]
		\arrow[dd,Rightarrow,shorten=8,shift right=8,"p_{{f}{,}{1_X}}"']
		\\
		\\
		p_Y.Ff.F1_X
		\arrow[ruu,"{p_f}.{1_{F1_X}}",bend left=30]
		\arrow[rr,"p_f"']
		&
		{}
		&
		Gf.G1_X.p_X
		\arrow[from=A,to=B,Rightarrow,shorten=10,"1.{p^X}"',shift left=2]
		&
		{}
	\end{tikzcd}  \begin{tikzcd}[column sep = 15, font=\fontsize{9}{6}]
		&
		G1_Y.p_Y.Ff
		\arrow[rdd,bend left=30,"{1_{G1_Y}}.{p_f}"]
		\arrow[dd,Rightarrow,shorten=8,"p_{{1_Y}{,}{f}}",shift left=8]
		\\
		\\
		p_Y.F1_Y.Ff
		\arrow[ruu,bend left=40,"{1_{p_Y}}.{1_{Ff}}"{name=A}]
		\arrow[ruu,bend right,"{p_{1_Y}}.{1_{Ff}}"'{name=B}]
		\arrow[rr,"p_f"']
		&
		{}
		&
		G1_Y.Gf.p_X
		\arrow[from=A,to=B,Rightarrow,shorten=10,"{p^Y}.1",shift right=2]
	\end{tikzcd}
	
\end{definition}

\begin{remark}
	The unit and counit of the adjoint equivalence $p_{f} \dashv p_{f}^{*}$ will be denoted $p_{\eta}^{f}$ and $p_{\varepsilon}^f$ respectively. There are also $3$-cell components $p_{\phi}^{*}$ depicted below left which are the mates of the $3$-cell components $p_{\phi}$. Finally, there are modification conditions for $p_{f}^{\eta}$ and $p_{f}^{\varepsilon}$, with the equation for $p_{f}^{\eta}$ depicted below. These modification conditions are derivable using the triangle identities of the adjunction $p_{f} \dashv p_{f}^{*}$, and the definition of $p_{\phi}^{*}$ via mates.  
	\\

	\noindent \begin{tikzcd}[column sep = 12, font=\fontsize{9}{6}]
		p_Y.Ff
		\arrow[rr,"1.F\phi"]
		\arrow[dd, leftarrow, "p_{f}^{*}"']
		&
		{}
		\arrow[dd,Rightarrow,shorten=15,"p_{\phi}^{*}"]
		&
		p_Y.Fg
		\arrow[dd,"p_{g}^{*}", leftarrow]
		\\
		\\
		Gf.p_X
		\arrow[rr,"G\phi.1"']
		&
		{}
		&
		Gg.p_X &{}
	\end{tikzcd} \begin{tikzcd}[row sep = 20, column sep = 15, font=\fontsize{9}{6}]
		&Gf.p_{X}\arrow[rd, "p^*_f"]
		\arrow[d, Leftarrow, shorten = 5, "p^{\eta}_f"]
		\\
		p_{Y}.Ff
		\arrow[ddd, "1.F\phi" description]
		\arrow[ru, "p_{f}"]
		\arrow[rr, "1"']
		&{}
		& p_Y .Ff\arrow[ddd, "1.F\phi" description]
		\\
		&=&&=
		\\
		\\
		p_{Y}.Fg \arrow[rr, "1"']
		&{}
		& p_{Y}.Fg
	\end{tikzcd}\begin{tikzcd}[column sep = 15, font=\fontsize{9}{6}]
		&Gf.p_{X}\arrow[rd, "p_{f}^{*}"]
		\arrow[ddd, "G\phi .1" description]
		\\
		p_{Y}.Ff
		\arrow[dd, Leftarrow, shorten = 10, "p_{\phi}", shift left = 10]
		\arrow[ddd, "1.F\phi"']
		\arrow[ru, "p_{f}"]
		&& p_{Y}.Ff
		\arrow[dd, Leftarrow, shorten = 10, "p_{\phi}^{*}"', shift right = 10]
		\arrow[ddd, "1.F\phi"]
		\\
		\\
		{}&Gg.p_{X}
		\arrow[rd, "p_{g}^{*}"]
		\arrow[d, Leftarrow, shorten = 5, "p_{g}^{\eta}"]
		&{}
		\\
		p_{Y}.Fg
		\arrow[ru, "p_{g}"]
		\arrow[rr, "1"']
		&{}
		& p_{Y}.Fg
	\end{tikzcd}
\end{remark}

\begin{definition}\label{degrees of strictness trinatural transformations}
	Call a trinatural transformation $p: F \rightarrow G$
	
	\begin{itemize}
		\item A \emph{pseudo-icon equivalence} if $F$ and $G$ agree on objects and its $1$-cell components are identities.
		\item \emph{locally strict} if its component pseudonatural equivalences $\left(f, \phi\right) \mapsto \left(p_{f}, p_{\phi}\right)$ are $2$-natural isomorphisms, hence if each $p_{\phi}$, $p_{\phi}^{*}$, $p_{f}^{\eta}$ and $p_{f}^{\varepsilon}$ are identities.
		\item \emph{strict}, or a \emph{$\mathbf{Gray}$-natural transformation} if all components are identities, except perhaps $p_{X}$ for objects $X \in \mathfrak{A}$.
		\item \emph{unital} if its unitors are identities.
		\item \emph{compositional} if its compositors are identities.
		\item \emph{semi-strict} if it is unital and compositional.
	\end{itemize}
\end{definition}

\begin{definition}\label{trimodification definition}
	Let $F, G: \mathfrak{A} \rightarrow \mathfrak{B}$ be $\mathbf{Gray}$-functors and let $p, q: F \rightarrow G$ be trinatural transformations. A \emph{trimodification} $\sigma: p \Rightarrow q$ consists of
	
	\begin{itemize}
		\item For every $X \in \mathfrak{A}$, a $2$-cell $\sigma_{X}: p_{X} \Rightarrow q_{X}$
		\item For every $f: X \rightarrow Y$, an invertible $3$-cell \begin{tikzcd}
			p_{Y}.Ff
			\arrow[rr, "p_{f}"]
			\arrow[dd, "\sigma_{Y}.1_{Ff}"']
			& {}\arrow[dd, Rightarrow, shorten = 10, "\sigma_{f}"]
			&
			Gf.p_{X}
			\arrow[dd, "1_{Gf}.\sigma_{X}"]
			\\
			\\
			q_{Y}.Ff \arrow[rr, "q_{f}"']
			&{}&
			Gf.q_{X} 
		\end{tikzcd} 
	\end{itemize}
	\noindent These data are subject to the following axioms.
	
	\begin{itemize}
		\item For every $X, Y \in \mathfrak{A}$, the assignation $f \mapsto \sigma_{f}$ defines an invertible modification, from the pseudonatural transformation below left to the pseudonatural transformation below right. Here ${\left(\Phi\right)}^{*}$ and ${\left(\Phi\right)}_{*}$ denote restriction and extension along $\Phi$ respectively. 
		
		$$\begin{tikzcd}[column sep = 10, row sep = 15, font=\fontsize{9}{6}]
			&&
			\mathfrak{A}\left(X{,}Y\right)
			\arrow[rrdd,"G"]
			\arrow[lldd, "F"']
			\\
			\\
			\mathfrak{B}\left(FX{,}FY\right)
			\arrow[rrdd,bend right,"{p_{Y}}_{*}"']
			&&&&
			\mathfrak{B}\left(GX{,} GY\right)
			\arrow[lldd,"p_{X}^{*}"'name=A]
			\arrow[lldd,bend left,"q_{X}^{*}" name=B]
			\arrow[llll,shorten=25,Leftarrow,"{p}^{X{,}Y}"]
			\\
			\\
			&&
			\mathfrak{B}\left(FX{,}GY\right)
			\arrow[from=A,to=B,"\sigma_{X}^{*}",Rightarrow,shorten=8]
		\end{tikzcd}\begin{tikzcd}[column sep = 10, row sep = 15, font=\fontsize{9}{6}]
			&&
			\mathfrak{A}\left(X{,}Y\right)
			\arrow[rrdd,"G"]
			\arrow[lldd, "F"']
			\\
			\\
			\mathfrak{B}\left(FX{,}FY\right)
			\arrow[rrdd, "{q_{Y}}_{*}"name = A]
			\arrow[rrdd,bend right,"{p_{Y}}_{*}"'name=B]
			&&&&
			\mathfrak{B}\left(GX{,} GY\right)
			\arrow[lldd,bend left,"q_{X}^{*}"]
			\arrow[llll,shorten=25,Leftarrow,"{p}^{X{,}Y}"]
			\\
			\\
			&&
			\mathfrak{B}\left(FX{,}GY\right)
			\arrow[from=A,to=B,"{\sigma_{Y}}_{*}",Leftarrow,shorten=8]
		\end{tikzcd}$$
		
		\noindent  This says that for every $\phi: f \Rightarrow g: X \rightarrow Y$ in $\mathfrak{A}$, the following equation holds in $\mathfrak{B}$. We will refer to this as the \emph{local modification condition} for $\sigma$.
	\end{itemize}
	
	\noindent\begin{tikzcd}[column sep = 12, font=\fontsize{9}{6}]
		&&p_{Y}.Ff
		\arrow[lldd, "1.F\phi"']
		\arrow[rrdd, "\sigma_{Y}.1" description]
		\arrow[dddd, Rightarrow, "{\left(\sigma_{Y}\right)}_{\left(F\phi\right)}", shorten = 30, shift right = 6]
		\arrow[rr, "p_f"]
		&&Gf.p_{X}
		\arrow[dd, Rightarrow, "\sigma_{f}", shorten = 10]
		\arrow[rrdd, "1.\sigma_{X}"]
		\\
		\\
		p_{Y}.Fg\arrow[rrdd, "\sigma_{Y}.1"']
		&&&&q_{Y}.Ff
		\arrow[lldd, "1.F\phi" description]
		\arrow[rr, "q_{f}"]
		\arrow[dd, Rightarrow, "q_{\phi}", shorten = 10]
		&&Gf.q_{X}\arrow[lldd, "G\phi.1_{q_{X}}"]&=
		\\
		\\
		&&q_{Y}.Fg
		\arrow[rr, "q_{g}"']
		&&Gg.q_{X}
	\end{tikzcd}\begin{tikzcd}[column sep = 12, font=\fontsize{9}{6}]
		&&p_{Y}.Ff
		\arrow[lldd, "1.F\phi"']
		\arrow[dd, Rightarrow, "p_{\phi}"', shorten = 10]
		\arrow[rr, "p_f"]
		&&Gf.p_{X}
		\arrow[lldd, "G\phi.1" description]
		\arrow[dddd, Rightarrow, "{\left(G\phi\right)}_{\left(\sigma_{X}\right)}"', shorten = 30, shift left = 6]
		\arrow[rrdd, "1.\sigma_{X}"]
		\\
		\\
		p_{Y}.Fg
		\arrow[rrdd, "\sigma_{Y}.1"']
		\arrow[rr, "p_{g}"]
		&&Gg.p_{X}
		\arrow[dd, Rightarrow, "\sigma_{g}"', shorten = 10]
		\arrow[rrdd, "1.\sigma_{X}" description]
		&&
		&&Gg.p_{X}
		\arrow[lldd, "G\phi.1"]
		\\
		\\
		&&Gf.q_{X}
		\arrow[rr, "q_{g}"']
		&&Gg.q_{X}
	\end{tikzcd}
	\begin{itemize}
		\item For every $X \in \mathfrak{A}$, the following equation, called the \emph{unit law}, holds in $\mathfrak{B}$.

		\begin{tikzcd}[column sep = 18, row sep = 15, font=\fontsize{9}{6}]
			p_X.F1_{X}
			\arrow[rr,bend left = 45, "1_{p_X}"]
			\arrow[dd,"\sigma_{X}.1"']
			&
			=
			&
			G1_{X}.p_{X}
			\arrow[dd,"1.\sigma_{X}"]
			\\
			&{}&&=
			\\
			q_{X}.F1_{X}
			\arrow[rr,bend left = 45,"1_{q_X}" {name = C}]
			\arrow[rr, bend right = 45, "q_{1_{X}}"' {name = D}]
			&
			{}
			&
			G1_{X}.q_X
			\arrow[from =C, to =D, Rightarrow, shorten = 10, "q^{X}"]
		\end{tikzcd}	\begin{tikzcd}[column sep = 18, row sep = 15, font=\fontsize{9}{6}]
			p_X.F1_{X}
			\arrow[rr,bend left = 45, "1_{p_X}" {name = A}]
			\arrow[rr, bend right = 45, "p_{1_{X}}"' {name = B}]
			\arrow[dd,"\sigma_{X}.1_{F1_{X}}"']
			&
			&
			G1_{X}.p_{X}
			\arrow[dd,"1_{G1_{X}}.\sigma_{X}"]
			\\
			&{}
			\arrow[d,Rightarrow, shorten = 5,"\sigma_{1_{X}}"]
			\\
			q_{X}.F1_{X}\
			\arrow[rr, bend right = 45, "q_{1_{X}}"']
			&
			{}
			&
			G1_{X}.q_X &{}
			\arrow[from =A, to =B, Rightarrow, shorten = 10, "p^X"]
		\end{tikzcd} 
		
		\item For every \begin{tikzcd}
			X \arrow[r, "f"] & Y \arrow[r, "g"] & Z
		\end{tikzcd} in $\mathfrak{A}$, the following equation, called the \emph{composition law}, holds in $\mathfrak{B}$.
		
		\begin{tikzcd}[row sep = 20, column sep = 15, font=\fontsize{9}{6}]
			&Gg.p_{Y}.Ff\arrow[rd, "1.p_{f}"]
			\arrow[d, Rightarrow, shorten = 5, "p_{g{,}f}"]
			\\
			p_{Z}.Fg.Ff
			\arrow[ddd, "\sigma_{Z}.1.1" description]
			\arrow[ru, "p_{g}.1"]
			\arrow[rr, "p_{gf}"']
			&{}\arrow[ddd, Rightarrow, shorten = 18, "\sigma_{gf}"]
			& Gg.Gf.p_{X}\arrow[ddd, "1.1.\sigma_{X}" description]
			\\
			&&&=
			\\
			\\
			q_{Z}.Fg.Ff \arrow[rr, "q_{gf}"']
			&{}
			& Gg.Gf.q_{X}
		\end{tikzcd}\begin{tikzcd}[column sep = 15, font=\fontsize{9}{6}]
			&Gg.p_{Y}.Ff\arrow[rd, "1.p_{f}"]
			\arrow[ddd, shorten = 10, "1_{Gg}.\sigma_{Y}.1" description]
			\\
			p_{Z}.Fg.Ff
			\arrow[dd, Rightarrow, shorten = 10, "\sigma_{g}.1", shift left = 10]
			\arrow[ddd, "\sigma_{Z}.1.1"']
			\arrow[ru, "p_{g}.1"]
			&& Gg.Gf.p_{X}
			\arrow[dd, Rightarrow, shorten = 10, "1.\sigma_{f}"', shift right = 10]
			\arrow[ddd, "1.1.\sigma_{X}"]
			\\
			\\
			{}&Gg.q_{Y}.Ff
			\arrow[rd, "1.q_{f}"]
			\arrow[d, Rightarrow, shorten = 5, "q_{g{,}f}"]
			&{}
			\\
			q_{Z}.Fg.Ff
			\arrow[ru, "q_{g}.1"]
			\arrow[rr, "q_{gf}"']
			&{}
			& Gg.Gf.q_{X}
		\end{tikzcd}
	\end{itemize}
	
	\noindent Call a trimodification
	
	\begin{itemize}
		\item \emph{strict} if its $3$-cell components are identities.
		\item \emph{costrict} if its $2$-cell components are identities.
	\end{itemize}
\end{definition}

\begin{remark}
	The $3$-cell component $\sigma_{f}$ of a trimodification $\sigma$ has a uniquely determined mate given by a $3$-cell $\sigma_{f}^{*}$ which can be produced by pasting $\sigma_{f}$ along $p_{f}^{\varepsilon}$ and $q_{f}^{\varepsilon}$. By mateship, $\sigma_{f}^{*}$ also satisfies conditions which are dual to those specified for $\sigma_{f}$.
	
	$$\begin{tikzcd}[font=\fontsize{9}{6}]
		p_{Y}.Ff
		\arrow[rr, leftarrow, "p_{f}^{*}"]
		\arrow[dd, "\sigma_{Y}.1_{Ff}"']
		& {}\arrow[dd, Rightarrow, shorten = 10, "\sigma_{f}^{*}"]
		&
		Gf.p_{X}
		\arrow[dd, "1_{Gf}.\sigma_{X}"]
		\\
		\\
		q_{Y}.Ff \arrow[rr, leftarrow, "q_{f}^{*}"']
		&{}&
		Gf.q_{X} 
	\end{tikzcd} $$
	
\end{remark}
\begin{definition}\label{Definition perturbation}
	Let $\sigma: p \Rightarrow q$ be a trimodification as in Definition \ref{trimodification definition} and let $\tau: p \Rightarrow q$ be another trimodification. A \emph{perturbation} $\sigma \Rrightarrow \tau$ consists of the assignation to every object $X \in \mathfrak{A}$, a $3$-cell $\Omega: \sigma_{X} \Rrightarrow \tau_{X}$ in $\mathfrak{B}$ such that the following equation holds for every $f: X \rightarrow Y$ in $\mathfrak{A}$.

	$$\begin{tikzcd}[column sep = 18, row sep = 20, font=\fontsize{9}{6}]
		p_Y.Ff
		\arrow[rr,bend left = 45, "\sigma_{Y}.1"]
		\arrow[dd,"p_f"']
		&
		{}
		\arrow[d,Rightarrow, shorten = 5,"\sigma_{f}"]
		&
		q_Y.Ff
		\arrow[dd,"q_{f}"]
		\\
		&{}&&=
		\\
		Gf.p_X
		\arrow[rr,bend left = 45,"1.\sigma_{X}" {name = C}]
		\arrow[rr, bend right = 45, "1.\tau_{X}"' {name = D}]
		&
		{}
		&
		Gf.q_X
		\arrow[from =C, to =D, Rightarrow, shorten = 10, shift right = 5, "1.\Omega_{X}"]
	\end{tikzcd}	\begin{tikzcd}[column sep = 18, row sep = 20]
		p_Y.Ff
		\arrow[rr,bend left = 45, "\sigma_{Y}.1" {name = A}]
		\arrow[rr, bend right = 45, "\tau_{Y}.1"' {name = B}]
		\arrow[dd,"p_f"']
		&
		&
		q_Y.Ff
		\arrow[dd,"q_{f}"]
		\\
		&{}
		\arrow[d,Rightarrow, shorten = 5,"\tau_f"]
		\\
		Gf.p_{X}\
		\arrow[rr, bend right = 45, "1.\tau_{X}"']
		&
		{}
		&
		Gf.q_{X} &{}
		\arrow[from =A, to =B, Rightarrow, shorten = 10, shift right = 5, "\Omega_{Y}.1"]
	\end{tikzcd} $$
	
\end{definition}

\subsection{$\mathbf{Gray}$-categories of $\left(3, k\right)$-transfors}\label{Subsection Gray categories of 3 k transfors}

\begin{remark}\label{Gray Functor Gray category with weak higher cells}
	We review the structure of the $\mathbf{Gray}$-category $\mathbf{Tricat}_\text{s}\left(\mathfrak{A}, \mathfrak{B}\right)$ whose objects are $\mathbf{Gray}$-functors, morphisms are trinatural transformations, $2$-cells are trimodifications and $3$-cells are perturbations.
	
	\begin{enumerate}
		\item Given $\mathbf{Gray}$-functors $F, G: \mathfrak{A} \rightarrow \mathfrak{B}$ and trinatural transformations $p, q: F \rightarrow G$, the hom-category $\mathbf{Tricat}_{s}\left(\mathfrak{A}, \mathfrak{B}\right)\left(F, G\right)\left(p, q\right)$ has objects given by trimodifications $\sigma: p \Rightarrow q$ and morphisms given by perturbations $\Omega: \sigma \Rrightarrow \tau$. Composition in this category is inherited component-wise from the hom-categories of $\mathfrak{B}$, and the identity perturbation on $\sigma$ has $3$-cell component on $X \in \mathfrak{A}$ given by the identity $1_{\sigma_{X}}$.
		\item Given $\mathbf{Gray}$-functors $F, G: \mathfrak{A} \rightarrow \mathfrak{B}$, we describe the hom-$2$-category $\mathbf{Tricat}_{s}\left(\mathfrak{A}, \mathfrak{B}\right)\left(F, G\right)$. \begin{itemize}
			\item Its objects are trinatural transformations $p: F \rightarrow G$, with the identity trimodification on $p$ having identity $2$-cell and $3$-cell components. 
			\item Its hom-categories between trinatural transformations $p$ and $q$ are given as described in Part (1).
			\item Given trimodifications \begin{tikzcd}
				p \arrow[r, "\sigma"] &q \arrow[r, "\tau"] & r
			\end{tikzcd}, their composite is the trimodification whose $2$-cell component on $X$ is given by \begin{tikzcd}
				{p}_{X} \arrow[r, "{\sigma}_{X}"] &{q}_{X} \arrow[r, "{\tau}_{X}"] & {r}_{X}
			\end{tikzcd} 
			and whose $3$-cell component on $f: X \rightarrow Y$ is given by the following pasting in the hom-$2$-category $\mathfrak{B}\left(FX, GY\right)$.
			
			$$\begin{tikzcd}[font=\fontsize{9}{6}]
				p_{Y}.Ff
				\arrow[dd, "p_{f}"']
				\arrow[rr, "\sigma_{Y}.1"]
				& {}\arrow[dd, Leftarrow, shorten = 10, "\sigma_{f}"]
				&
				q_{Y}.Ff \arrow[dd, "q_{f}"] \arrow[rr, "\tau_{Y}.{1}"]
				&{}\arrow[dd, shorten = 10, "\tau_{f}", Leftarrow]
				&r_{Y}.Ff\arrow[dd, "r_{f}"]
				\\
				\\
				Gf.p_{X} 
				\arrow[rr, "1.\sigma_{X}"']
				&{}&
				Gf.q_{X}\arrow[rr, "1.\tau_{X}"']
				&{}&Gf.r_{X}
			\end{tikzcd}$$ 
			\item Given horizontally composable perturbations 
			\begin{tikzcd}
				p
				\arrow[rr, bend left, "\sigma" name = A]
				\arrow[rr, bend right, "\sigma '"' name = B]
				&& q
				\arrow[rr,bend left, "\tau" name = C]
				\arrow[rr, bend right, "\tau '"' name = D]
				&& r
				\arrow[from = A, to = B, Rightarrow, shorten = 5, "\Sigma"]
				\arrow[from = C, to = D, Rightarrow, shorten =5, "\Omega"]
			\end{tikzcd} their horizontal composite $\Omega * \Sigma: \tau. \sigma \Rrightarrow \tau'. \sigma'$ has component on $X$ given by the horizontal composite of their components $\Omega_{X}*\Sigma_{X}$ in the hom-$2$-category $\mathfrak{B}\left(FX, GX\right)$.
		\end{itemize} 
		\item The identity trinatural transformation on a $\mathbf{Gray}$-functor $F$ is the strict trinatural transformation whose $1$-cell components on objects are also identities.
		\item Given $\mathbf{Gray}$-functors $F, G, H: \mathfrak{A} \rightarrow \mathfrak{B}$, the composition $2$-functor
		
		$$\circ: \mathbf{Tricat}_{s}\left(\mathfrak{A}, \mathfrak{B}\right)\left(G, H\right) \otimes \mathbf{Tricat}_{s}\left(\mathfrak{A}, \mathfrak{B}\right)\left(F, G\right) \rightarrow \mathbf{Tricat}_{s}\left(\mathfrak{A}, \mathfrak{B}\right)\left(F, H\right)$$
		
		\noindent is defined in the following way. \begin{enumerate}
			\item Given a pair of composable trinatural transformations \begin{tikzcd}
				F \arrow[r, "p"] & G \arrow[r, "q"] & H
			\end{tikzcd}, their composite will have \begin{itemize}
				\item $1$-cell component on $X$ given by ${\left(q \circ p\right)}_{X}:=$ \begin{tikzcd}
					FX \arrow[r, "p_{X}"] & GX \arrow[r, "q_{X}"] & HX
				\end{tikzcd}
				\item $2$-cell component on $f: X \rightarrow Y$ given by the $2$-cell depicted below, with the rest of the adjoint equivalence constructed similarly.
				$${\left(q \circ p\right)}_{f}:= \begin{tikzcd}
					q_{Y}.p_{Y}.Ff \arrow[r, "1.p_{f}"] & q_{Y}.Gf.p_{X} \arrow[r, "q_{f}.1"] & Hf.q_{X}.p_{X}
				\end{tikzcd}$$
				\item $3$-cell component ${\left(q \circ p\right)}_{\phi}$ on $\phi: f \Rightarrow g: X \rightarrow Y$ given by the following pasting in the hom-$2$-category $\mathfrak{B}\left(FX, HY\right)$, with the $3$-cell component ${\left(q\circ p\right)}_{\phi}^{*}$ defined similarly.
				
				$$\begin{tikzcd}[column sep = 12, font=\fontsize{9}{6}]
					q_{Y}.p_{Y}.Ff
					\arrow[dd,"1.1.F\phi"']
					\arrow[rr,"1.p_f"]
					&
					{}
					\arrow[dd,Leftarrow,shorten=15,"1.p_\phi"]
					&
					q_{Y}.Gf.p_{X} \arrow[rr, "q_{f}.1"]
					\arrow[dd,"1.G\phi .1"]
					&{}\arrow[dd, Leftarrow, shorten = 15, "q_{\phi}.1"]
					& Hf.q_{X}.p_{X}\arrow[dd, "H\phi .1.1"]
					\\
					\\
					q_{Y}p_{Y}.Fg
					\arrow[rr,"1.p_{g}"']
					&
					{}
					&
					q_{Y}.Gg.p_X \arrow[rr, "q_{g}.1"']
					&{}
					& Hg.q_{X}.p_{X}
				\end{tikzcd}$$
				\item Unitor at $X$ given by the following horizontal composite in the hom-$2$-category $\mathfrak{B}\left(FX, HX\right)$. 
				
				$$\begin{tikzcd}[font=\fontsize{9}{6}]
					q_{X}p_{X}
					\arrow[rr, bend left, "1.1_{p_{X}}" name = A]
					\arrow[rr, bend right, "1.p_{1_{X}}"' name = B]
					&& q_{X}.p_{X}
					\arrow[rr,bend left, "1_{q_{X}}.1" name = C]
					\arrow[rr, bend right, "q_{1_{X}}.1"' name = D]
					&& q_{X}.p_{X}
					\arrow[from = A, to = B, Rightarrow, shorten = 5, "1.p^{X}"]
					\arrow[from = C, to = D, Rightarrow, shorten =5, "q^{X}.1"]
				\end{tikzcd}$$
				\item Compositor at \begin{tikzcd}
					X \arrow[r, "f"] & Y \arrow[r, "g"] & Z
				\end{tikzcd} given by the following pasting in the hom-$2$-category $\mathfrak{B}\left(FX, HZ\right)$.

				$$\begin{tikzcd}[font=\fontsize{9}{6}, column sep = 30]
					&q_{Z}.p_{Z}.Fg.Ff
					\arrow[d, "1.p_{g}.1"']
					\arrow[ddr, "1.p_{gf}", bend left, shift left = 3]
					\\
					&q_{Z}.Gg.p_{Y}.Ff
					\arrow[r, shorten = 15, Rightarrow, "1.p_{g{,}f}"]
					\arrow[ld, "q_{g}.1.1"']
					\arrow[rd, "1.1.p_{f}"'] &{}
					\\
					Hg.q_{Y}.p_{Y}.Ff
					\arrow[rr, Rightarrow, shorten = 40, "{\left(q_{g}\right)}_{\left(p_{f}\right)}"]
					\arrow[rd, "1.1.p_{f}"'] &{}& q_{Z}.Gg.Gf.p_{X}
					\arrow[ldd, "q_{gf}.1", bend left, shift left = 3]
					\arrow[ld, "q_{g}.1.1"']
					\\
					&Hg.q_{Y}.Gf.p_{X}
					\arrow[r, Rightarrow, shorten = 15, "q_{g{,}f}.1"]
					\arrow[d, "1.q_{f}.1"']
					&{}
					\\
					& Hg.Hf.q_{X}.p_{X}
				\end{tikzcd}$$
			\end{itemize}
			\item Consider \begin{tikzcd}
				F\arrow[rr, "p"]
				&& G
				\arrow[rr,bend left, "q" name = C]
				\arrow[rr, bend right, "q '"' name = D]
				&& H
				\arrow[from = C, to = D, Rightarrow, shorten =5, "\tau"]
			\end{tikzcd}, where $p, q$ and $q'$ are trinatural transformations and $\tau$ is a trimodification. Then the whiskering is given by the trimodification whose $2$-cell component on $X$ is given by the whiskering of $\tau_{X}$ with $p_{X}$, and whose $3$-cell component on $f: X \rightarrow Y$ is given by the pasting in the hom-$2$-category $\mathfrak{B}\left(FX, HX\right)$ depicted below left.
			\item Consider \begin{tikzcd}
				F
				\arrow[rr, bend left, "p" name = A]
				\arrow[rr, bend right, "p '"' name = B]
				&& G\arrow[rr, "q"]
				&& H
				\arrow[from = A, to = B, Rightarrow, shorten = 5, "\sigma"]
			\end{tikzcd}, where $p, p'$ and $q$ are trinatural transformations and $\sigma$ is a trimodification. Then the whiskering is given by the trimodification whose $2$-cell component on $X$ is given by the whiskering of $\sigma_{X}$ with $q_{X}$, and whose $3$-cell component on $f: X \rightarrow Y$ is given by the pasting in the hom-$2$-category $\mathfrak{B}\left(FX, HX\right)$ depicted below right.
			\item The whiskering of perturbations $\Omega$ with trinatural transformations $p$ on either side is given by the component-wise whiskering of the $3$-cell $\Omega_{X}$ in $\mathfrak{B}$ with the $1$-cell $p_{X}$ in $\mathfrak{B}$.
			\item Let \begin{tikzcd}
				F
				\arrow[rr, bend left, "p" name = A]
				\arrow[rr, bend right, "p '"' name = B]
				&& G
				\arrow[rr,bend left, "q" name = C]
				\arrow[rr, bend right, "q '"' name = D]
				&& H
				\arrow[from = A, to = B, Rightarrow, shorten = 5, "\sigma"]
				\arrow[from = C, to = D, Rightarrow, shorten =5, "\tau"]
			\end{tikzcd} be a pair of interchangeable trimodifications. Then their interchanger $\tau_{\sigma}$ is the perturbation whose component at an object $X$ is given by the interchanger in $\mathfrak{B}$ of $\tau_{X}$ and $\sigma_{X}$.
			
		\end{enumerate}
	\end{enumerate}
	
	$$\begin{tikzcd}[column sep = 8, font=\fontsize{9}{6}]
		q_{Y}.p_{Y}.Ff
		\arrow[dd,"1.\sigma_{Y}.1"description]
		\arrow[rr,"1.p_f"]
		&
		{}
		\arrow[dd,Rightarrow,shorten=15,"1.\sigma_f", shift right = 5]
		&
		q_{Y}.Gf.p_{X} \arrow[rr, "q_{f}.1"]
		\arrow[dd,"1.1.\sigma_{X}"description]
		&{}\arrow[dd, Rightarrow, shorten = 15, "{\left(q_{f}\right)}_{\left(\sigma_{X}\right)}", shift right = 5]
		& Hf.q_{X}.p_{X}\arrow[dd, "1.1.\sigma_{X}"description]
		\\
		\\
		q_{Y}{p'}_{Y}.Ff
		\arrow[rr,"1.{p'}_{f}"']
		&
		{}
		&
		q_{Y}.Gf.{p'}_{X} \arrow[rr, "q_{f}.1"']
		&{}
		& Hf.q_{X}.{p'}_{X}&{}
	\end{tikzcd}\begin{tikzcd}[column sep = 8, font=\fontsize{9}{6}]
		q_{Y}.p_{Y}.Ff
		\arrow[dd,"\tau_{Y}.1.1"description]
		\arrow[rr,"1.p_f"]
		&
		{}
		\arrow[dd,Rightarrow,shorten=15,shift right = 4, "{\left(\tau_{Y}\right)}_{\left(p_{f}\right)}"]
		&
		q_{Y}.Gf.p_{X} \arrow[rr, "q_{f}.1"]
		\arrow[dd,"\tau_{Y}.1.1"description]
		&{}\arrow[dd, Rightarrow, shorten = 15, "\tau_{f}.1"]
		& Hf.q_{X}.p_{X}\arrow[dd, "1.\tau_{X}.1"description]
		\\
		\\
		{q'}_{Y}{p}_{Y}.Ff
		\arrow[rr,"1.{p}_{f}"']
		&
		{}
		&
		{q'}_{Y}.Gf.{p}_{X} \arrow[rr, "{q'}_{f}.1"']
		&{}
		& Hf.{q'}_{X}.{p}_{X}
	\end{tikzcd}$$
\end{remark}

\begin{lemma}\label{tautological condition for failure of semi-stricts to be closed under composition}
	Suppose \begin{tikzcd}
		F \arrow[r, "p"] & G \arrow[r, "q"] & H
	\end{tikzcd} is a composable pair of semi-strict trinatural transformations between $\mathbf{Gray}$-functors from $\mathfrak{A}$ to $\mathfrak{B}$.
	
	\begin{enumerate}
		\item The composite $q \circ p$ is unital.
		\item The composite $q \circ p$ is compositional if and only if the following whiskered interchanger is the identity. Note that the unlabelled regions strictly commute by compositionality of $p$ and $q$.

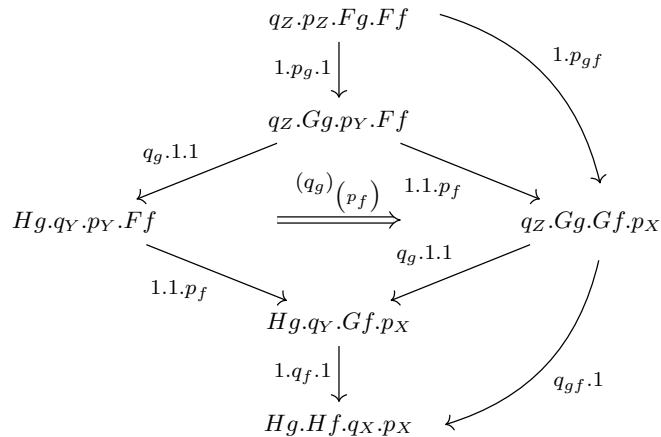
		$$\begin{tikzcd}[font=\fontsize{9}{6}, column sep = 30]
			&q_{Z}.p_{Z}.Fg.Ff
			\arrow[d, "1.p_{g}.1"']
			\arrow[ddr, "1.p_{gf}", bend left, shift left = 3]
			\\
			&q_{Z}.Gg.p_{Y}.Ff
			\arrow[ld, "q_{g}.1.1"']
			\arrow[rd, "1.1.p_{f}"'] &{}
			\\
			Hg.q_{Y}.p_{Y}.Ff
			\arrow[rr, Rightarrow, shorten = 40, "{\left(q_{g}\right)}_{\left(p_{f}\right)}"]
			\arrow[rd, "1.1.p_{f}"'] &{}& q_{Z}.Gg.Gf.p_{X}
			\arrow[ldd, "q_{gf}.1", bend left, shift left = 3]
			\arrow[ld, "q_{g}.1.1"']
			\\
			&Hg.q_{Y}.Gf.p_{X}
			\arrow[d, "1.q_{f}.1"']
			&{}
			\\
			& Hg.Hf.q_{X}.p_{X}
		\end{tikzcd}$$
		\item $q\circ p$ is not semi-strict in general.
	\end{enumerate}
\end{lemma}

\end{document}